\documentclass[11pt]{amsart}
\usepackage{latexsym,amsmath,amssymb,graphicx}

\def\urlfont{\DeclareFontFamily{OT1}{cmtt}{\hyphenchar\font='057}
              \normalfont\ttfamily \hyphenpenalty=10000}
\def\pplogo{\vbox{\kern-\headheight\kern -15pt
\halign{##&##\hfil\cr&{
\ppnumber}\cr\rule{0pt}{2.5ex}&\ppdate\cr} }} \makeatletter
\def\ps@firstpage{\ps@empty \def\@oddhead{\hss\pplogo}%
  \let\@evenhead\@oddhead 
}
\def\maketitle{\par
 \begingroup
 \def\thefootnote{\fnsymbol{footnote}}
 \def\@makefnmark{\hbox
 to 0pt{$^{\@thefnmark}$\hss}}
 \if@twocolumn
 \twocolumn[\@maketitle]
 \else \newpage
 \global\@topnum\z@ \@maketitle \fi\thispagestyle{firstpage}\@thanks
 \endgroup
 \setcounter{footnote}{0}
 \let\maketitle\relax
 \let\@maketitle\relax
 \gdef\@thanks{}\gdef\@author{}\gdef\@title{}\let\thanks\relax}
\renewcommand{\thetable}{\@Alph\c@table}
\makeatother

\DeclareFontFamily{OT1}{rsfs10}{}
\DeclareFontShape{OT1}{rsfs10}{m}{n}{ <-> rsfs10 }{}
\DeclareMathAlphabet{\mathscript}{OT1}{rsfs10}{m}{n}

 \suppressfloats
\def\ppnumber{\vbox{\baselineskip16pt }}
\def\ppdate{Revised July, 2002}
\date{} 

\title[Large $N$ dualities and transitions in geometry]{Large $N$ dualities and transitions in geometry}

\author[]{A. Grassi and  M. Rossi}

\address{Department of Mathematics, University of Pennsylvania,
Philadelphia, PA 19104} \email{grassi@math.upenn.edu}

\address{Dipartimento di Matematica, Universit\`a di Torino,
via Carlo Alberto 10, 10123 Torino} \email{michele.rossi@unito.it}

\thanks{Research partially supported by the  Institute for Advanced Study, and by National Science
Foundation grants  DMS-9706707, DMS-0074980, DMS-9729992 and by
Italian MIUR's grants. We thank the Institute for Advanced Study,
  the University of Pennsylvania and the Universit\`a di Torino
  for hospitality during  various stages of this project.}

\def\P{{\mathbb{P}}}
\def\pl{\mathbb{P}^1}

\def\p2{\mathbb{P}^2}
\def\p3{\mathbb{P}^3}
\def\p4{\mathbb{P}^4}

\def\co{\mathcal{O}}

\def\su{\operatorname{SU}}

\def\Z{\mathbb{Z}}

\def\C{\mathbb{C}}
\def\R{\mathbb{R}}
\def\M{\mathbf{M}}
\def\Q{\mathbb{Q}}
\theoremstyle{plain}
\newtheorem{theorem}{Theorem}[section]
\newtheorem{proposition}[theorem]{Proposition}
\newtheorem{corollary}[theorem]{Corollary}
\newtheorem{conjecture}[theorem]{Conjecture}
\newtheorem{remark}[theorem]{Remark}
\newtheorem{lemma}[theorem]{Lemma}
\theoremstyle{remark}
\newtheorem{example}[theorem]{Example}
\theoremstyle{definition}
\newtheorem{definition}[theorem]{Definition}
\theoremstyle{definitionproposition}
\newtheorem{definitionproposition}[theorem]{Definition/Proposition}
\newcommand{\cy}{Calabi--Yau }
\newcommand{\ka}{K\"{a}hler }

\begin{document}

\begin{abstract}
 Survey article
based on lectures given by the first author in May 2001 during 4th
SIGRAV and SAGP2001 Graduate School.

 The focus of these lectures
is Gopakumar--Vafa's insight
 that ``Large $N$ dualities'' (relating gauge theories and
 closed strings) are  realized, in certain cases, by ``transition
 in geometry".
In their pivotal 1998 example, the gauge theory is $SU(N)$
Chern-Simons theory on $S^3$, for large $N$,  and the transition
 is the  {\it ``conifold"} transition between  two Calabi--Yau
varieties. Much progress has been made to support Gopakumar and
Vafa's conjecture, including the lift of the transition to a
transformation between
  7-manifolds with $G_2$ holonomy. In another direction,
  this set up leads us to consider the uncharted territory of
   ``open Gromov-Witten invariants".
  The lectures, hence the
notes, were prepared for an audience of beginning graduate
students, in mathematics and physics, whom we hope to get
interested in this subject.
 Because most of the material
presented in these lectures comes from the
 physics literature, we aimed to build a bridge for the mathematicians
 towards the
 physics papers on the subject.\\

\end{abstract}

\maketitle

 In  1974 't Hooft
conjectured that large $N$ gauge theories are dual to closed
string theory. In 1998, Gopakumar and Vafa conjectured that
$SU(N)$ Chern--Simons theory on $S^3$  is dual to II--$A$ string
theory (with fluxes) compactified on a certain local Calabi--Yau
manifold $Y$, where the geometry of $Y$ is the key to the duality.

 It is in fact possible to do a topological surgery on $Y$
(a birational contraction followed by a complex deformation in
algebraic geometry) to obtain another Calabi--Yau $\widehat{Y}$; it
turns out that $\widehat{Y} \cong T^*S^3$.  $Y$ and $\widehat{Y}$
 are said to be related by a ``{\it geometric conifold
transition}". By previous work of Witten, Chern--Simons theory on
$S^3$ is equivalent to II--$A$ on $\widehat{Y}$, with $SU(N)$
D-branes wrapped on $S^3$.

 Evidence for the conjecture comes by comparing the
partition function for the Chern--Simons theory on $S^3$ and the
partition function for II--$A$ on $Y$.
 The corresponding mathematical quantities are certain
 topological invariants of $S^3$
and Gromov-Witten invariants on $Y$; knot invariants on $S^3$ and
 ``open Gromov-Witten invariants" on $Y$.
The ``open Gromov-Witten invariants"  should  ``count" maps of
Riemann surfaces with
 boundary to $Y$.
 We use quotation marks, as they  are not defined; yet, in this particular
 case (and, as it turns out many other cases) it is possible to make some  working assumptions
 and compute invariants.
 There is still an ambiguity, but as it turns out
 there is also an ambiguity on the Chern--Simons side, and the ambiguities
on  both sides match.

  The topic of the last lecture in Como was the
  strategy to prove the conjecture, proposed by Acharya, Atiyah, Maldacena and Vafa,
  by lifting the II--$A$ theories to $M-$theory compactified on
  $7$-dimensional manifolds with $G_2$ holonomy.
\vskip 0.1in

 Section \ref{GVC} contains the core
 of Gopakumar and Vafa's conjecture and the work of 't Hooft and Witten
leading to it; we also present  the evidence supporting the
conjecture and its mathematical implications. In Section 4, we
present the strategy of Acharya, Atiyah, Maldacena and Vafa   and include
some basics on spaces with $G_2$ holonomy.

 The first Section describes in detail the geometry of the conifold
transition between two manifolds (which are local Calabi--Yau),
because the local geometry is the key to the duality. We also include
two sections on transitions between Calabi-Yau threefolds and
their significance in algebraic geometry and the physics of string
theory. In  Section 2 we present  some background on
 Chern--Simons theory.
\vskip 0.2in
 The lectures, hence the notes, were prepared for an
audience of beginning graduate students, in mathematics and
physics, whom we hoped to get interested in this subject.
 Because most of the material presented in these lectures comes from the
 physics literature, we aimed to build a bridge for the mathematicians
 towards the
 physics papers on the subject.
 On one hand, we tried to make these lectures self-contained and did not assume much knowledge beyond the
 first/second year courses.
 On the other, we  thought it was important to outline links between these lectures
 and other research topics in string theory and mathematics,
 even when these were not essential to the main motif of the lectures.
 In these cases, we just gave statements, without necessarily defining all the terms involved.

 \vskip 0.2in

We gloss over the notion of wrapped D-branes and Lagrangian
submanifolds, as these were discussed in A. Lerda and K. Fukaya's
lectures, as well as many aspects of conformal field theory, the
topic of Y. Stanev's lectures.
 There is no discussion of II--$A$ theory itself, partly because of
 time constraints, partly because II--$A$, II--$B$ theories
 and Gromov-Witten invariants have
recently been in the spot light, thanks to the celebrated ``mirror
 symmetry".
 \vskip 0.1in
 Many of the results presented in these lectures appeared in preprint
 form, or were announced, while the lectures were prepared and
 given. Other related papers appeared afterwards; we do not
 discuss these papers, as  the notes
 closely follow the lectures.

The second author attended the lectures and at the end wrote
completely sections \ref{HAMILTON}, \ref{cli} and the Appendices,
which were only sketched in the lectures.

\vskip 0.2in

The first author would like to thank the organizers of the
 4th SIGRAV Graduate School on contemporary
relativity and gravitational physics and 2001 school on Algebraic
Geometry and Physics (SAGP2001) for the opportunity to give these
talks. We also would like to thank D. Freed, S. Katz, J. Maldacena
and N. Seiberg for kindly explaining their work. Thanks are also
due to B. Agboola, B. S. Acharya, R. Donagi, S. Garbiero,  D. Harbater, P. Horja, K. Karu, D.
Morrison, B. Ovrut, J. Talvacchia, K. Uhlenbeck, I. Zharkov and
especially L. Traynor, for many useful conversations.
 We are very grateful to A. Greenspoon for his helpful comments
 on a previous draft.
\vskip 0.2in

 A.G. is  much indebted to D. E. Diaconescu, who patiently
answered many  questions on various topics concerning these
lectures.

 \tableofcontents

\section{ Geometry and topology of transitions}

 The focus of these lectures is Gopakumar-Vafa's insight
 that ``Large $N$ dualities'' (relating gauge theories and
 closed strings) are  realized, in certain cases, by ``transition
 in geometry".
In their pivotal example \cite{gvc} the gauge theory is $SU(N)$
Chern-Simons theory on $S^3$, for large $N$,  and the transition
 is {\it ``the conifold transition"} between  two Calabi--Yau
varieties $\widehat Y \supset S^3$ and $Y$. Their conjecture is
discussed in Section  \ref{GVC}, while here we describe in detail
the geometry of the conifold transition between two varieties.
The local geometry is in fact the key to the duality.

We also include two sub-sections on transitions between Calabi-Yau
threefolds and their significance in algebraic geometry and the
physics of string theory.

 $Y$ and $\widehat Y$ are  local Calabi-Yau's, i.e.
 open neighborhoods in  Calabi-Yau manifolds.
  The Calabi-Yau condition is needed to preserve the supersymmetry of the
 physical (II--$A$) string theory:

\begin{definition}
\label{cym}
\textit{A Calabi--Yau manifold}%
\index{Calabi-Yau, manifold} is a smooth $n$--dimensional complex algebraic
manifold with trivial canonical bundle, i.e. $\Omega _{Y}^{n}\cong \mathcal{O%
}_{Y}$\textit{\ and such that}
\[
H^{j}\left( \mathcal{O}_{Y}\right) =0\ \ \forall j,\ \ \ 0<j<n.
\]
\end{definition}
 It can be verified that hypersurfaces of degree $d+1$ in $\P^{d}$ are
$(d-1)$--Calabi-Yau manifolds.
Elliptic curves and $K3$ surfaces are the $one$ and $two$-dimensional
Calabi-Yau manifolds.

This definition  of \textit{Calabi-Yau variety} is the most common in
the algebraic geometry literature: it is the natural
generalization of that of a $K3$ surface. It is worthwhile to keep
in mind that there are other, non--equivalent, definitions of a
Calabi-Yau threefold; we will discuss a definition, which is
relevant in the physics context, and its equivalence to the
following one in (\ref{cyf}), Section \ref{mt}. Note also that the
current definition of $K3$ is different from the one originally
used by Weil (see for example \cite{Barth1984}). For a nice
presentation of some of the different definitions and implications
among them, see \cite{Joyce 2000}.

In the three--dimensional case it is first possible to have \textit{%
transitions}%
\index{transition} between topologically different Calabi-Yau manifolds:

\begin{definition}
(\cite{CoxKatz1999}, \cite{Morrison 1999}) Let $Y$ be a Calabi-Yau
threefold and $\phi :Y\longrightarrow \overline{Y}$ be a
bimeromorphic contraction onto a normal variety. If there exists a
complex deformation (\textit{smoothing}) of $\overline{Y}$
 to a smooth Calabi-Yau threefold $\widehat{Y}$ then the process
from $Y$ to $\widehat{Y}$ is called \textit{\ a transition}.
\end{definition}

This concept plays an important role both in algebraic geometry
and in superstring theory as we will see later. The following
transition, the \textit{conifold transition}, is the focus of the
work of Vafa and collaborators and of these lectures; in
\ref{trans-ag} we briefly discuss other transitions of Calabi-Yau
manifolds. This example is based on  Clemens' construction
\cite{Clemens1983} and reported in \cite{Greene1995} (see also
\cite{CoxKatz1999}, example 6.2.4.1).

\begin{example}
\label{excon} (Conifold transition%
\index{conifold, transition}) Let $%
\overline{Y}\subset\P^{4}\left(x_{0}:\ldots :x_{4}\right) $  be the generic quintic threefold
containing the plane $\pi $ defined by $%
x_{3}=x_{4}=0$. It is the hypersurface defined by the equation
\[
x_{3}g\left( x_{0},\ldots ,x_{4}\right) +x_{4}h\left( x_{0},\ldots
,x_{4}\right) =0
\]
where $g,h$ are generic homogeneous polynomials of degree 4 (sections in $H^{0}\left(
 \mathcal{O}_{\P ^{4}}\left( 4\right) \right) $).  $%
\overline{Y}$ is singular precisely at the sixteen points defined by the equations:
\[
x_{3}=x_{4}=g=h=0.
\]
 We will see in \ref{conifold topology} that the topology of the
variety around each singular point is that of a real cone, hence the name
\textit{conifold}. The local equation defining each singularity is
that of a {\it node} (see also Appendix \ref{sings} and equation (\ref{node}) after the definition
\ref{a1t}):
\begin{equation}
z_{1}z_{3}+z_{2}z_{4}=0 \ \ \subset \C^4.
\end{equation}
Now consider the   threefold $Y\subset $ $\P%
^{4}\times \P^{1}$ defined by the equations:

\begin{equation}\label{bly}
\left\{
  \begin{array}{ll} y_{0}g\left( x_{0},\ldots ,x_{4}\right) +y_{1}h\left( x_{0},\ldots ,x_{4}\right) =0\\
 y_{0}x_{4}-y_{1}x_{3}=0,
  \end{array} \right.
\end{equation}
\noindent with $[y_0, y_1] \in \P^1$. It can be directly verified
that $Y$ is smooth (or use Bertini's theorem); then $ \phi
:Y\longrightarrow \overline{Y}
$  is an isomorphism outside the sixteen nodes of $%
\overline{Y}$ and their inverse images in $Y$, which are sixteen
copies of $\P^{1}$s. $Y$ is a birational resolution of
$\overline{Y}$ (see Appendix \ref{sings}); $\phi$ is also called a
``small blow up" of $Y$, because the inverse images of points are
complex curves and not complex surfaces.
  In
particular $K_{Y}\sim \phi ^{*}(K_{\overline{Y}})\sim \mathcal{O}_{Y}$, that is, $%
\phi $ is a crepant resolution (see \ref{sings}). Moreover
\[
h^{1,0}\left( Y\right) =h^{2,0}\left( Y\right) =h^{1,0}\left( \overline{Y}%
\right) =h^{2,0}\left( \overline{Y}\right) =0;
\]
then $Y$ is a Calabi-Yau threefold with
\[
h^{1,1}\left( Y\right) =h^{1,1}\left( \overline{Y}\right) +1=2.
\]
Note also that all the contracted $\P^{1}$'s are on the same
extremal ray of the Mori cone $\overline{NE}\left( Y\right) ,$
(see \ref{moricone}) i.e. $\phi $ cannot be factored in other
contractions. $\phi $ is called a primitive contraction of type
$I$ (see \ref{trans-ag}). \noindent On the other hand
$\overline{Y}\subset \P^{4}$ can be deformed to the generic
quintic threefold $\widehat{Y}\subset \P^{4}$
which is again a Calabi-Yau. The process of going from $Y$ to $%
\widehat{Y}$ is a (primitive) extremal transition of type $I$. We
will see in \ref{conifold topology} that the topology of these
singularities is that of a node: this transition is often called
the \textit{conifold transition}.

By Clemens' topological analysis one can see that $Y$ and
$\widehat{Y}$ do not have the same topology. See subsection
\ref{conifold topology} and theorem \ref{Clemens thm} for more
details.
\end{example}

\vskip 0.2in \subsection{The local topology of a conifold transition\label%
{conifold topology}}~ \vskip 0.1in

 Here we analyze the local geometry
and topology of a conifold transition $Y$ to $\widehat{Y}$
presented in example \ref{excon}.

\begin{definition}\label{a1t} A threefold singularity defined by the equation
$$x^{2}+y^{2}+z^{2}+v^{2}=0$$
is called a node (nodal singularity).  (See Appendix \ref{sings}.)
\end{definition}

By a change of coordinates, the equation of the node can be
rewritten as:
\begin{equation}
z_{1}z_{3}+z_{2}z_{4}=0,  \label{node}
\end{equation}
via  the affine transformation
\begin{equation}
\begin{array}{l}
x=z_{1}+iz_{3} \\
y=z_{3}+iz_{1} \\
z=z_{2}+iz_{4} \\
v=z_{4}+iz_{2}.
\end{array}
\label{trans}
\end{equation}

The singularities of example \ref{excon} are nodes.

\begin{example} \textbf{The conifold, revisited.}

The original threefold $\overline{Y}\subset \P^{4}$ is given by
the equation:
\[
x_{3}g\left( x_{0},\ldots ,x_{4}\right) +x_{4}h\left( x_{0},\ldots
,x_{4}\right) =0
\]
By a linear projective transformation we may assume the point
$P=\left( 1:0:\ldots :0\right) $ to be one of the sixteen singular
points of $\overline{Y}$ and localize our analysis in a
neighborhood $\overline{U}$ of $P$. By intersecting $\overline{Y}$
with the affine open subset of $\P^{4}$ defined by $x_{0}\neq 0$
we get the local equation of $\overline{U}\subset
\C^{4}$%
\[
z_{3}\widetilde{g}\left( z_{1},\ldots ,z_{4}\right) +z_{4}\widetilde{h}%
\left( z_{1},\ldots ,z_{4}\right) =0
\]
where $z_{i}:=x_{i}/x_{0}$ for $i=1,\ldots 4$,
$\widetilde{g}:=g/x_{0}^{4}$ and $\widetilde{h}:=h/x_{0}^{4}$.
Since $g$ and $h$ are generic we may
assume $\widetilde{g}$ and $\widetilde{h}$ to be smooth maps $\C%
^{4}\longrightarrow \C$ submersive at the origin (i.e. at $P\in
\overline{U}$) and by the inverse function theorem we have locally
\begin{eqnarray*}
\widetilde{g}\left( z_{1},\ldots ,z_{4}\right) &=&z_{1} \\
\widetilde{h}\left( z_{1},\ldots ,z_{4}\right) &=&z_{2}
\end{eqnarray*}
up to a suitable analytic change of coordinates (this is the well known \textit{%
local submersion theorem}).
\end{example}

\begin{theorem}
\label{Clemens thm}(\cite{Clemens1983}, lemma 1.11)
\begin{enumerate}
\item \textit{Let} $\overline{U}$ \textit{be the neighborhood of a threefold nodal singularity,\\ then}
$\overline{U}$ \textit{ is a real cone over } $S^2 \times S^3$.
\item  \textit{Let } ${U}$ \textit{be a neighborhood of the strict transform
of a node in } $Y$\textit{, then\\} $U\cong  D^ 4 \times S^2
\subset  \mathbb C^{2}\times S^{2}.$  \\
{\it Furthermore} $
\mathcal{N}_{U|\P^{1}}\cong \mathcal{O}_{\P^{1}}\left( -1\right)
\oplus \mathcal{O}_{\P^{1}}\left( -1\right) $.
\item  \textit{Let } $%
\widehat{U}$ \textit{be the deformed neighborhood of a node,
then}\newline $\widehat{U}\cong  D^3 \times S^3 \subset
T^{*}S^{3}\cong \mathbb R^{3}\times S^{3}$. \textit{In particular
the non--trivial } $S^{3}$\textit{ is the vanishing cycle of }
$\widehat{U}$ \textit{and it is locally embedded as a Lagrangian
submanifold in} $T^{*}S^{3}.$
\item \textit{The
conifold transition is a local surgery which replaces a tubular
neighborhood}
$D^{4}\times S^{2}$ \textit{of the exceptional fiber} $\P_{\C%
}^{1}\cong S^{2}$ \textit{in} $U$ \textit{by} $S^{3}\times D^{3}$
\textit{to obtain a smoothing} $\widehat{U}$\ of $\overline{U}$.
{\it In particular} $U$ {\it and} $\widehat{U}$ \textit{ are
topologically different.}
 \textit{ This is the classical surgery
between two manifolds with the same boundary.}
\item  \textit{More generally, there are relations between the
Betti numbers of the Calabi--Yau manifolds} $Y$ {\it and }
$\widehat{Y}$
 \textit{ as in example \ref{excon}.}
\end{enumerate}
\begin{figure}[h]
     \begin{center}
     \scalebox{1}{\includegraphics{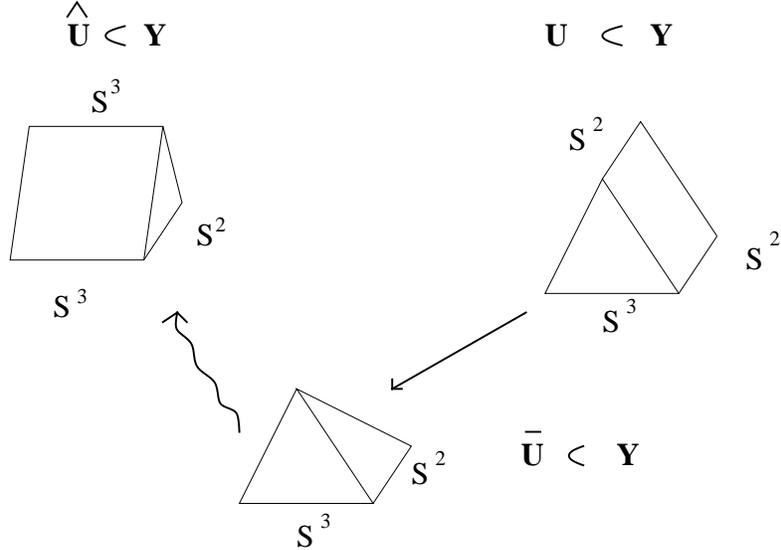}}
     \caption{The topology of the conifold transition}\label{toptrans}
     \end{center}
     \end{figure}
\end{theorem}

The invariants discussed in the rest of the paper are determined
by the local geometry around the singular locus, so we identify
(sometimes perhaps too freely) the Calabi-Yau manifolds $\widehat{Y}$
 and $Y$ with the affine varieties $\R^3 \times S^3$ and
$\R^4 \times S^2$ containing the local neighborhoods $\widehat{U}$
and $U$.

\noindent The following proof of the theorem is a review of what is
explained in the first section of \cite{Clemens1983} and also
\cite{CdO}.

\vskip 0.2in

\noindent {\it Proof:}

\noindent \textbf{(i)} As we have seen in (\ref{node}),
the local equation of a threefold $\overline{U}$ with a nodal singularity at the origin is:
\begin{equation}
z_{1}z_{3}+z_{2}z_{4}=0,
\end{equation}
Consider now the affine transformation
\begin{equation}
\begin{array}{l}
w_{1}=\left( z_{1}+z_{3}\right) /2 \\
w_{2}=i\left( -z_{1}+z_{3}\right) /2 \\
w_{3}=\left( z_{2}+z_{4}\right) /2 \\
w_{4}=i\left( -z_{2}+z_{4}\right) /2
\end{array}
\label{affinity}
\end{equation}
and    set $%
w_{j}=u_{j}+iv_{j}$; we can now identify $\overline{U}$ with the subset $%
\overline{V}\subset \R^{8}$ defined by the equations:
\begin{equation}
\begin{array}{l}
\sum_{j=1}^{4}u_{j}^{2}-\sum_{j=1}^{4}v_{j}^{2}=0 \\
\sum_{j=1}^{4}u_{j}v_{j}=0.
\end{array}
\label{real-node}
\end{equation}
Note now that there is a diffeomorphism
\[
\overline{V}\setminus \left\{ \left( 0,\ldots ,0\right) \right\}
\cong \left( \R^{4}\setminus \left\{ \left( 0,\ldots ,0\right)
\right\} \right) \times S^{2}
\]
where $S^{2}$ is the unit sphere in $\R^{3}$. In fact for every
positive real number $\rho $ we can consider the radius $\rho $
hypersphere $S_{\rho }^{7}\subset $ $\R^{8}$ and the section
$V_{\rho }:=S_{\rho }\cap \left( \overline{V}\setminus \left\{
\left( 0,\ldots ,0\right) \right\} \right) $. Clearly we get
\[
\overline{V}\setminus \left\{ \left( 0,\ldots ,0\right) \right\}
=\coprod_{\rho \in \R_{>0}}V_{\rho }
\]
On the other hand $V_{\rho }$ has equations
\[
\begin{array}{l}
\sum_{j=1}^{4}u_{j}^{2}=\sum_{j=1}^{4}v_{j}^{2}=\frac{\rho ^{2}}{2} \\
\sum_{j=1}^{4}u_{j}v_{j}=0
\end{array}
\]
Hence $V_{\rho }\cong S^{3}\times S^{2}$ since the fiber over a fixed point $%
\left( u_{1}^{o},\ldots ,u_{4}^{o}\right) \in S_{\rho
/\sqrt{2}}^{3}$ is given by the subset of $\R^{4}\left(
v_{1},\ldots ,v_{4}\right) $ defined by
\[
\begin{array}{l}
\sum_{j=1}^{4}v_{j}^{2}=\frac{\rho ^{2}}{2} \\
\sum_{j=1}^{4}u_{j}^{o}v_{j}=0
\end{array}
\]
which is clearly a $S^{2}$. Therefore
\begin{equation}
\coprod_{\rho \in \R_{>0}}V_{\rho }\cong \left( \R_{>0}\times
S^{3}\right) \times S^{2}\cong \left( \R^{4}\setminus \left\{
\left( 0,\ldots ,0\right) \right\} \right) \times S^{2}
\label{p-blow}
\end{equation}
and $\overline{U}\cong \overline{V}$ \textit{identifies with the
real cone over} $S^{3}\times S^{2}$.

\vskip 0.1in

\noindent $\bullet $ {\textbf{(ii) The blown up conifold, the
small resolution of a  nodal singularity.} }

Motivated by formula (\ref{bly}), we consider
 the standard projection\\
  $\phi: \C^4 \times \P^1 \to \C^4$
 and its restriction to the open smooth threefold $U \subset \C^4 \times \P^1$
 defined by:
\begin{equation}
\begin{array}{l}
y_{0}z_{4}-y_{1}z_{3}=0 \\
y_{0}z_{1}+y_{1}z_{2}=0,
\label{blow-node}
\end{array}
\end{equation}
with $[y_0, y_1] \in \P^1$.
 $\phi_{|U}=\varphi:  U\longrightarrow
\overline{U}$.
Recall that
 $\overline{U} $  is defined by the equation $z_1 z_3 + z_2 z _4=0$ and
 has a  nodal threefold singularity at the origin.
$\varphi $ induces an isomorphism between the open sets $U\setminus \phi ^{-1}\left(
P\right) \cong \overline{U}\setminus \left\{ \left( 0,\ldots
,0\right) \right\} \cong \overline{V}\setminus \left\{ \left(
0,\ldots ,0\right) \right\} $.
As in the previous, compact example, $U \to \overline{U}$ is
a birational resolution of $\overline{U}$ (see Appendix \ref{sings}).

\vskip 0.1in

This ``small resolution" of $\overline{U}$ was obtained by ``blowing up" the
plane $z_3=z_4=0$; by blowing up the plane $z_3=z_2=0$ we would
have  another small resolution $U_+$ isomorphic to $U$
outside the locus of the exceptional curves.  $U_+$ is called the
{\it flop} of $U$ and the birational transformation
\begin{equation}\label{cyflop}
 U \leftarrow \cdots \rightarrow U_+
\end{equation}
the ``flop".
 By analogy the transformation in
Section \ref{mt} will also be called a flop.

In particular we then have a diffeomorphism
\begin{equation}
U\setminus \phi ^{-1}\left( P\right) \cong \left( \R^{4}\setminus
\left\{ \left( 0,\ldots ,0\right) \right\} \right) \times S^{2}
\label{punctured blow}
\end{equation}
and we want to extend it to the exceptional fiber $\phi
^{-1}\left( P\right) \cong \P^{1}\cong S^{2}$ to give a
diffeomorphism
\begin{equation}
U\cong \R^{4}\times S^{2}  \label{local-blow}
\end{equation}
In order to construct it observe that under the affine transformation $%
\left( \text{\ref{affinity}}\right) $ and the above identification
$\C ^{4}\left( w_{1},\ldots ,w_{4}\right) \cong \R^{8}\left(
u_{1},\ldots ,u_{4},v_{1},\ldots ,v_{4}\right) $ the neighborhood
$U$ is sent diffeomorphically onto the subset of $\R^{8}\times
{\P_\C }^{1} $ defined by
\begin{equation}
\begin{array}{l}
y_{0}u_{3}+y_{0}v_{4}-y_{1}u_{1}-y_{1}v_{2}+i\left(
y_{0}v_{3}-y_{0}u_{4}-y_{1}v_{1}+y_{1}u_{2}\right) =0 \\
y_{0}u_{1}-y_{0}v_{2}+y_{1}u_{3}-y_{1}v_{4}+i\left(
y_{0}v_{1}+y_{0}u_{2}+y_{1}v_{3}+y_{1}u_{4}\right) =0
\end{array}
\label{eqs}
\end{equation}
Hence the fiber over a fixed point $\left( y_{0}:y_{1}\right) ^{o}\in \P%
_{\C}^{1}$ is a $\R^{4}\subset \R^{8}$ ensuring the existence of
the diffeomorphism $\left( \text{\ref{local-blow}}\right) $ up to
possibly shrinking $U$. Moreover by splitting $y_{0}$ and $y_{1}$
into real and imaginary parts the equations $\left(
\text{\ref{eqs}}\right) $ reduce to the following matricial form:
\[
\mathbf{v}=A\mathbf{u}
\]
where $\mathbf{u}$ and $\mathbf{v}$ are vectors whose entries are given by $%
u_{j}$ and $v_{j}$ respectively and $A$ is an antisymmetric matrix
uniquely determined by the fixed projective point $\left(
y_{0}:y_{1}\right) ^{o}$. Since outside of the origin the
coordinates $u_{j}$ and $v_{j}$ have to satisfy the equations
$\left( \text{\ref{real-node}}\right) $ this suffices to show that
the restriction of the diffeomorphism $\left( \text{\ref
{local-blow}}\right) $ to $U\setminus \phi ^{-1}\left( P\right) $
gives precisely the diffeomorphism $\left( \text{\ref{punctured
blow}}\right) $.

\noindent Note that $U$ can be identified with the total space of
the normal bundle $\mathcal{N}_{U|\P^{1}}$, which is a holomorphic
vector bundle of rank 2 over $\P^{1}$. By the Grothendieck theorem
(see for instance \cite{Ok}) we have the splitting

\[
\mathcal{N}_{U|\P^{1}}\cong \mathcal{O}_{\P^{1}}\left(
d_{1}\right) \oplus \mathcal{O}_{\P^{1}}\left( d_{2}\right)
\]
for some $d_{1},d_{2}\in \Z$. The local equations $\left(
\text{\ref{blow-node}}\right) $ allows us to determine those
integers. In fact we can choose two local charts on $S^{2}\cong
\P^{1}\left( y_{0}:y_{1}\right) $ around the north and the south
poles respectively. Say $\tau :=y_{0}/y_{1}$ and $\sigma
:=y_{1}/y_{0}$ are the two local coordinates on $\P^{1}$. Lifting
these charts to $\mathcal{N}_{U|\P^{1}}$ we can choose the two
local parameterizations
\[
\left( \tau ,z_{1}\right) \oplus \left( \tau ,z_{4}\right) \quad
,\quad \left( \sigma ,-z_{2}\right) \oplus \left( \sigma
,z_{3}\right) .
\]
Look at the fibre over a fixed point $\left( y_{0}:y_{1}\right)
=\left( \tau :1\right) =\left( 1:\sigma \right) $ in the gluing of
the charts. Since here
$\sigma =\tau ^{-1}$ by the local equations $\left( \text{\ref{blow-node}}%
\right) $ we get
\begin{eqnarray*}
-z_{2} &=&\sigma ^{-1}z_{1}=\tau z_{1} \\
z_{3} &=&\sigma ^{-1}z_{4}=\tau z_{4}
\end{eqnarray*}
which means that the transition functions $\tau ^{-d_{1}},\tau
^{-d_{2}}\in
\C^{*}=GL\left( 1,\C\right) $ are given by $\tau $, i.e. $%
d_{1}=d_{2}=-1$.

\vskip 0.1in \noindent  \textbf{(iii) The deformed conifold as a
symplectic manifold.}

Consider the (real) 1--parameter family of local smoothings
$\widehat{U}_{t}$ of $\overline{U}$ defined by
\begin{equation}
\begin{array}{l}
\sum_{j=1}^{4}u_{j}^{2}-\sum_{j=1}^{4}v_{j}^{2}=t \\
\sum_{j=1}^{4}u_{j}v_{j}=0
\end{array}
\ ,\ t\in \R_{>0}  \label{smoothing}
\end{equation}
Note that the generic quintic hypersurface $\widehat{Y}\subset
\P^{4}$ smoothing $\overline{Y}$ in the example \ref{excon} can be
chosen to admit local equations as in $\left(
\text{\ref{smoothing}}\right) $ for some real $t_{o}>0$ since the
real 1--dimensional arc parametrized by $t$ can be chosen
transversely with respect to the Zariski closed subset of singular
quintic hypersurfaces and connecting $\overline{Y}$ to
$\widehat{Y}$. Consider now the map
\[
\R^{8}\left( u_{1},\ldots ,u_{4},v_{1},\ldots ,v_{4}\right)
\longrightarrow \R^{8}\left( q_{1},\ldots ,q_{4},p_{1},\ldots
,p_{4}\right)
\]
defined by setting
\begin{equation}
\begin{array}{l}
q_{j}=\frac{u_{j}}{\sqrt{t+\sum_{i}v_{i}^{2}}} \\
p_{j}=v_{j}
\end{array}
\label{local-smooth}
\end{equation}
For every $t>0$ it maps $\widehat{U}_{t}$ diffeomorphically onto the\textit{%
\ cotangent bundle }$T^{*}S^{3}\cong S^{3}\times \R^{3}$\textit{\
to the unit sphere }$S^{3}\subset \R^{4}\left( q_{1},\ldots
,q_{4}\right) $ embedded in $\R^{8}$ as follows:
\begin{equation}
\begin{array}{l}
\sum_{j=1}^{4}q_{j}^{2}=1 \\
\sum_{j=1}^{4}q_{j}p_{j}=0
\end{array}
\label{ctg-bundle}
\end{equation}
Note that the 3--cycle $S_{t}\subset \widehat{U}_{t}$ described in $\R%
^{8}$ by
\[
\begin{array}{l}
\sum_{j=1}^{4}u_{j}^{2}=t \\
v_{1}=\ldots =v_{4}=0
\end{array}
\]
which vanishes when $t=0$, is diffeomorphically sent onto the
unit sphere $S^{3}\subset T^{*}S^{3}$.

\vskip 0.1in

The canonical\textit{\ }symplectic form \index{form, canonical
symplectic} given by
\[
\omega :=d\vartheta
\]
where $\vartheta :=\sum_{j=1}^{4}p_{j}dq_{j}$ is the
\textit{Liouville form of }$\R^{8}$, \index{form, Liouville}
induces a vanishing symplectic form on $S^{3}$ since this sphere
is described in $T^{*}S^{3}$ by $p_{1}=\ldots =p_{4}=0$ (locally
only three of these equations are needed). This shows that $S^{3}$\textit{\ }%
is a \textit{Lagrangian subvariety} \index{Lagrangian, subvariety}
\textit{\ of }$T^{*}S^{3}$:

\begin{definition}
A subvariety $Y\subset X$  is called \textit{Lagrangian} if $\dim
Y=\left( 1/2\right) \dim X$ and the symplectic form $\omega $ of $X$ vanishes
on every tangent vector to $Y$ i.e.
\[
\forall p\in Y,\ \forall u,v\in T_{p}Y\quad \omega \left(
u,v\right) =0
\]
\end{definition}

\noindent The same is then true for the vanishing cycle
$S_{t}\subset \widehat{U}_{t}$.

\vskip 0.3in

\noindent $\bullet$  {\textbf{ (iv) The local description of the
conifold transition.} Consider the diffeomorphism:
\begin{equation}\label{alfa}
\alpha :( \R^{4} ( \mathbf{u} ) \setminus \mathbf{0} ) \times \mathbb %
R^{4} ( \mathbf{v} ) \longrightarrow ( \mathbb R ^{4} ( \mathbf{q} )
\setminus {\ \mathbf{0} } ) \times \mathbb R ^{4} ( \mathbf{p} )
\end{equation}
given by
\[
\begin{array}{l}
q_{j}=\frac{u_{j}}{\sqrt{\sum_{i}u_{i}^{2}}} \\
p_{j}=v_{j}\sqrt{\sum_{i}u_{i}^{2}}.
\end{array}
\]
Note that, by $\left( \text{\ref{real-node}}\right) $ and $\left( \text{%
\ref{ctg-bundle}}\right) $, $\alpha$  restricts to a
diffeomorphism
\begin{equation}
U\setminus \phi ^{-1}\left( P\right) \cong \left( \R^{4}\setminus
\left\{ \mathbf{0}\right\} \right) \times S^{2}\stackrel{\alpha }{\cong }%
S^{3}\times \left( \R^{3}\setminus \left\{ \mathbf{0}\right\}
\right)
\end{equation}
In particular the fiber over a fixed point $\mathbf{u}^{o}\in \R%
^{4}\setminus \left\{ \mathbf{0}\right\} $ such that
$\sum_{i}\left( u_{i}^{o}\right) ^{2}=\rho ^{2}$, which is the
2--sphere $S_{\rho
}^{2}\subset $ $\R^{4}\left( \mathbf{v}\right) $ given by $%
\sum_{j}v_{j}^{2}-\rho ^{2}=\sum_{j=1}^{4}u_{j}^{o}v_{j}=0$, is
diffeomorphically sent onto the fiber over the fixed point $\mathbf{q}%
^{o}=\alpha \left( \mathbf{u}^{o}\right) $, which is the 2--sphere
$S_{\rho
^{2}}^{2}\subset $ $\R^{4}\left( \mathbf{p}\right) $ given by $%
\sum_{j}p_{j}^{2}-\rho ^{4}=\sum_{j=1}^{4}q_{j}^{o}p_{j}=0$.
Calling $D^{n}$ the closed unit ball in $\R^{n}$, this means
that $\alpha $\textit{\ restricts to give a diffeomorphism }
\begin{equation}
\left( D^{4}\setminus \left\{ \mathbf{0}\right\} \right) \times S^{2}%
\stackrel{\alpha }{\cong }S^{3}\times \left( D^{3}\setminus \left\{ \mathbf{0%
}\right\} \right)
\end{equation}
\textit{which reduces to the identity on their boundaries
}$S^{3}\times S^{2} $.\textit{\ } Hence recalling $\left(
\text{\ref{local-blow}}\right) $ we can cut out the interior of a
$D^{4}\times S^{2}$ around the exceptional fibre $\phi ^{-1}\left(
P\right) $ in $U$ and paste in by $\alpha $ the interior of a
$S^{3}\times D^{3}$ to get $\widehat{U}_{t}$ for some $t>0$.

\vskip 0.2in

\vskip 0.1in

\noindent $\bullet $ {\textbf{(v) The Betti numbers.} } \vskip
0.1in If $\overline{Y}$ has $N$ nodes (and no other singular
points) and $\delta $
is the number of linearly independent vanishing cycles in the smoothing $%
\widehat{Y}$, we get the following relationship between the Betti
and the Euler numbers of $Y$ and $\widehat{Y}$:
\begin{equation}
\begin{array}{l}
b^{3}\left( Y\right) =b^{3}\left( \widehat{Y}\right) -2\delta \\
b^{2}\left( Y\right) +b^{4}\left( Y\right) =b^{2}\left(
\widehat{Y}\right)
+b^{4}\left( \widehat{Y}\right) +2\left( N-\delta \right) \\
\chi \left( Y\right) =\chi \left( \widehat{Y}\right) +2N
\end{array}
\label{topology}
\end{equation}
(see \cite{Clemens1983} and \cite{WvG} for detailed proofs). Note
that by the Calabi-Yau condition the first equation above gives
the following relationship between the Hodge numbers of $Y$ and
$\widehat{Y}$:
\[
h^{2,1}\left( Y\right) =h^{1,2}\left( Y\right) =h^{2,1}\left( \widehat{Y}%
\right) -\delta =h^{1,2}\left( \widehat{Y}\right) -\delta. \ \
\diamondsuit
\]

\vskip 0.2in
 The invariants discussed in the rest of the paper are determined
 by the local geometry around the conifold locus, so we
 identify the local Calabi-Yau's  $Y, \ \widehat{Y}$ and
 $\overline{Y}$ with the local neighborhoods $U, \ \widehat{U}$ and $\overline
 U$.

\phantom{xx}

 \subsection{\label{trans-ag} Transitions of Calabi-Yau threefolds}~ \vskip 0.1in

\noindent Let $Y$ and $\overline{Y}$ be projective Calabi-Yau manifolds and $%
\phi $ a birational contraction. See  Appendix \ref{sings}
for the definitions of the different types of singularities used in
this section.

\begin{definition}
$\phi :Y\to \overline{Y}$ is a \textit{primitive }contraction
\index{contraction, primitive} if it cannot be further factored
into birational morphisms of normal varieties.
\end{definition}

A non--primitive Calabi-Yau contraction may be factored into a
composite of primitive contractions (see \cite{Wilson1989}), so, without
loss of generality we can consider $\phi $ to be primitive. In
this case the pull--back $\phi
^{*}H $ of an ample divisor $H$ on $%
\overline{Y}$ will cut the Mori cone (see \ref{moricone}) $\overline{NE}%
\left( Y\right) $ along an \textit{extremal face}. Such
contractions are also called \textit{extremal} \index{contraction,
extremal} and the associated transitions \textit{\ primitive
extremal transitions.}
\index{transition, primitive}%
\index{transition, extremal}

\begin{definition}
\cite{Wilson1992}  A primitive contraction is:
\begin{itemize}
\item  \textit{of type $I$}%
\index{contraction, of type $I$} if the exceptional locus $E$ of the
associated primitive contraction $\phi $ is composed of finitely
many curves,

\item  \textit{of type $II$}%
\index{contraction, of type $II$} if $\phi $ contracts a divisor
down to a point,

\item  \textit{of type $III$}%
\index{contraction, of type $III$} if $\phi $ contracts a divisor
down to a curve.
\end{itemize}
\end{definition}

In the first case $\phi \left( E\right) $ is composed of a finite
number of isolated singularities, each with a small resolution.
Since $Y$ is smooth these singularities are necessarily terminal
of index 1 and therefore cDV points. In the second case $E$ must
be irreducible and more precisely it is a del Pezzo surface (see
\cite{Reid1980}); $\phi \left( E\right) $ is a canonical singular point
of index 1\textit{.}

\noindent In the third case $E$ is again an irreducible surface
contracted
down to a curve $\phi \left( E\right) $ of canonical singularities for $%
\overline{Y}$. In particular if $\phi $ is crepant then $E$ is a
conic bundle over the curve $\phi \left( E\right) $ which is a
smooth curve of (generically $cA_{1}$ or $cA_2$) cDV points (see
\cite{Reid1980} and \cite{Wilson1992}, theorem 2.2).

The simplest example of a non--trivial transition of type $I$ is the
conifold transition\textit{\ } of example \ref{excon}, i.e. a
transition allowing only isolated simple double points (nodes) for
$\overline{Y}$. In fact these singularities can (at least locally)
be smoothed. The following results also hold:

\begin{theorem}
(Friedman \cite{Friedman1986}) \textit{\ If} $\phi $ \textit{\ is
of type} $I$ \textit{\ and the singularity is an ordinary double
point, then} $\overline{Y}$ \textit{\ is smoothable unless} $\phi
$ \textit{\ is the contraction of a single} $\P^{1}$ \textit{\ to
an ordinary double point.}
\end{theorem}

\begin{theorem}(Altmann, Gross, Schlessinger)
(\cite{alt1997, Gross1997(a), Gross1997(b), Schlessinger1971}

\begin{itemize}
\item  \textit{\ If }$\phi $ \textit{\ is of type} $II$ \textit{\ and} $\overline{Y}$
\textit{\ is} $\Q$\textit{--factorial, then} $\overline{Y}$
\textit{\ is smoothable unless} $E\cong \P^{2}$ \textit{\ or}
$E\cong \Bbb{F}_{1}$

\item  \textit{\ If} $\phi $ \textit{\ is of type} $III$ \textit{\ and} $\overline{Y}$
\textit{\ is} $\Bbb{Q}$\textit{--factorial, then} $\overline{Y}$
\textit{\ is smoothable unless} $\phi (E)\cong \P^{1}$ \textit{\
and} $E^3={7}, \ {8}$.
\end{itemize}
\end{theorem}

\vskip 0.1in

After Clemens' work (see \ref{conifold topology}), Reid
suggested that the birational classes of Calabi-Yau threefolds
would \textit{fit together into one irreducible family} (see \cite{Reid1987}). In fact
he speculated that transitions may connect a general Calabi-Yau
threefold to a non--K\"{a}hler analytic threefold with trivial
canonical class, Betti number $b_{2}=0$ and diffeomorphic to a
connected sum of $N$ copies of $S^{3}\times S^{3}$, where $N$ is
arbitrarily large. This conjecture is usually known as
\textit{Reid's fantasy}. There exist various pieces of evidence for this
conjecture (the Calabi-Yau web: see e.g. \cite{Avram1996},
\cite{Chiang1996}).

\phantom{xx}

 \subsection{Transitions and mirror symmetry}~ \vskip 0.1in
\index{mirror, symmetry} \noindent Assume that there exists a
transition
from $Y_{1}$ to $%
\widehat{Y}_{1}$, factorizing through a birational contraction
$\phi :Y_{1}\longrightarrow \overline{Y}_{1}$; assume also that
the \textit{mirror partners}
\index{mirror, partners} ${Y}_{2}$ of $%
Y_{1}$ and $%
\widehat{Y}_{2}$ of $\widehat{Y}_{1}$ exist (see, for example,
\cite{Morrison 1999}).

\noindent It is believed that the mirror partners $\widehat{Y}_2$
and ${Y}_2$ are also connected by a transition, which factorizes
through a birational contraction  $\phi^{\circ }
:\widehat{Y}_{2}\longrightarrow \overline{Y}_{2}$; the transition
between $\widehat{Y}_2$ and ${Y}_2$ is often called the ``reverse
transition". \index{transition, reverse} It is not known if this
conjecture holds; see for example \cite{Batyrev1998}, for the case
of the conifold transition.

 \begin{figure}[h]
     \begin{center}
     \scalebox{1}{\includegraphics{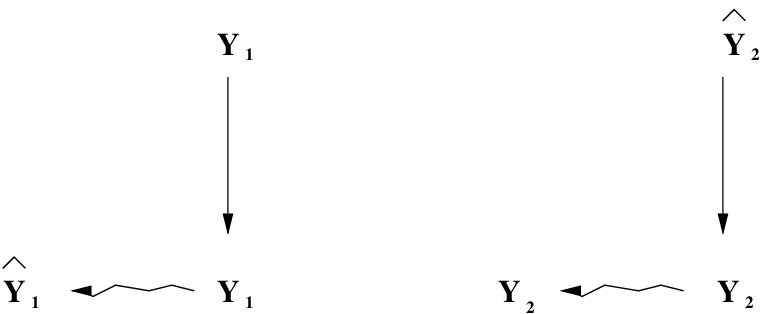}}
  \end{center}
     \end{figure}

\noindent The mirror symmetry exchanges the Hodge numbers
$h^{1,2}$ (representing the dimension of the complex moduli space)
with $h^{1,1}$ (the K\"{a}hler moduli space) of the Calabi-Yau
mirror partners; this exchange is consistent with a partner mirror
transition as we have seen in subsection \ref{conifold
topology}. \cite{GP} outlined an heuristic approach to
``continuously'' extending mirror symmetry to all the Calabi-Yau
threefolds belonging to the same connected component of the web
generated by conifold transitions. Actually \textit{if transitions
would connect to each other all Calabi-Yau threefolds}, which is a
rough version of Reid's fantasy, then it could give an
approach to establish mirror symmetry for all of them.

\noindent In the examples studied by \cite{CdOFKM} and \cite{Morrison 1999} $%
Y_{1}$, ${\widehat{Y}}_{1}$ and their mirrors are related by a
primitive contraction of type $III$ (see Appendix
\ref{Candelas'ex}).

\phantom{xx}

\subsection{Transitions, black holes etc.} ~ \vskip 0.1in

\index{black hole}The transitions among Calabi--Yau manifolds are
crucial also in the context of string theory, as they connect two
topologically distinct compactifications of a 10--dimensional type
II string theory (to 4--dimensional string vacua). Since, in spite
of the small number of consistent 10--dimensional string theories,
their Calabi-Yau compactifications give rise to a multitude of
4--dimensional topologically distinct string vacua, transitions may
result to be the suitable mathematical tool which is able to
restore a concept of \textit{uniqueness }in compactified string
theory when mirror symmetry and a version of Reid's fantasy (the
Calabi-Yau web) is assumed. The physical interpretation would then
be that \textit{two 4--dimensional topologically distinct string
vacua may be connected to each other by means of a black hole
condensation. } \index{black hole, condensation} This is the work
of \cite{Greene1995}, \cite{Strominger1995}.

\noindent Strominger gave a physical explanation of how to resolve
the conifold singularities of the moduli space of classical string
vacua by means of massless Ramond--Ramond (RR) black holes
\index{RR, black hole} (see Appendix \ref{blackholes}).

\noindent In \cite{Greene1995} the transformation of a massive
black hole into a massless one at the conifold model is called
\textit{condensation.} Not only conifold transitions have a
physical counterpart. For example a similar interpretation
involves transitions of type $II$ in the context of string--string
duality (see \cite{kmp}, \cite{Berglund1995}, \cite{Berglund1997}).

\noindent Transitions of Calabi--Yau manifolds also have a role in
5--dimensional supersymmetric theories (see for example
\cite{Morrison--Seiberg 1997}, \cite{dkv}).

\section{Chern--Simons theory}

\index{Chern--Simons, theory}

We discuss some basics of classical Chern-Simons theory (following
\cite{CS} and\cite{df}) and of its quantum version (following
\cite{Witten1989} and \cite{Labastida1999}).

 The first evidence for the conjecture comes from comparing the
 expansion of the Chern-Simons partition function (with and without knots),
 so the last section is dedicated to the computational
 aspects and link invariants.
 We start with a quick
review of the mathematical background for Chern--Simons theory,
principal bundles and connections:  Appendix \ref{diffdef}
contains more details.

\vskip 0.2in
 Let $\pi :P\rightarrow M$ be a principal $G$-bundle
with $G$ acting on the right (see definition \ref{principal bdl}).
In particular, for any $m\in M$, $\pi ^{-1}(m)\cong G$. The
differential of this map gives an isomorphism between the tangent
space $\pi ^{-1}(m)$ to each fiber at a point $p \in \pi
^{-1}(m)$:
\[
d\pi :T_{p}\pi ^{-1}\left( m\right) \stackrel{\cong }{\longrightarrow }%
T_{id}G\cong \mathfrak{g}
\]
Let $TP$ denote the tangent bundle of $P$:
\begin{definition}
The \textit{vertical bundle on }$P$
\index{bundle, vertical} is the vector sub--bundle $\mathcal{V}P$ of $TP$
given by $\ker \left( d\pi \right) $; that is, for every $p\in P$
\[
\mathcal{V}_{p}P:=\ker \left[ d_{p}\pi :T_{p}P\longrightarrow T_{\pi \left(
p\right) }M\right]
\]
\end{definition}

\noindent Then the vertical bundle $\mathcal{V}P$ associated with
the principal $G$--bundle $\left( P,\pi \right) $ is a vector
bundle whose standard fibre is the Lie algebra $\mathfrak{g}$
associated with $G$ (see remark \ref{g-fibre}).\textit{\ }

A connection is an infinitesimal version of a $G$-equivariant family of
sections of $\pi: P \to M$.

\begin{definition}
\label{def conn}A \textit{connection on a principal }$G$\textit{--bundle }$%
\left( P,\pi \right) $
\index{connection} is a vector sub-bundle $\mathcal{H}P$ of $TP$ such that
\begin{equation}
TP=\mathcal{H}P\oplus \mathcal{V}P  \label{splitting}
\end{equation}
and for every $p\in P$ and $\sigma \in G$%
\begin{equation}
d_{p}R\left( \sigma \right) \left( \mathcal{H}_{p}P\right) =\mathcal{H}%
_{p\sigma }P  \label{G-invariance}
\end{equation}
where $R$ is the right action of $G$ on $P$ (see definition \ref{LRactions}).
\end{definition}

\begin{definition}
\label{conn.forms}

\begin{enumerate}
\item  The \textit{connection form of a connection} $\mathcal{H}P$
\index{form, connection} is the $\mathfrak{g}$--valued 1--form
$A\in \Omega ^{1}\left( P,\mathfrak{g}\right) $ such that, for
every $p\in P$ and $u\in T_{p}P
$%
\begin{equation}
A_{p}u:=\left( d_{id}\lambda _{p}\right) ^{-1}\left(
\mathcal{V}_{p}u\right) \in T_{id}G\cong \mathfrak{g}
\label{g-connection}
\end{equation}
where $\lambda _{p}:G\stackrel{\cong }{\longrightarrow }\pi ^{-1}\left( \pi
\left( p\right) \right) \subset P$ is the diffeomorphism given by $\lambda
_{p}\left( \sigma \right) :=p\sigma $. It is a characteristic form of the
connection $\mathcal{H}P$ since $\mathcal{H}P=\ker A$ (see proposition \ref
{connections}).

\item  The \textit{curvature form of a connection }$\mathcal{H}P$
\index{form, curvature} is the $\mathfrak{g}$--valued 2--form
$\Omega \in \Omega ^{2}\left( P,\mathfrak{g}\right) $ defined by:
\begin{equation}
\Omega _{p}\left( u,v\right) :=-A_{p}\left[ U,V\right] _{p},\ \ \ \forall
p\in P,\ u,v\in T_{p}P  \label{g-curvature}
\end{equation}
where $U,V$ are any horizontal vector fields on $P$ extending the horizontal
parts $\mathcal{H}_{p}u$ and $\mathcal{H}_{p}v$ of $u$ and $v$ respectively
(recall the splitting $\left(
\text{\ref{split}}\right) $).
\end{enumerate}
\end{definition}

\begin{definition}
A \textit{gauge transformation of }$P$
\index{gauge transformation} is an automorphism $\varphi $ of $P$ which
induces the identity map on the base manifold $M$.
\end{definition}

\noindent Gauge transformations on $P$ form a group $\mathcal{G}_{P}$,
\index{gauge transformations, group of} and $\left(
\text{\ref{gauge on connection}}\right) $ defines an action of $\mathcal{G}%
_{P}$ on the affine space of connections $\mathcal{A}_{P}$ (see proposition
\ref{connections}).

\begin{definition}
Let $\gamma :I:=\left[ 0,1\right] \longrightarrow M$ be a loop with base
point $m\in M$ and let $\widetilde{\gamma }_{p}:I\longrightarrow P$ be the
unique \textit{horizontal lift of }$\gamma $ \textit{with initial point }$%
p\in P$, i.e. such that
\[
d\widetilde{\gamma }_{p}\left( TI\right) \subset \mathcal{H}P\text{\quad and}%
\quad \widetilde{\gamma }_{p}\left( 0\right) =p
\]
Define a diffeomorphism of the fibre $\pi ^{-1}\left( m\right) $ by
\begin{equation}
\begin{array}{cccc}
h_{\gamma }: & \pi ^{-1}\left( m\right)  & \longrightarrow  & \pi
^{-1}\left( m\right)  \\
& p & \longmapsto  & \widetilde{\gamma }_{p}\left( 1\right)
\end{array}
\label{holgamma}
\end{equation}
Then:
\begin{equation}
\text{Hol}_{\mathcal{H}P}\left( m\right) :=\{h_{\gamma }:\gamma \text{ is a
loop based at }m\}  \label{Hol(m)}
\end{equation}
is a group (with the composition of morphisms), called \textit{the holonomy
group of the connection }$\mathcal{H}$\textit{$P$}
\index{holonomy group} \textit{at }$m\in M$.

If the base manifold $M$ is connected all these groups are isomorphic by $%
\left(
\text{\ref{connect-hol}}\right) $. Then Hol$_{\mathcal{H}P}$ is called
\textit{the holonomy group of the connection }$\mathcal{H}$\textit{$P$.}
\end{definition}

\noindent Note that for every $p\in P$ it is possible to identify Hol$_{%
\mathcal{H}P}\left( \pi \left( p\right) \right) $ with the subgroup of $G$%
\begin{equation}
G_{\mathcal{H}P}\left( p\right) :=\left\{ \sigma _{\gamma }\left( p\right)
\in G:h_{\gamma }\left( p\right) =p\sigma _{\gamma }\left( p\right) \text{
and }h_{\gamma }\in \text{Hol}_{\mathcal{H}P}\left( \pi \left( p\right)
\right) \right\}  \label{hol-sbgroup}
\end{equation}
If $p,q\in \pi ^{-1}\left( m\right) $ then $G_{\mathcal{H}P}\left( p\right) $%
 and $G_{\mathcal{H}P}\left( q\right) $ are conjugate subgroups and they
coincide if $p$ and $q$ can be joined by a horizontal curve in $P$.

\begin{definition}
\textit{\ The restricted holonomy group of the connection }$\mathcal{H}$%
\textit{$P$}
\index{holonomy group, restricted} \textit{at }$m\in M$
\begin{equation}
H_{\mathcal{H}P}^{\left( o\right) }\left( m\right) \subset
\text{Hol}_{\mathcal{H}P}\left( m\right)   \label{restr.holonomy}
\end{equation}
is defined by considering homotopically trivial loops based at $m$.
\end{definition}

\noindent As before, if $M$ is connected we can define the restricted
holonomy group $H_{\mathcal{H}P}^{\left( o\right) }\subset \ $Hol$_{\mathcal{%
H}P}$. Moreover for every $p\in P$ we can identify the restricted holonomy
subgroup $H_{\mathcal{H}P}^{\left( o\right) }\left( \pi \left( p\right)
\right) $ with a suitable subgroup $G_{\mathcal{H}P}^{\left( o\right)
}\left( p\right) \subset G_{\mathcal{H}P}\left( p\right) \subset G$.

\phantom{xx}

\subsection{Classical Chern--Simons action}~ \vskip 0.1in

 Let us assume the base manifold $M=\pi \left(
P\right) $ to be a smooth and compact 3--manifold. Let
$\mathcal{A}_{P}$ be the affine space of all possible
connections on $P$ and choose $A \in \mathcal{A}_{P}$ with
associated connection $\mathcal{H}P=\ker A $. If $\Omega \in
\Omega ^{2}\left( P,\mathfrak{g}\right) $ is the
$\mathfrak{g}$--valued curvature 2--form of the chosen connection
then
\[
\Omega \wedge \Omega \in \Omega ^{4}\left( P,\mathfrak{g}\otimes
\mathfrak{g}\right)
\]

\begin{definition}
\textit{The Chern--Weil }
\index{form, Chern--Weil}%
\index{Chern--Weil form} $4$\textit{--form associated with the Killing form }%
$\left\langle \ ,\ \right\rangle $ (see definition \ref{killing}) is\textit{%
\ }$\left\langle \Omega \wedge \Omega \right\rangle \in $ $\Omega ^{4}\left(
P\right) $.
\end{definition}

\begin{definition}
\label{CSform}\label{CSform-def} \textit{A Chern--Simons form}
\index{form, Chern--Simons}%
\index{Chern--Simons, form} is an anti--derivative $\alpha \in \Omega
^{3}\left( P\right) $ of $\left\langle \Omega \wedge \Omega \right\rangle $.
\end{definition}

\begin{proposition}
Let $\alpha :=\left\langle A\wedge \Omega \right\rangle -%
\frac{1}{6}\left\langle A\wedge \left[ A,A\right] \right\rangle $. Then:

\begin{enumerate}
\item  $d\alpha =$\textit{\ }$\left\langle \Omega \wedge \Omega
\right\rangle $,

\item  \textit{\ if }$\varphi $\textit{\ is a gauge transformation of }$P$,
\begin{equation}
\left( \delta \varphi \right) \alpha =\alpha -\frac{1}{6}\left\langle \phi
\wedge \left[ \phi ,\phi \right] \right\rangle +d\left\langle \left(
Ad_{\sigma _{\varphi }^{-1}}\circ A\right) \wedge \phi \right\rangle
\label{gauge on CSform}
\end{equation}
\textit{where }$\delta $ \textit{is the codifferential, }$\sigma _{\varphi }$%
\textit{\ is associated with }$\varphi $\textit{\ as in }$\left( \text{\ref
{gauge element}}\right) $\textit{, }$\phi :=$\textit{\ }$\left( \delta
\sigma _{\varphi }\right) \left( \delta \lambda \right) A$ \textit{and }$%
\left( \delta \lambda \right) A$\textit{\ is the Maurer--Cartan form of the
connection }$\mathcal{H}P$ \textit{as defined in }$\left( \text{\ref
{MCconnection}}\right) $.

\item  \textit{If }$\alpha ^{\prime }$\textit{\ is a Chern-Simons form, the }%
$3$\textit{--form }$\left( \delta \varphi \right) \alpha ^{\prime }-\alpha
^{\prime }+\frac{1}{6}\left\langle \phi \wedge \left[ \phi ,\phi \right]
\right\rangle $ \textit{is exact.}
\end{enumerate}
\end{proposition}

\noindent The proof follows directly by the definition \ref{CSform-def} of $%
\alpha $ and by the gauge action on connections $\left( \text{\ref{gauge on
connection}}\right) $. By $\left( \text{\ref{gauge on curvature}}\right) $
and the $Ad$--invariance (see $\left( \text{\ref{Ad}}\right) $) of the
Killing form \textit{the Chern--Weil form }$\left\langle \Omega \wedge
\Omega \right\rangle $\textit{\ is gauge invariant. } Moreover:

\begin{proposition}
$\alpha ^{\prime }-\left( \delta \varphi \right) \alpha ^{\prime }$ \textit{%
defines a cohomology class}
\[
\left( \delta \sigma _{\varphi }\right) \Phi _{A}\in H^{3}\left( P,\R%
\right) ,
\]
\textit{which is independent of the chosen Chern--Simons form} $\alpha
^{\prime }$. \textit{We can also assume that}
\begin{equation}
\rho \Phi _{A}\in H^{3}\left( G,\Bbb{Z}\right)   \label{Zcondition}
\end{equation}
\textit{for a suitable real number} $\rho $.
\end{proposition}

\noindent In fact, the 3--form $\alpha ^{\prime }-\left( \delta
\varphi \right) \alpha ^{\prime }$ is closed for every gauge
transformation $\varphi $ and any Chern--Simons form $\alpha
^{\prime }$. Also it is the image by the codifferential $\delta
\sigma _{\varphi }$ of the cohomology class $\Phi _{A}\in
H^{3}\left( G,\R\right) $ associated with the closed 3--form
\[
\frac{1}{6}\left\langle \left( \delta \lambda \right) A\wedge \left[ \left(
\delta \lambda \right) A,\left( \delta \lambda \right) A\right]
\right\rangle \in \Omega ^{3}\left( G\right)
\]
Note that the choice of $\rho \in \R$ depends only on the connection $%
\mathcal{H}P$.

\begin{definition}
If there exist a global section
\[
s:M\longrightarrow P,
\]
the \textit{Chern--Simons Lagrangian on }$M$
\index{Lagrangian, Chern--Simons}
\index{Chern--Simons, Lagrangian} \textit{is the }$3$\textit{--form }
\begin{equation}
\mathcal{L}\left( A,s\right) :=\rho \left( \delta s\right) \alpha \in \Omega
^{3}\left( M\right)   \label{CSlagrangian}
\end{equation}
and the associated \textit{Chern--Simons action }
\index{Chern--Simons, action} is obtained by integrating it over $M$%
\begin{equation}
S\left( \mathcal{L}\right) :=\int_{M}\mathcal{L}\left( A,s\right)
\label{CSaction}
\end{equation}
\end{definition}

\begin{remark}
\begin{enumerate}
\item  The existence of a section means that $P$ is \textit{parallelizable},
which is the case for example when $G$ is simply connected (see
\cite{df}, lemma 2.1 for a proof of this fact.)

\item  By Stokes' theorem \textit{the Chern--Simons action }$S$\textit{\
does not depend on the choice of the Chern--Simons form }$\alpha $ \textit{%
when }$M$\textit{\ is assumed to be without boundary.}

\item  For any gauge transformation $\varphi $, the 3--form $\mathcal{L}%
\left( A,s\right) -\left( \delta \varphi \right) \mathcal{L}\left(
A,s\right) $\textit{\ }defines the integral cohomology class
\[
\rho \delta \left( \sigma _{\varphi }\circ s\right) \Phi _{A}\in H^{3}\left(
M,\Bbb{Z}\right)
\]
hence
\begin{equation}
S\left( \mathcal{L}\right) -S\left( \left( \delta \varphi \right) \mathcal{L}%
\right) =\rho \int_{M}\delta \left( \sigma _{\varphi }\circ s\right) \Phi
_{A}\in \Bbb{Z}  \label{quantization}
\end{equation}

\item  For the particular case $G=SU\left( 2\right) $ the integral bilinear
forms on $\mathfrak{g}=\mathfrak{su}_{2}$ are parameterized by
$k\in \Bbb{Z}$ as follows:
\[
\forall X,Y\in \mathfrak{su}_{2}\quad \left\langle X,Y\right\rangle _{k}=%
\frac{k}{8\pi ^{2}}\text{ tr}\left( XY\right) .
\]
Then the real coefficient in $\left( \text{\ref{Zcondition}}\right) $ can be
given by $\rho :=\left( 8\pi ^{2}\right) ^{-1}$ and the Chern--Simons
Lagrangian $\left( \text{\ref{CSlagrangian}}\right) $ becomes
\[
\mathcal{L}\left( A,s\right) =\frac{1}{8\pi ^{2}}\text{ tr}\left( A^{\prime
}\wedge dA^{\prime }+\frac{2}{3}A^{\prime }\wedge A^{\prime }\wedge
A^{\prime }\right)
\]
where $A^{\prime }:=\left( \delta s\right) A$ (see section 6 in
\cite{df}). This is the typical shape of a Chern--Simons
Lagrangian usually adopted in the
physics literature, although the gauge group $G$ is more general than $%
SU\left( 2\right) $.
\end{enumerate}
\end{remark}

\begin{proposition}
\textit{The} \textit{\ }Chern--Simons action
\begin{equation}
S\left[ A\right] :=\exp \left( ik2\pi S\left( \mathcal{L}\right) \right)
\label{expCS}
\end{equation}
\textit{is well defined and gauge invariant, where} $k\in \Bbb{Z}$ \textit{%
is called} the level of the theory\textit{. Furthermore,} $S\left[ A\right] $
\textit{depends only on the choice of the gauge equivalence class of
connections} $\left[ A\right] \in \mathcal{A}_{P}/\mathcal{G}_{P}$, \textit{%
where} $\mathcal{G}_{P}$ \textit{acts on} $\mathcal{A}_{P}$ \textit{as in}
$\left( \text{\ref{gauge on connection}}\right) $.
\end{proposition}

\noindent In fact any two sections of $P$ are related by a gauge
transformation and the assumption $\left( \text{\ref{Zcondition}}\right) $
holds.

\noindent From the physical point of view it is relevant to point out the
\textit{quantization law }expressed by $\left( \text{\ref{quantization}}%
\right) $ and $\left( \text{\ref{expCS}}\right) $. The real factor
$\rho $ defined in $\left( \text{\ref{Zcondition}}\right) $ may be
considered to be a \textit{normalizing factor of the Killing form
of }$\mathfrak{g}$. Then we can write $\left(
\text{\ref{quantization}}\right) $ as:
\[
S\left( \mathcal{L}\right) -S\left( \left( \delta \varphi \right) \mathcal{L}%
\right) =\int_{M}\delta \left( \sigma _{\varphi }\circ s\right) \Phi _{A}\in
\Bbb{Z}.
\]
We can also relate any gauge transformation $\varphi $ with a map $%
M\rightarrow G$ by taking $\sigma _{\varphi }\circ s$. In this way we get an
immersion of the gauge group $\mathcal{G}_{P}$ into the group of maps from $%
M $ to $G$. $\int_{M}\delta \left( \sigma _{\varphi }\circ
s\right) \Phi _{A} $ is called the \textit{winding number of the
gauge transformation } \index{gauge transformation, winding number
of} $\varphi $\textit{. }Since this number is homotopically
invariant it plays the role of counting homotopy classes of gauge
transformations, giving two relevant consequences:

\begin{enumerate}
\item  the Chern--Simons action $\left(
\text{\ref{CSaction}}\right) $ is invariant under any gauge transformation
homotopically equivalent to the identity,

\item  as in Dirac's famous work on magnetic monopoles, the integer $k$
in $\left( \text{\ref{expCS}}\right) $ turns out to be closely related to the
central charge of the theory. Moreover, in the quantum field theory defined
by the following partition function $\left( \text{\ref{partition}}\right) $ $%
k^{-1}$ is proportional, for large $k$, to the square $\lambda $ of the
\textit{coupling constant }of the theory (see $\left( \text{\ref{lambda}}%
\right) $).
\end{enumerate}

\phantom{xx}

\subsection{Chern--Simons quantum field theory}~ \vskip 0.1in
\begin{definition}
The \textit{Chern--Simons partition function}
\index{Chern--Simons, partition function} is the Feynman integral of the
Chern--Simons action $\left(
\text{\ref{expCS}}\right) $ taken over all the gauge equivalence classes of
connections:
\begin{equation}
Z\left( M\right) :=\int_{\mathcal{A}_{P}/\mathcal{G}_{P}}S\left[ A\right]
D\left[ A\right] .  \label{partition}
\end{equation}
\end{definition}

\noindent This defines the \textit{Chern-Simons quantum field
theory} \index{Chern--Simons, QFT} (see for example
\cite{Deligne1999}) whose fields are precisely the elements of
$\mathcal{A}_{P}/\mathcal{G}_{P}$.

\begin{definition}
Let $K$ be a knot in $M$, i.e. an embedding of the circle $S^{1}$ and $R$ a
representation of $G$. \textit{The Wilson line} $W_{K}^{R}$
\index{Wilson line, functional} is the functional
\begin{equation}
W_{K}^{R}:\mathcal{A}_{P}/\mathcal{G}_{P}\longrightarrow \R
\label{Wilson}
\end{equation}
where $W_{K}^{R}\left[ A\right] :=$tr$_{R}\left( h_{K}\right) $ and $h_{K}$
is the holonomy around $K$.
\end{definition}

\noindent Note that the real number tr$_{R}\left( h_{K}\right) $ is well
defined for any representation $R$ of $G$. $K$ can be thought as a closed
loop in $M$; for every point $m\in K$ we obtain an element $h_{K}\in $Hol$_{%
\mathcal{H}P}\left( m\right) $ as in $\left(
\text{\ref{holgamma}}\right) $. If $M$ is connected $h_{K}$ does not depend
on the choice of $m\in K$ since we can proceed as in $\left( \text{\ref
{connect-hol}}\right) $ to obtain $h_{K}\in $Hol$_{\mathcal{H}P}$. By $%
\left( \text{\ref{hol-sbgroup}}\right) $ $h_{K}$ defines a conjugacy class
in $G$.

\noindent The Wilson line are metric independent (i.e. covariant) and gauge
invariant functionals of the fields; they are then \textit{observables}
\index{Wilson line, observable} of the theory.

\noindent Since tr$_{R}\left( h_{K}\right) $ is gauge invariant, we define:

\begin{definition}
The \textit{unnormalized expectation value} is formally assigned by the
Feynman integral
\begin{equation}
Z\left( M;K,R\right) :=\int_{\mathcal{A}_{P}/\mathcal{G}_{P}}S\left[
A\right] W_{K}^{R}D\left[ A\right]   \label{unnorm.exp.value}
\end{equation}
and its \textit{expectation value}
\index{Wilson line, expectation value} is given by
\begin{equation}
\left\langle W_{K}^{R}\right\rangle :=Z\left( M;K,R\right) /Z\left( M\right)
\label{exp.value}
\end{equation}
\end{definition}

\noindent If we now consider a link $L$ in $M$, i.e. the union of $r\geq 1$
oriented and non--intersecting knots $\left\{ K_{i}\right\} _{i=1}^{r}$ in
the oriented manifold $M$ and a collection of irreducible representations $%
\mathcal{R}:=\left\{ R_{i}\right\} _{i=1}^{r}$ of $G$, one for each knot $%
K_{i}$, we have:

\begin{definition}
The \textit{correlation function }
\index{Chern--Simons, correlation function} of our quantum field theory is
\begin{equation}
Z\left( M;L,\mathcal{R}\right) :=\int_{\mathcal{A}_{P}/\mathcal{G}%
_{P}}S\left[ A\right] \prod_{i=1}^{r}W_{K_{i}}^{R_{i}}D\left[ A\right]
\label{correlation}
\end{equation}
\end{definition}

\phantom{xx}

\subsection{\label{HAMILTON}The Hamiltonian
formulation of the Chern--Simons QFT (following Witten's canonical
quantization)} ~ \vskip 0.1in

\index{Hamiltonian, formulation of Chern--Simons QFT}Although the
mathematical definitions of path integrals in $\left(
\text{\ref{partition}}\right) $, $\left( \text{\ref{unnorm.exp.value}}%
\right) $ and $\left( \text{\ref{correlation}}\right) $ are quite delicate,
the explicit integrals are calculated in \cite{Witten1989}. Witten first uses the
\textit{stationary--phase approximation }%
\index{stationary--phase approximation}
\index{approximation, stationary--phase}in the ``classical limit'' $%
k\rightarrow \infty $ and then%
\index{quantization, Witten's canonical} \textit{canonical
quantization}. Here we present the basic ideas of this second
method. A very useful and pleasant reference on the
argument is \cite{Atiyah1990a}, to which we refer the reader for a
deeper understanding. We will not discuss the stationary--phase
approximation since it lies outside the aim of the present work,
although its relevance is fundamental in giving the confirmation
that the partition functions introduced by the Feynman approach in the
previous section are the same as those we will evaluate in the next section
by the Hamiltonian approach: see the first part of section 2 in
\cite{Witten1989} and section 7.2 in \cite{Atiyah1990a}

The main purpose in QFT of a Feynman path integral is to provide a
relativistically invariant approach, since this is a fundamental property of
the \textit{Lagrangian density }which in our case is expressed by the
Chern--Simons action $\left(
\text{\ref{CSaction}}\right) $ multiplied by $2\pi k$. If we want to
focus on a time--evolution in the theory we have to break the relativistic
symmetry by constructing a \textit{time--evolution operator} $\exp \left(
itH\right) $ in a certain ``Hilbert'' space $\mathcal{H}$ representing the
space of physical states%
\index{Hilbert space, of physical states}. The generator $H$ is the \textit{%
Hamiltonian operator }of the theory%
\index{Hamiltonian operator}
\index{time--evolution operator}. In general there are formal rules which
allows one to produce the space $\mathcal{H}$ and the Hamiltonian $H$ of a QFT
whose partition function is known.

\noindent In the case of Chern--Simons QFT the space--time is
represented by the 3--manifold $M$. We can separate out space and
time by ``cutting'' $M$ along a surface $\Sigma $. Near the cut
$M$ looks like $\Sigma \times \R $, giving us the desired
separation of space and time. Let us then reduce to considering the
particular case $M=$ $\Sigma \times \R$ which can be treated by
means of \textit{canonical quantization} to construct the physical
space $\mathcal{H=H}\left( \Sigma \right) $ of the Chern--Simons
theory \textit{quantized on }$\Sigma $. More precisely this means
to ``quantize'' the space of classical solutions, which are the
critical fields of the Chern--Simons action $\left(
\text{\ref{CSaction}}\right) $.

\begin{proposition}
\textit{The space of classical solutions of Chern--Simons theory is the
subspace of gauge equivalence classes of flat connections in }$\mathcal{A}%
_{P}/\mathcal{G}_{P}$\textit{\, which can be naturally identified with the
following}
\[
\mathcal{M}_{M}:=\hom \left( \pi _{1}\left( M\right) ,G\right) /G
\]
\textit{where }$G$\textit{\ acts by conjugation (}See \cite{df},
proposition 3.5 for more details).
\end{proposition}

\noindent The statement follows by $\left( \text{\ref{CSlagrangian}}\right) $
and the fact that $\alpha $ is by definition an anti--derivative of $%
\left\langle \Omega \wedge \Omega \right\rangle $. In fact
\begin{equation}
dS\left( \mathcal{L}\left( A,s\right) \right) =0\quad \Longleftrightarrow
\quad \Omega =0  \label{Euler-Lagrange}
\end{equation}
i.e. the latter is the%
\index{Euler--Lagrange equation, of Chern--Simons theory}
\index{Chern--Simons, classical solutions} \textit{Euler--Lagrange equation }%
of the classical Chern--Simons theory whose solutions are given by \textit{%
flat connections}. See \cite{df}, proposition 3.1 for details on
differentiating. Note that by $\left(
\text{\ref{struct.eqn}}\right) $ this Euler--Lagrange equation
involves only first order derivatives of the fields. This is a
peculiarity of Chern--Simons gauge theory together with the
independence of the choice of the metric. Since the restricted
holonomy subgroups $\left( \text{\ref {restr.holonomy}}\right) $
of a flat connections are always trivial it is possible to define
a morphism
\[
\pi _{1}\left( M\right) \longrightarrow \ \text{Hol}_{\mathcal{H}P}
\]
(see e.g. \cite{Poor 1981}, proposition 2.40). By recalling $\left( \text{\ref
{hol-sbgroup}}\right) $ we actually get a morphism from $\pi _{1}\left(
M\right) $ to $G$ which is well defined up to conjugation. On the contrary a
similar equivalence class of morphisms suffices to determine a flat
connection on $P$.

\noindent Since we are in the particular case $M=\Sigma \times \R$
our space of classical solutions reduces to
\begin{equation}
\mathcal{M}_{\Sigma }:=\hom \left( \pi _{1}\left( \Sigma \right) ,G\right) /G
\label{class.sol.}
\end{equation}
This space is not dependent on the time variable described by $\R$
implying that we actually have no time--evolution in our theory
i.e. \textit{we have no dynamics and all is purely topological}:
\textit{hence the Hamiltonian }$H$\textit{\ must be trivial.}

\noindent The following result allows one to ``quantize''
$\mathcal{M}_{\Sigma }$:

\begin{theorem}
(\cite{Narasimhan--Seshadri 1965}, \cite{d}) \textit{The space of
classical solutions }$\mathcal{M}_{\Sigma }$\textit{\ is
homeomorphic to the moduli space} $M_{\tau }$ \textit{of
holomorphic} $G$\textit{--bundles over the Riemann surface}
$\Sigma _{\tau }$ \textit{obtained by the }choice\textit{\ of a
complex structure} $\tau $ \textit{on} $\Sigma $. \textit{On
}$M_{\tau }
$ \textit{we have a natural choice} \textit{of a holomorphic line bundle} $L$%
. \textit{The finite dimensional complex vector space}
\begin{equation}
\mathcal{H}_{\tau }^{k}\left( \Sigma \right) :=H^{0}\left( M_{\tau
},L^{\otimes k}\right)   \label{Hilbert}
\end{equation}
\textit{of global holomorphic sections of} $L^{\otimes k}$ \textit{gives the
Hilbert space of the quantized theory at level} $k$%
\index{Hilbert space, of Chern--Simons QFT}.
\end{theorem}

\noindent When $G=SU\left( N\right) $ the moduli space $M_{\tau }$ turns out
to be a projective algebraic variety. Hence we have the natural choice $L:=%
\mathcal{O}_{M_{\tau }}\left( 1\right) $ i.e. the line--bundle associated
with the hyperplane section. Otherwise, when $G$ is more general, the choice
of the complex structure $\tau $ on $\Sigma $ gives a natural complex
structure on the infinite dimensional affine space $\mathcal{A}_{P}$. The
moduli space $M_{\tau }$ can then be identified with the%
\index{symplectic quotient} \textit{symplectic quotient }$\mathcal{A}_{P}//%
\mathcal{G}_{P}$\textit{\ }(see \cite{Atiyah1990a}, chapter 4, for
a definition) under the action $\left( \text{\ref{gauge on
connection}}\right) $ of the gauge group $\mathcal{G}_{P}
$ (see \cite{AB} for the details). On $\mathcal{A}_{P} $ the \textit{%
Quillen line--bundle }$\mathcal{L}$\textit{\ }(see \cite{Quillen 1986}),
\index{bundle, Quillen's line--}
\index{Quillen line--bundle}whose curvature is $-2\pi i$ times the
K\"{a}hler form of $\mathcal{A}_{P}$, descends to give a well--defined
line--bundle $L$ on $M_{\tau }$.

\noindent The crucial point now is that the vector space $%
\mathcal{H}_{\tau }^{k}\left( \Sigma \right) $ apparently
\textit{depends on the choice
of the complex structure} $\tau $ \textit{on} $\Sigma $, which goes against the desired
general covariance of our theory. Actually $\mathcal{H}_{\tau }^{k}\left(
\Sigma \right) $ varies holomorphically with $\tau $ giving rise to a
holomorphic vector bundle over the moduli space of compact Riemann surfaces
of fixed genus which turns out to admit a canonical \textit{projectively
flat connection }%
\index{connection, projectively flat}\textit{which permits one to
identify the fibers up to a scalar factor. }This fact can be
proved in several ways, as described in chapter 6 of
\cite{Atiyah1990a}. See also \cite{Hitchin 1990} and
\cite{Axelrod1991} for more details.

\noindent The choice $\left(
\text{\ref{Hilbert}}\right) $  then gives rise to a%
\index{modular functor} \textit{modular functor}
\begin{equation}
\Sigma \longrightarrow \mathcal{H}^{k}\left( \Sigma \right)
\label{mod.functor}
\end{equation}
in the spirit of a \textit{rational conformal field theory}
\index{field theory, rational conformal}
\index{rational conformal field theory}as defined in \cite{Segal1988}: such a
functor is well defined \textit{up to a scalar factor. }It is a particular
case of a%
\index{TQFT} \index{field theory, topological quantum}
\textit{topological quantum field theory.} Let us now briefly
recall what it is as axiomatized in \cite{Atiyah1989}. The
interested reader may also consider chapter 2 in
\cite{Atiyah1990a} and appendix B.6 in \cite{CoxKatz1999} for
some short reviews on the subject and \cite{Quinn 1995} for a broader
treatment.

\begin{definition}
\label{TQFT}(\textit{Axiomatic TQFT}) A $\left( d+1\right) $--dimensional
topological quantum field theory is a functor $Z$ which associates

\begin{itemize}
\item  with each compact oriented $d$--dimensional manifold $\Sigma $ a
finite--dimensional complex vector space $Z_{\Sigma }$,

\item  with each compact oriented $\left( d+1\right) $--dimensional manifold
$M$ whose boundary is $\partial M=\Sigma $ a vector $Z\left( M\right) \in
Z_{\Sigma }$,
\end{itemize}

and which satisfies the following axioms:

\begin{enumerate}
\item  \textit{(Involutory) }if $%
\overline{\Sigma }$ denotes $\Sigma $ with the opposite orientation and $%
Z_{\Sigma }^{*}$ denotes the dual vector space of $Z_{\Sigma }$ then
\[
Z_{\overline{\Sigma }}=Z_{\Sigma }^{*}
\]

\item  \textit{(Multiplicativity) }if $\amalg $ denotes the disjoint union
of $d$--manifolds then
\[
Z_{\Sigma _{1}\amalg \Sigma _{2}}=Z_{\Sigma _{1}}\otimes Z_{\Sigma _{2}}
\]

\item  \textit{(Associativity) }if $\partial M_{1}=\overline{\Sigma }%
_{1}\amalg \Sigma _{2}$, $\partial M_{2}=\overline{\Sigma }_{2}\amalg \Sigma
_{3}$ and $M=M_{1}\cup _{\Sigma _{2}}M_{2}$ is the gluing of $M_{1}$ and $%
M_{2}$ along $\Sigma _{2}$ then
\[
Z\left( M\right) =Z\left( M_{2}\right) \circ Z\left( M_{1}\right)
\]
where by the previous axioms
\begin{eqnarray*}
Z\left( M_{1}\right)  &\in &Z_{\Sigma _{1}}^{*}\otimes Z_{\Sigma _{2}}=\hom
_{\C}\left( \Sigma _{1},\Sigma _{2}\right)  \\
Z\left( M_{2}\right)  &\in &Z_{\Sigma _{2}}^{*}\otimes Z_{\Sigma _{3}}=\hom
_{\C}\left( \Sigma _{2},\Sigma _{3}\right)  \\
Z\left( M\right)  &\in &Z_{\Sigma _{1}}^{*}\otimes Z_{\Sigma _{3}}=\hom _{%
\C}\left( \Sigma _{1},\Sigma _{3}\right)
\end{eqnarray*}

\item  \textit{(Unit) }if the empty set is considered as a compact $d$%
--dimensional oriented manifold then
\[
Z_{\emptyset }=\C
\]

\item  \textit{(Identity) }if $I$ denotes the oriented interval $\left[
0,1\right] $ let us consider the product $\left( d+1\right) $--manifold $%
\Sigma \times I$ whose boundary is $\partial \left( \Sigma \times I\right) =$
$\overline{\Sigma }\amalg \Sigma $; then
\[
Z\left( \Sigma \times I\right) =\Bbb{I}\in \hom _{\C}\left( \Sigma
,\Sigma \right)
\]
where $\Bbb{I}$ is the identity endomorphism of $\Sigma $.
\end{enumerate}
\end{definition}

Let us now come back to the Hamiltonian formulation of Chern--Simons quantum
field theory. In $\left( \text{\ref{mod.functor}}\right) $ we defined a
correspondence
\[
Z:\Sigma \longmapsto Z_{\Sigma }:=\mathcal{H}^{k}\left( \Sigma \right)
\]
between a compact surface $\Sigma \subset M$ and the finite dimensional
complex vector space of ``physical states'' of the level $k$ theory
quantized along $\Sigma $ by ``canonical quantization''%
\index{quantization, canonical}. This turns out to give a TQFT giving the
Hamiltonian interpretation of the partition function unrigorously expressed
by the path integral in $\left(
\text{\ref{partition}}\right) $. Precisely, by writing
\begin{eqnarray}
M &=&M_{1}\cup _{\Sigma }M_{2}  \label{chopping} \\
\partial M_{1} &=&\overline{\emptyset }\amalg \Sigma   \nonumber \\
\partial M_{2} &=&\overline{\Sigma }\amalg \emptyset   \nonumber
\end{eqnarray}
axioms 1,2,3 and 4 give
\begin{equation}
Z\left( M\right) =Z\left( M_{2}\right) \circ Z\left( M_{1}\right) \in \hom _{%
\C}\left( \C,\C\right) =\C  \label{hamiltonian partition}
\end{equation}
This is the mathematically well defined evaluation of the  partition
function. \index{partition function, Hamiltonian formulation of} It is
completely topological and the scalar indeterminacy in defining $Z_{\Sigma
}$ does not influence its value: actually $Z\left( M\right) $ does not
even depend on the choice of $\Sigma $ since $\partial M=\emptyset $ and $%
Z\left( M\right) \in Z_{\emptyset }$.

In order to perform an analogous Hamiltonian interpretation of the
correlation function\index{correlation function, Hamiltonian formulation of}
$Z\left( M;L,\mathcal{R}\right) $ ``defined'' by the path integral in $\left(
\text{\ref{correlation}}\right) $ we have to relativize the definition of
the TQFT $Z$ to the triple $\left( M,L,\mathcal{R}\right) $ given by a
3--manifold $M$ and a link $L\subset M$ marked by a collection of
irreducible representations $\mathcal{R}$ of $G$. Let us assume $L$ to be
transverse to $\partial M=\Sigma $ so that it gives a collection $\partial L$
of signed points in $\Sigma $. Moreover we can mark $\partial L$ by a
collection $\partial \mathcal{R}$ of irreducible representations of $G$
induced by representations in $\mathcal{R}$. Let us write
\begin{equation}
\partial \left( M,L,\mathcal{R}\right) =\left( \Sigma ,\partial L,\partial
\mathcal{R}\right)   \label{tern boundary}
\end{equation}
and then relativize $Z$ by defining it as a functor which associates

\begin{itemize}
\item  with each $d$--dimensional triple $\left( \Sigma ,\partial L,\partial
\mathcal{R}\right) $ a finite--dimensional complex vector space $Z_{\left(
\Sigma ,\partial L,\partial \mathcal{R}\right) }$,

\item  with each $\left( d+1\right) $--dimensional triple $\left( M,L,\mathcal{%
R}\right) $, whose boundary is as in $\left( \text{\ref{tern boundary}}%
\right) $, a vector $Z\left( M;L,\mathcal{R}\right) \in Z_{\left( \Sigma
,\partial L,\partial \mathcal{R}\right) }$,
\end{itemize}

\noindent and which satisfies the axioms 1, 2, 3, 4, 5 of definition \ref
{TQFT}. The crucial point now is to relativize $\left( \text{\ref
{mod.functor}}\right) $ to give an analogous definition of $Z_{\left( \Sigma
,\partial L,\partial \mathcal{R}\right) }$. Recall that by $\left( \text{\ref
{hol-sbgroup}}\right) $ the choice of a point $p\in \partial L\subset \Sigma
=\partial M$ determines \textit{a conjugacy class }in $G$. Since $p$ is
marked by an irreducible representation in $\partial \mathcal{R}$ the order
of such a conjugacy class turns out to be the level $k$. Hence the
collection $\partial L$ of marked points in $\Sigma $ gives rise to a set $%
C_{\partial L}:=\left\{ C_{p}\right\} _{p\in \partial L}$ of conjugacy classes
of order $k$ in $G$. Let us denote by
\[
\hom _{\partial L}\left( \pi _{1}\left( \Sigma \setminus \partial L\right)
,G\right)
\]
the set of morphisms $\pi _{1}\left( \Sigma \setminus \partial L\right)
\longrightarrow G$ sending a homotopy class of loops around $p\in \partial L$
into the conjugacy class $C_{p}$. Factoring out by conjugation leads to the
space
\begin{equation}
\mathcal{M}_{\left( \Sigma ,\partial L,\partial \mathcal{R}\right) }:=\hom
_{\partial L}\left( \pi _{1}\left( \Sigma \setminus \partial L\right)
,G\right) /G  \label{wilson class.sol}
\end{equation}
which is the analogue of $\mathcal{M}_{\Sigma }$ as defined in
$\left( \text{\ref{class.sol.}}\right) $. The quantization of
$\mathcal{M}_{\left( \Sigma ,\partial L,\partial
\mathcal{R}\right) }$ now proceeds in the same way since the
results of \cite{Narasimhan--Seshadri 1965} and \cite{d} can be applied
in this case too.

\begin{theorem}
\textit{The space }$\mathcal{M}_{\left( \Sigma ,\partial L,\partial \mathcal{%
R}\right) }$ \textit{is homeomorphic to a moduli space} $M_{\tau }^{\left(
k\right) }$ \textit{of holomorphic} $G$\textit{--bundles over the Riemann
surface} $\Sigma _{\tau }$ \textit{obtained by the choice of a complex
structure} $\tau $ \textit{on} $\Sigma $. \textit{On this space we have a
natural choice for a line bundle} $L_{k}$ \textit{whose holomorphic sections
give the quantization at level }$k$\textit{\ i.e. }
\begin{equation}
\mathcal{H}_{\tau }^{k}\left( \Sigma ,\partial L,\partial \mathcal{R}\right)
:=H^{0}\left( M_{\tau }^{\left( k\right) },L_{k}\right)
\label{Wilson-Hilbert}
\end{equation}
\end{theorem}

\noindent Note that the introduction of Wilson lines also makes the moduli
spaces $M_{\tau }^{\left( k\right) }$ dependent on the level $k$. As above
the finite dimensional complex vector space defined in (\ref{Wilson-Hilbert})%
\index{Hilbert space, of Chern--Simons QFT} varies holomorphically with
$\tau $ and gives rise to a projectively flat holomorphic vector bundle over
the moduli space of compact Riemann surfaces of fixed genus. Up to a scalar
factor we have obtained the desired relativized modular functor
\[
Z:\left( \Sigma ,\partial L,\partial \mathcal{R}\right) \longmapsto
Z_{\left( \Sigma ,\partial L,\partial \mathcal{R}\right) }:=\mathcal{H}%
^{k}\left( \Sigma ,\partial L,\partial \mathcal{R}\right)
\]
Note that an evaluation of the expectation value\index{Wilson
line, expectation value} $\left\langle
W_{L}^{\mathcal{R}}\right\rangle $ defined by applying $\left(
\text{\ref{exp.value}}\right) $and $\left(
\text{\ref{correlation}}\right) $ needs to fix once and for all the
undefined scalar factor. It can be realized by the choice of a
\textit{framing }(see definition \ref{framing}) for every knot
composing the link $L$: here we shall not enter into details about
by referring to \cite{Witten1989} and \cite{Atiyah1990b} for a
long their treatment. In the next section we will consider the problem for the
particular case in which $L$ is the unknotted knot.

\phantom{xx}
 \subsection{Computability and link invariants}~ \vskip 0.1in

\label{cli}

Let  $M$ be as in $\left( \text{\ref{chopping}}\right) $. By
$\left( \text{\ref{hamiltonian partition}}\right) $ and axiom 1 in
definition \ref {TQFT} we get
\begin{equation}
Z\left( M\right) =\left( \chi _{1},\chi _{2}\right)  \label{coupling}
\end{equation}
where $\chi _{1},\chi _{2}\in Z_{\Sigma }$. Similarly if we consider a
Wilson observable $W_{L}^{\mathcal{R}}$ on $M$ we get
\begin{equation}
Z\left( M;L,\mathcal{R}\right) =\left( \psi _{1},\psi _{2}\right)
\label{W-coupling}
\end{equation}
where $\psi _{1},\psi _{2}\in Z_{\left( \Sigma ,\partial L,\partial \mathcal{%
R}\right) }$.

\noindent \textit{These are the fundamental relations allowing the effective
computation of} $Z\left( M\right)$, $Z\left( M;L,\mathcal{R}\right) $ \textit{%
and }$\left\langle W_{L}^{\mathcal{R}}\right\rangle $, \textit{essentially by
connecting them with the link invariants of} $L$ \textit{in} $M$.

\noindent In the present section, following \cite{Witten1989}, we compute some of
those quantities when $M=S^{3}$ and $G=SU\left( N\right) $.

\begin{proposition}
\textit{Assume} $M=S^{3}$ \textit{and} $G=SU\left( N\right) $. \textit{Then
the expectation value }$\left\langle W_{L}^{\mathcal{R}}\right\rangle $%
\textit{\ of any Wilson observable can be inductively evaluated like a }%
Jones polynomial\index{Jones polynomial}\textit{\ }$V_{L}\left( q\right) $
\textit{in the variable }
\begin{equation}
q:=\exp \left( \frac{2\pi i}{N+k}\right)   \label{q}
\end{equation}
\textit{by applying the} skein relation $\left( \text{\ref{vev-skein}}%
\right) $ \textit{and the }mirror property $\left( \text{\ref{vev-mirror}}%
\right) $\textit{, when }$L$\textit{\ is considered in the }standard framing%
\index{standard framing, of a knot/link} \textit{and }$\mathcal{R}$\textit{\
is assigned by choosing the defining }$N$\textit{--dimensional
representation }$R$\textit{\ of }$SU\left( N\right) $\textit{\ for every
knot composing }$L$\textit{. In particular, if }$L$\textit{\ is the un-knot }$K$\textit{\ }
\begin{equation}
\left\langle W_{K}^{R}\right\rangle =\frac{q^{\frac{N}{2}}-q^{-\frac{N}{2}}}{%
q^{\frac{1}{2}}-q^{-\frac{1}{2}}}=\frac{\sin \left( \frac{N\pi }{N+k}\right)
}{\sin \left( \frac{\pi }{N+k}\right) }  \label{uvev}
\end{equation}
\textit{Moreover\ }
\begin{equation}
Z\left( S^{3}\right) =\left( k+N\right) ^{-N/2}\sqrt{\frac{k+N}{N}}%
\prod_{j=1}^{N}\left\{ 2\sin \left( \frac{j\pi }{k+N}\right) \right\} ^{N-j}
\label{S3partition}
\end{equation}
\textit{and}
\begin{eqnarray}
Z\left( S^{3};K,R\right)  &=&\frac{2}{\left( k+N\right) ^{N/2}}\sqrt{\frac{%
k+N}{N}}\sin ^{N-2}\left( \frac{\pi }{k+N}\right)   \label{unknot-partition}
\\
&&\sin \left( \frac{N\pi }{k+N}\right) \prod_{j=2}^{N-1}\left\{ 2\sin \left(
\frac{j\pi }{k+N}\right) \right\} ^{N-j}  \nonumber
\end{eqnarray}
\textit{\ }
\end{proposition}

\noindent Jones polynomials were firstly defined in \cite{Jones 1985} and
then generalized in \cite{Jones 1987} as a particular case of a
two--variable polynomial associated with a link by means of the
Ocneanu trace of a Hecke algebra representation of its braid
group. See also sections 1.3 and 1.4 in \cite{Atiyah1990a} and
section 2 in \cite{Labastida1999} for quick, but aimed at our
purpose, surveys on the argument.\index{Jones polynomial}

\begin{definition}
Denote by $L_{n}$ a link whose planar projection admits $n$ normal
crossings and by $L_{n+}$ and $L_{n-}$ those links admitting $n+1$ normal
crossings composed of the previous $n$ and by a further crossing which is an
\textit{over-crossing} or an \textit{under-crossing}, respectively. Given a
link $L\subset S^{3}$ the \textit{Jones polynomial} $V_{L}\left( q\right) $
is a Laurent polynomial in the variable $q^{\frac{1}{2}}$ inductively
defined by the\index{skein relation} \textit{skein relation}
\begin{equation}
\left( q^{\frac{1}{2}}-q^{-\frac{1}{2}}\right) V_{L_{n}}\left( q\right) -q^{%
\frac{N}{2}}V_{L_{n+}}\left( q\right) +q^{-\frac{N}{2}}V_{L_{n-}}\left(
q\right) =0  \label{Jones-skein}
\end{equation}
and the \textit{mirror property}
\begin{equation}
V_{L}\left( q\right) =V_{L^{\prime }}\left( q^{-1}\right)   \label{mirror}
\end{equation}
where $L^{\prime }$ is the mirror image%
\index{mirror image, of a knot/link}
\index{mirror property, of knot/link invariants} of the link $L$.
\end{definition}

\noindent To fix ideas start by considering the case in which $L$ is given
by two unlinked and unknotted circles $K_{1},K_{2}$ and $\Sigma $ is a
2--sphere $S^{2}$ which separates the two components of $L$ without cutting
any of them. Hence we get
\begin{eqnarray*}
&&Z_{\left( \Sigma ,\partial L,\partial \mathcal{R}\right) }=Z_{\Sigma
}=Z_{S^{2}} \\
&&\psi _{1}=Z\left( M_{1};K_{1},R_{1}\right) \\
&&\left( \quad ,\psi _{2}\right) =Z\left( M_{2};K_{2},R_{2}\right)
\end{eqnarray*}
Since $\dim _{\C}Z_{S^{2}}=1$, all the vectors $\chi _{1},\chi
_{2},\psi _{1},\psi _{2}$ are multiples of the same vector. By $\left( \text{%
\ref{coupling}}\right) $ and $\left( \text{\ref{W-coupling}}\right) $ this
gives
\begin{eqnarray*}
Z\left( M;L,\mathcal{R}\right) \cdot Z\left( M\right) &=&\left( \psi
_{1},\psi _{2}\right) \left( \chi _{1},\chi _{2}\right) \\
&=&\left( \psi _{1},\chi _{2}\right) \left( \chi _{1},\psi _{2}\right)
=Z\left( M;K_{1},R_{1}\right) \cdot Z\left( M;K_{2},R_{2}\right)
\end{eqnarray*}
whose quotient by $Z\left( M\right) ^{2}$ is
\begin{equation}
\left\langle W_{L}^{\mathcal{R}}\right\rangle =\left\langle
W_{K_{1}}^{R_{1}}\right\rangle \left\langle W_{K_{2}}^{R_{2}}\right\rangle
\label{multiplicativity}
\end{equation}
By iterating such a relation for an arbitrary collection of unlinked and
unknotted Wilson lines $L=\left\{ K_{i}\right\} _{i=1}^{r}$ we obtain that
\begin{equation}
\left\langle W_{L}^{\mathcal{R}}\right\rangle =\prod_{i=1}^{r}\left\langle
W_{K_{i}}^{R_{i}}\right\rangle  \label{vevs-multiplicativity}
\end{equation}
A first consequence of such a multiplicativity on expectation values of
unlinked and unknotted Wilson lines is that $\left\langle
W_{K}^{R}\right\rangle \neq 0$ for an unknotted Wilson line otherwise we
would have a Chern--Simons theory which does not distinguish a knot from a link!

\noindent Let us now consider four marked points $\left\{ p_{j}\right\}
_{j=1}^{4}$ on $\Sigma =S^{2}$. They may be obtained either as the
transversal section of the unlinked and unknotted link $L_{0}=\left\{
K_{1},K_{2}\right\} $ ( $S^{2}$ cuts two points on both $ K_{1}$ and $K_{2}$)
or as a section of the
two links $L_{+},L_{-}$ given by the two oriented knots whose planar normal
crossings projection gives a figure eight ($S^{2}$ cuts two points on both the
circles composing the figure eight): $L_{+}$ has an over--crossing while $L_{-}$ an under--crossing. If
we assume that the same representation $R$ of $G$ is associated with every
knot composing these links we may arrange the four points to give
\begin{eqnarray}
\left( \Sigma ,\partial L_{0},\partial \mathcal{R}_{0}\right)  &=&\left(
\Sigma ,\partial L_{+},\partial \mathcal{R}_{+}\right) =\left( \Sigma
,\partial L_{-},\partial \mathcal{R}_{-}\right)  \\
&=&\left( S^{2},\left\{ p_{j}\right\} _{j=1}^{4},\left\{ R,R,\overline{R},%
\overline{R}\right\} \right) =:\mathcal{H}  \label{boundary}
\end{eqnarray}
If we have the decomposition
\[
R\otimes R=\bigoplus_{h=1}^{s}E_{h}
\]
where $E_{h}$ is an irreducible representation of $G$, it turns out that
\begin{equation}
d:=\dim _{\C}\mathcal{H}\leq s  \label{dim_bound}
\end{equation}
and we get $d=s$ for large $k$ (see \cite{Witten1989}, section 3). In particular
if $G=SU\left( N\right) $ and $R$ is the defining $N$--dimensional
representation, then $s=2$ and
\begin{equation}
d=\left\{
\begin{array}{ll}
1 & \text{\quad if\quad }k=1 \\
2 & \quad \text{otherwise}
\end{array}
\right.   \label{SUdim_bound}
\end{equation}
For $i=1,2$ let us call $M_{i}^{0},M_{i}^{+},M_{i}^{-}$ the two
pieces cut by $S^{2}$ in the three different cases. Note that the
exterior pieces may be assumed to be
\begin{equation}
M_{1}^{0}=M_{1}^{+}=M_{1}^{-}=:M_{1}  \label{exterior}
\end{equation}
while the interior pieces $M_{2}^{0},M_{2}^{+},M_{2}^{-}$ may be thought to
be related by a diffeomorphism on the boundary exchanging two of the four
marked points. As in $\left( \text{\ref{W-coupling}}\right) $ the four
pieces $M_{1},M_{2}^{0},M_{2}^{+},M_{2}^{-}$ determine four vectors
\[
\psi _{1},\psi _{2}^{0},\psi _{2}^{+},\psi _{2}^{-}\in \mathcal{H}
\]
whose products evaluate the associated partition functions. Actually these
vectors are not known but the dimensional bound $\left( \text{\ref{dim_bound}%
}\right) $ may give rise to relations among them and their
products which results in being similar to the defining relations of
some link invariants. In particular when $G=SU\left( N\right) $
and all the knots are associated with the defining
$N$--dimensional representation, the dimensional bound $\left(
\text{\ref{SUdim_bound}}\right) $ allows one to conclude that $\psi
_{2}^{0},\psi _{2}^{+},\psi _{2}^{-}$ are linearly dependent and
so there must exist $\alpha ,\beta ,\gamma \in \C$ such that
\begin{equation}
\alpha \left( \psi _{1},\psi _{2}^{0}\right) +\beta \left( \psi _{1},\psi
_{2}^{+}\right) +\gamma \left( \psi _{1},\psi _{2}^{-}\right) =0
\label{vectorial skein}
\end{equation}
Hence the same relation can be established on the associated correlation
functions as follows:
\begin{equation}
\alpha Z\left( M;L_{0},\mathcal{R}_{0}\right) +\beta Z\left( M;L_{+},%
\mathcal{R}_{+}\right) +\gamma Z\left( M;L_{-},\mathcal{R}_{-}\right) =0
\label{skein0}
\end{equation}
It actually gives a recursive relation among links $L_{n},L_{n+}$ and $L_{n-}
$. In fact we can always cut these links by an $S^{2}$ leaving outside all
the first $n$ crossings: its interior then again gives $%
M_{2}^{0},M_{2}^{+},M_{2}^{-}$, respectively. Since $\alpha ,\beta ,\gamma $
depend only on the three vectors $\psi _{2}^{0},\psi _{2}^{+},\psi _{2}^{-}$%
, $\left( \text{\ref{vectorial skein}}\right) $ does not depend on $\psi _{1}$
and we again get
\begin{equation}
\alpha Z\left( M;L_{n},\mathcal{R}_{n}\right) +\beta Z\left( M;L_{n+},%
\mathcal{R}_{n+}\right) +\gamma Z\left( M;L_{n-},\mathcal{R}_{n-}\right) =0
\label{skein}
\end{equation}
We can then assume $\alpha \neq 0$, otherwise $\left( \text{\ref{skein}}%
\right) $ would imply that up to a scalar factor we can exchange an
over--crossing for an under--crossing i.e. every knot could be untied and our
Chern--Simons theory would not distinguish topologically non--equivalent
observables!

\noindent Since $M=S^{3}$ it is possible to continuously deform $L_{+}$ and $%
L_{-}$ to an oriented circle $K$ by applying a\index{Reidmeister move}
\textit{Reidemeister moving} i.e. a transformation induced on the planar
image with normal crossings of a knot in $S^{3}$ by a homeomorphism applied
to the original spatial knot (see \cite{Reidemeister1933}). By $\left( \text{\ref
{skein0}}\right) $ we can then write
\[
\alpha Z\left( M;\left\{ K_{1},K_{2}\right\} ,\left\{ R,R\right\} \right)
+\left( \beta +\gamma \right) Z\left( M;K,R\right) =0
\]
Divide by $Z\left( M\right) $ and recall $\left( \text{\ref{multiplicativity}%
}\right) $ to get
\[
\alpha \left\langle W_{K}^{R}\right\rangle \left\langle
W_{K}^{R}\right\rangle +\left( \beta +\gamma \right) \left\langle
W_{K}^{R}\right\rangle =0
\]
Since $\left\langle W_{K}^{R}\right\rangle \neq 0$ we obtain
\begin{equation}
\left\langle W_{K}^{R}\right\rangle =-\frac{\beta +\gamma }{\alpha }
\label{unknot-vev}
\end{equation}
Then by the knowledge of $\alpha ,\beta ,\gamma $, $\left( \text{\ref{skein}}%
\right) $ allows to inductively determine $\left\langle W_{L}^{\mathcal{R}%
}\right\rangle $ for every $L$ once we know a relation linking $\left\langle
W_{L}^{\mathcal{R}}\right\rangle $ and $\left\langle W_{L^{\prime }}^{%
\mathcal{R}^{\prime }}\right\rangle $ .

\noindent To determine $\alpha ,\beta ,\gamma $ let us concentrate
on the boundary diffeomorphisms relating
$M_{2}^{0},M_{2}^{+},M_{2}^{-}$. We can pass from $L_{+}$ to
$L_{0}$ by exchanging two of the four marked points on the
boundary $S^{2}$. Let us denote by
\[
f:M_{2}^{+}\longrightarrow M_{2}^{0}
\]
this ``half--monodromy'' diffeomorphism. Note that
\[
f\circ f:M_{2}^{+}\longrightarrow M_{2}^{-}
\]
since exchanging again the same two points we pass from $L_{0}$ to $L_{-}$.
By functoriality of TQFT we get an induced isomorphism $Z\left( f\right) \in
\ $ Aut $\left( \mathcal{H}\right) $ such that
\begin{equation}
\psi _{2}^{-}=Z\left( f\right) \psi _{2}^{0}=Z\left( f\right) ^{2}\psi
_{2}^{+}  \label{bdry-trasf}
\end{equation}
Since $Z\left( f\right) $ must satisfy its characteristic equation we get
the relation
\begin{equation}
\psi _{2}^{-}-\ \left( \text{tr\ }Z\left( f\right) \right) \psi
_{2}^{0}+\left( \det Z\left( f\right) \right) \psi _{2}^{+}=0
\label{char.eq.}
\end{equation}
which allows us to completely determine $\alpha ,\beta ,\gamma $ from the
knowledge of the eigenvalues of $Z\left( f\right) $. The latter are
calculated when $M=S^{3}$ in \cite{Moore-Sieberg 1988}. By comparing $\left( \text{%
\ref{vectorial skein}}\right) $ and $\left( \text{\ref{char.eq.}}\right) $
and setting $q$ as in $\left( \text{\ref{q}}\right) $ we can rewrite $%
\left( \text{\ref{skein}}\right) $ for $M=S^{3}$ as follows:
\begin{eqnarray}
\left( q^{\frac{1}{2}}-q^{-\frac{1}{2}}\right) Z\left( M;L_{n},\mathcal{R}%
_{n}\right) -q^{\frac{1}{2N}}Z\left( M;L_{n+},\mathcal{R}_{n+}\right)  &&
\label{TQFTskein} \\
+q^{-\frac{1}{2N}}Z\left( M;L_{n-},\mathcal{R}_{n-}\right)  &=&0  \nonumber
\end{eqnarray}
Hence by $\left( \text{\ref{unknot-vev}}\right) $ the expectation value for
the unknotted Wilson line\index{Wilson line, expectation value} is given by
\begin{equation}
\left\langle W_{K}^{R}\right\rangle =\frac{q^{\frac{1}{2N}}-q^{-\frac{1}{2N}}%
}{q^{\frac{1}{2}}-q^{-\frac{1}{2}}}  \label{TQFTuvev}
\end{equation}
This value does not coincide with $\left( \text{\ref{uvev}}\right) $ since the
relation $\left( \text{\ref{TQFTskein}}\right) $ is similar but not equal to
the skein relation $\left( \text{\ref{Jones-skein}}\right) $. The reason from
such a discrepancy must be found in the implicit \textit{framing }choice%
\textit{\ }we used to write $\left( \text{\ref{vectorial skein}}\right) $,
which is not the same as the \textit{standard framing} used in knot theory.

\begin{definition}
\label{framing}A \textit{framing} of a knot%
\index{framing, of a knot}
\index{framed knot}
$K$ is a closed curve $K_{f}$ obtained as a small deformation of $K$ along a
normal vector field direction. The pair $\left( K,K_{f}\right) $ is called a
\textit{framed knot.}
\end{definition}

\noindent At the end of subsection \ref{HAMILTON} we noted that the
evaluation of a Wilson observable expectation value $\left\langle W_{L}^{%
\mathcal{R}}\right\rangle $ needs to fix once and for all the undefined scalar
factors which occur in the projective definition of the Hamiltonian
quantities via TQFT. Actually by making assumptions $\left( \text{\ref
{boundary}}\right) $ and $\left( \text{\ref{exterior}}\right) $ we did a
particular choice of those scalar factors which does not coincide with the
canonical choice usually adopted for knots in $S^{3}$ by requiring that the
\textit{Gauss self--linking number} \index{self--linking, Gauss number}is
trivial for every knot\textit{\ }(see \cite{Witten1989} section 2.1 for the
definition; see also \cite{mv2001} section 3 for a recent discussion
of the problem in connection with the concept of a framed knot): this is
what is usually meant by the \textit{standard framing} \index{standard framing,
of a knot/link}of a knot.

\noindent Note that the coefficient associated with the unknotted unlinked $%
L_{0}$ is $q^{1/2}-q^{-1/2}$ both in $\left( \text{\ref{TQFTskein}}\right) $
and in $\left( \text{\ref{Jones-skein}}\right) $. Since by $\left( \text{\ref
{bdry-trasf}}\right) $ we pass from $\psi _{2}^{0}$ to $\psi _{2}^{-}$ by
applying $Z\left( f\right) $ while its inverse $Z\left( f\right) ^{-1}$
allows us to pass to $\psi _{2}^{+}$ we can argue that
\[
q^{-\frac{N}{2}}q^{\frac{1}{2N}}=\left( q^{\frac{N}{2}}q^{-\frac{1}{2N}%
}\right) ^{-1}=\exp \left( \pi i\frac{\left( 1-N^{2}\right) }{N\left(
N+k\right) }\right)
\]
\textit{is the factor expressing the framing change through the
half--monodromy} $f$. It follows that, by adopting the standard framing, the
expectation value $\left( \text{\ref{TQFTuvev}}\right) $ of the unknotted
Wilson line must be rewritten as in $\left( \text{\ref{uvev}}\right) $.
Although the skein relations $\left( \text{\ref{TQFTskein}}\right) $ and $%
\left( \text{\ref{Jones-skein}}\right) $ are not the same, the
``polynomials'' defined by the former also satisfy the mirror property
\begin{equation}
\left\langle W_{L}^{\mathcal{R}}\right\rangle \left( q\right) =\left\langle
W_{L^{\prime }}^{\mathcal{R}^{\prime }}\right\rangle \left( q^{-1}\right)
\label{vev-mirror}
\end{equation}
We can then conclude that the skein relation
\begin{equation}
\left( q^{\frac{1}{2}}-q^{-\frac{1}{2}}\right) \left\langle W_{L_{n}}^{%
\mathcal{R}_{n}}\right\rangle -q^{\frac{N}{2}}\left\langle W_{L_{n+}}^{%
\mathcal{R}_{n+}}\right\rangle +q^{-\frac{N}{2}}\left\langle W_{L_{n-}}^{%
\mathcal{R}_{n-}}\right\rangle =0  \label{vev-skein}
\end{equation}
and the mirror property $\left( \text{\ref{vev-mirror}}\right) $ allow us to
inductively express in the standard framing the expectation value $%
\left\langle W_{L}^{\mathcal{R}}\right\rangle $ of any Wilson observable in $%
S^{3}$, when $G=SU\left( N\right) $ and all the representations associated
with knots are the defining $N$--dimensional ones.

\noindent Note that when we fix $N=2$ the unique variable is the level $k$
of the theory while when $N$ is general $\left\langle W_{L}^{\mathcal{R}%
}\right\rangle $ can be interpreted also like a\index{HOMFLY
polynomial} \textit{HOMFLY polynomial} (see \cite{fy} for the
definition of this two--variable polynomial invariant of links).

\noindent The skein relation $\left( \text{\ref{vev-skein}}\right) $ cannot
evaluate the partition function $Z\left( S^{3}\right) $ and consequently the
correlation function of any Wilson observable. Their evaluation follows by
generalizing the previous procedure to every three--manifold $M$.

\begin{definition}
Let $K\subset S^{3}$ be an unknotted circle and $T$ a tubular neighborhood
of $K$, i.e. a solid torus centered in $K.$ Then
\[
S^{3}=\left( S^{3}\setminus T\right) \cup _{\Sigma }T
\]
where $\Sigma :=\partial T$ is a two--dimensional torus. If before the
gluing we apply a diffeomorphism on the boundary $\partial T$ then the
gluing will give us a new three-manifold $M$ which is said to be obtained by
$S^{3}$ \index{surgery, on a knot}\textit{after a surgery on the knot} $K$.
\end{definition}

\begin{proposition}
\label{surgery}\textit{Any three--manifold} $M$ \textit{can be obtained by} $%
S^{3}$ \textit{up to a finite number of} surgeries on knots. \textit{Hence
the partition functions and expectation values on a general} $M$ \textit{can
be evaluated by those on} $S^{3}$ \textit{once it is known how the
repeated surgeries act on these quantities and on the knot framings}.
\end{proposition}

\noindent An important application of this proposition is given by the
manifold
\[
M:=S^{2}\times S^{1}
\]
If we think of $S^{3}$ as the compactification by a point of
$\R^{3}$ and of $K$ as the unit circle in the plane $z=0$,
consider the following surgery on $K$. Let $\Sigma $ be a
two--dimensional torus around $K$
invariant under an inversion of $\R^{3}$: the tubular neighborhood of $%
K $ is the interior $T_{1}$ of $\Sigma $. Note that the exterior $%
T_{2}=S^{3}\setminus T_{1}$ is a solid torus too and we get
\begin{equation}
S^{3}=T_{1}\cup _{\Sigma }T_{2}  \label{S3}
\end{equation}
On the other hand if $T_{1},T_{2}$ are thought of as two solid tori
which can be identified by a translation of $\R^{3}$ we get
\begin{equation}
S^{2}\times S^{1}=T_{1}\cup _{\Sigma }T_{2}  \label{S2xS1}
\end{equation}
since $T_{i}=D_{i}\times S^{1},\Sigma =S^{1}\times S^{1}$ and $%
S^{2}=D_{1}\cup _{S^{1}}D_{2}$. $\left( \text{\ref{S3}}\right) $ and $\left(
\text{\ref{S2xS1}}\right) $ differ simply by the diffeomorphism applied on
the boundary $\Sigma $ to glue the solid tori $T_{i}$: in the former it is
given by an inversion while in the latter by a translation.

\noindent This example is important because $Z\left( S^{2}\times S^{1};L,%
\mathcal{R}\right) $ can be obtained by the TQFT axioms easier than $Z\left(
S^{3};L,\mathcal{R}\right) $. Then we get a method to evaluate our partition
functions on $S^{3}$, which is the main ingredient of  Witten's proof of a
conjecture of Verlinde (see \cite{Verlinde1988}) already proved in
\cite{Moore-Sieberg 1988}.
In \cite{Verlinde1988} it is shown how to canonically get a basis $\left\{
v_{0},\ldots ,v_{t-1}\right\} $ of $Z_{\Sigma }$ \index{Verlinde basis}after
the choice of a homology basis $\left\{ \gamma _{1},\gamma _{2}\right\} $
for $H_{1}\left( \Sigma ,\Bbb{Z}\right) $: calling $T$ the interior of $%
\Sigma $ the first basis vector $v_{0}$ is chosen to give $Z\left( T\right)
\in Z_{\Sigma }$. The two solid tori $T_{1},T_{2}$ giving $S^{2}\times S^{1}$
in $\left( \text{\ref{S2xS1}}\right) $ are two identical copies of $T$
identified by a translation. This gives
\begin{equation}
v_{0}=Z\left( T_{2}\right) \quad ,\quad \left( v_{0},\quad \right) =Z\left(
T_{1}\right) \quad ,\quad \left( v_{0},v_{0}\right) =Z\left( S^{2}\times
S^{1}\right)   \label{S2xS1partition}
\end{equation}
On the other hand if we think of $\Sigma $ as in $\left( \text{\ref{S3}}%
\right) $ the inversion of $\R^{3}$ acts on $H_{1}\left( \Sigma ,\Bbb{Z}%
\right) $ by sending
\begin{eqnarray}
\gamma _{1} &\longmapsto &-\gamma _{1}  \label{homolgy-inv} \\
\gamma _{2} &\longmapsto &\gamma _{2}  \nonumber
\end{eqnarray}
Let $\tau =a+ib$ be the complex number in the\index{Siegel, upper
half--plane} Siegel upper half--plane
\[
\Bbb{H}:=\left\{ \tau \in \C:\ \text{Im\ }\left( \tau \right)
>0\right\}
\]
representing the isomorphism class of the complex torus $\Sigma $. The
transformation induced on $\Bbb{H}$ by the inversion acts as follows:
\[
\tau =a+ib\longmapsto \frac{1}{\left| \tau \right| ^{2}}\left( -a+ib\right)
=-\tau ^{-1}
\]
It is the modular transformation represented by the element
\[
S=\left(
\begin{array}{ll}
0 & -1 \\
1 & 0
\end{array}
\right)
\]
in the modular group\index{modular group}
\[
\Gamma :=SL\left( 2,\Bbb{Z}\right) /\left\{ \pm I\right\}
\]
where $I=\left(
\begin{array}{ll}
1 & 0 \\
0 & 1
\end{array}
\right) $. $\Gamma $ acts on $\Bbb{H}$ by setting
\[
\left(
\begin{array}{ll}
a & b \\
c & d
\end{array}
\right) \tau =\left( a\tau +b\right) \left( c\tau +d\right) ^{-1}
\]
Since the isomorphism classes of complex tori are parametrized by the\textit{%
\ modular curve} $\Gamma \setminus \Bbb{H}$ \index{modular curve}it turns
out that the inversion realizes a diffeomorphism of $\Sigma $ which
preserves the complex structure (see the first chapter in \cite{Silverman} for
further details about and a careful construction of the quotient $\Gamma
\setminus \Bbb{H}$). It induces an isomorphism on $Z_{\Sigma }$ which can be
represented on the Verlinde basis by a complex $t\times t$ matrix $S_{i}^{j}$
such that
\[
v_{i}=\sum_{j}S_{i}^{j}v_{j}
\]
Therefore by $\left( \text{\ref{S3}}\right) $ and $\left( \text{\ref
{S2xS1partition}}\right) $ we get
\[
Z\left( S^{3}\right) =\left( v_{0},\sum_{j}S_{0}^{j}v_{j}\right)
=\sum_{j}S_{0}^{j}\left( v_{0},v_{j}\right)
\]
This formula gives an effective evaluation of $Z\left( S^{3}\right) $ since
the numbers $g_{ij}:=\left( v_{i},v_{j}\right) $ and the matrix $S_{i}^{j}$
are given by the knowledge of the Verlinde basis of $Z_{\Sigma }$. Hence by
setting $S_{i,j}:=\sum_{k}S_{i}^{k}g_{jk}$ we get
\[
Z\left( S^{3}\right) =S_{0,0}
\]
When $G=SU\left( N\right) $ we obtain the following result:
\begin{equation}
S_{0,0}=\left( k+N\right) ^{-N/2}\sqrt{\frac{k+N}{N}}\prod_{j=1}^{N}\left\{
2\sin \left( \frac{j\pi }{k+N}\right) \right\} ^{N-j}  \label{S00}
\end{equation}
allowing us to conclude $\left( \text{\ref{S3partition}}\right) $.By recalling $%
\left( \text{\ref{uvev}}\right) $ we are able to write $Z\left(
S^{3};K,R\right) $ as in $\left( \text{\ref{unknot-partition}}\right) $
for the unknotted knot $K$ in the defining $N$--dimensional representation $R
$ of $SU\left( N\right) $.

\section{\label{GVC} The Gopakumar--Vafa conjecture}

This section discusses the conjecture itself, its origin, and
its relation to geometric transitions.
 We also present supporting evidence, which leads to the
 uncharted territory of ``open Gromov-Witten invariants".

We start with the original observation of Gopakumar and Vafa
(by comparing the partition functions) and show in this first part how Witten's interpretation of the Chern--Simons theory as
an open string theory (see \cite{Witten1992}) provides the tools for the geometric interpretation
of the duality.

 There is no discussion of II--$A$ theory itself, partly because of
 time constraints, partly because II--$A and $II--$B$ theories
 and (closed) Gromov-Witten invariants have
recently been in the spot light, thanks to the celebrated ``mirror
 symmetry" and its enumerative predictions (see for
 example \cite{CoxKatz1999}).

\begin{conjecture}(Gopakumar--Vafa)\cite{gva}
\label{GVconjecture}%
\index{Gopakumar--Vafa conjecture}
\index{conjecture, of Gopakumar and Vafa} (Notation as in \ref{conifold
topology}): \textit{The} $SU\left( N\right) $ \textit{Chern--Simons
theory on }$S^{3}\subset \widehat{Y}:=T^{*}S^{3}$ \textit{of level }$k$%
\textit{\ is equivalent,  for large }$N$,\textit{ to a type II--$A$ closed string theory (with
fluxes) on the local
Calabi-Yau }$Y:=\mathcal{O}_{\P^{1}}\left( -1\right) \oplus \mathcal{O}%
_{\P^{1}}\left( -1\right) .$
\end{conjecture}
\noindent (The language used here  reflects the reformulation of the conjecture
given in \cite{ov} rather than the original one.)

\begin{theorem}
\label{Witten92}\cite{Witten1992}\textit{\ Let }
$
\widehat{Y}=T^{*}L $
\textit{\ be a local Calabi--Yau threefold. Then there exist topological
string theories with }$\widehat{Y}$\textit{\ as target space, such that
their open sectors are exactly equivalent to a QFT on }$L$.
\end{theorem}

\begin{conjecture}
\label{GVWconjecture}(Gopakumar--Vafa after Witten)%
\index{Gopakumar--Vafa conjecture} \index{conjecture, of Gopakumar
and Vafa}\textit{\ A topological open string theory of type
II--}$A$\textit{\ on }$\widehat{Y}:=T^{*}S^{3}$\textit{\ with
}$N$\textit{\ D6--branes wrapped around the base }$S^{3}$\textit{\
is equivalent, for large }$N$\textit{, to a type II--}$A$\textit{\
closed string theory on the local Calabi-Yau
}$Y:=\mathcal{O}_{\P^{1}}\left(
-1\right) \oplus \mathcal{O}_{\P^{1}}\left( -1\right) $\textit{\ with }$%
N$\textit{\ units of }$2$\textit{--form Ramond--Ramond flux}\index{RR, flux}%
\textit{\ through the exceptional }$S^{2}$\textit{.}

\noindent \textit{The transition}\index{transition}\textit{\ from }$%
Y$\textit{\ to }$\widehat{Y}$\textit{\ realizes the geometrical model of a
physical closed/open duality among string theories of type II--$A$. That is, the transition from }$Y$\textit{\ to }$\widehat{Y}$\textit{%
\ realizes the geometrical model of a physical duality relating a particular
type II--$A$ closed string theory on }$Y$\textit{\ and the }$SU\left( N\right)
$\textit{ Chern--Simons\ QFT on the Lagrangian submanifold }$S^{3}$\textit{%
\ of }$\widehat{Y}$\textit{\ for large }$N$\textit{.}

\end{conjecture}

\noindent This formulation of  conjecture \ref{GVconjecture} is
already given in \cite{gvc}; see also \cite{ov}. See
\cite{Vafa2001} for the correspondence between D6--branes and units
of RR flux.

Witten's work is more general: he proposes a string theory interpretation of the
Chern--Simons $U\left( N\right) $ gauge theory on a real
three--dimensional Lagrangian submanifold $L$ of a complex
Calabi-Yau threefold $\widehat{Y}$ and also extends beyond the hypothesis
\begin{equation}\label{ctg.hyp.}
\widehat{Y}=T^{*}L.
\end{equation} We refer to Appendix \ref{Witten} for more
details.

\vskip 0.2in
\noindent {\it Sketch of the proof: how Theorem \ref{Witten92} implies \ref{GVconjecture} $\leftrightarrow$ \ref{GVWconjecture}.}

Witten constructs  an``$A$--twisted sigma model''on $\widehat{Y}$.
In particular he consider maps $\phi $ from a Riemann surface $\Sigma $ with
boundary $\partial \Sigma $, to the target space $Y$, (i.e. $\phi $ is a
bosonic field of the open sector of this $A$--model) satisfying some
conditions. The most important assumption is that
\begin{equation}
\phi \left( \partial \Sigma \right) \subset L,  \label{bdr.cond.}
\end{equation}
There are also boundary conditions, involving derivatives of $\phi $
along the components of $\partial \Sigma $ and the fermionic fields. These
conditions are needed to preserve the fermionic symmetry but they do not
enter directly in the geometric picture (see section 3.1 in \cite{Witten1992} for
more details). If $Y=T^{*}L$ the weak coupling limit%
\index{limit, low energy (weak coupling)}
of the abstract string Lagrangian reduces exactly to the Lagrangian of a QFT
on $L$, that is, \textit{\ ``there are neither perturbative corrections nor
instanton corrections''} (see definition \ref{instanton})\textit{.} \textit{%
In the }$A$\textit{--twisted case such a limit turns out to be exactly a
Chern--Simons }$U\left( N\right) $\textit{ gauge theory}.

Gopakumar and Vafa observed that the above boundary conditions may be
expressed in terms of D--branes (see A. Lerda's lectures in the same volume)
by saying that Witten's open string theory is an\index{sigma model,
$A$--twisted} $A$\textit{--model topological open string theory with }$N$%
\textit{\ topological D}$6$\textit{--branes wrapped on }$L. \ \ \diamondsuit$

\phantom{xx}

\subsection{Matching of the free energies}~ \vskip
0.1in

In the next two subsections, we review the evidence for the conjectures
\ref{GVconjecture} and \ref{GVWconjecture}. The first evidence is given by
the \textit{matching of the ``free energies''} (or equivalently\textit{\
partition functions) }for the theories involved in the conjecture. The
second one is given by \textit{\ comparisons of the expectation values of
observables }in the two theories.

\vskip 0.2in

\begin{theorem}
\textit{\ The genus }$g$\textit{\ contribution to the perturbative
expansion of the free energy }$\left(
\text{\ref{top.exp.}}\right) $\textit{\ of the Chern--Simons theory on }$%
S^{3}$ \textit{coincides with the genus }$g$\textit{\ contribution
to the free energy of the closed string theory on
}$\mathcal{O}_{\P^{1}}\left( -1\right) \oplus $
$\mathcal{O}_{\P^{1}}\left( -1\right) $.
\end{theorem}

\noindent $\bullet$ {\bf The Chern--Simons side.}

\begin{definition} Let $Z\left(
S^{3}\right) $ be the partition function given by $\left( \text{\ref
{S3partition}}\right) $. Set
\begin{equation}
F\left( S^{3}\right) =-\log Z\left( S^{3}\right).   \label{free
energy}
\end{equation}
\end{definition}

\begin{proposition}
(\cite{'t Hooft 1974}, \cite{Periwal 1993}) \textit{For large }$N,$\textit{\ the free
energy }$\left( \text{\ref{free energy}}\right) $\textit{\ of a }$SU\left(
N\right) $\textit{ gauge Chern--Simons QFT on }$S^{3}$\textit{\ can be
expanded as follows }
\begin{equation}
F\left( S^{3}\right) =\sum_{g \geq 0}\mathcal{F}_{g}\left( \tau \right) N^{2-2g}.  \label{top.exp.}
\end{equation}

\textit{Here  }%
\begin{equation}
\lambda :=\frac{2\pi }{k+N}  \label{lambda}
\end{equation}
\textit{is the Chern--Simons coupling constant,  } $\tau :=\lambda
N$ \textit{ the 't Hooft coupling constant.\index{'t Hooft,
coupling constant} \textit{ The weak--coupling limit }$\lambda
\rightarrow 0$, $N\rightarrow +\infty $ \textit{leaves constant the
}'t Hooft coupling constant.}
\end{proposition}

\noindent \textit{ Sketch of the proof:} The statement follows by observing that, in the ``double line
notation'', Feynman diagrams contributing to the free energy $F$ may be
thought of as a sort of ``triangulation'' of a compact, connected
topological surface given by an admissible subdivision of the topological
surface into polygons and disks. The latter occur as the internal planar
regions of loops in Feynman diagrams: they have to be understood like
polygons admitting two edges and two vertices. 't Hooft observed that the
contribution due to a Feynman diagram is proportional to $\lambda
^{e-v}N^{h-l}$ where $l$ is the number of diagram loops (\textit{quark loops
}in 't Hooft's notation) and $e,v,h$ the number of edges, vertices and faces
respectively, in the induced ``triangulation''. Since a diagram loop
increases $h$ by 1 and $e,v$ by 2, the contribution due to a Feynman diagram
without loops and admitting $h^{\prime }=h-l$ faces is proportional to $%
\lambda ^{e-v}N^{h-l}$ as well. The Euler characteristic formula
\[
2-2g=h-e+v
\]
allows one to conclude that the Feynman diagrams' contributions to the free
energy $F$ can be labeled by the genus $g$ of the topological surface and
the number of faces $h$ of the induced ``triangulation''. The associated
contribution is then proportional to $\lambda ^{2g-2+h}N^{h}$ to get
\[
F=\sum_{g}\left( \sum_{h}C_{g,h}\lambda ^{2g-2+h}N^{h}\right)
\]
where $C_{g,h}$ are suitable coefficients computed by Periwal. If we now
consider the \index{limit, low energy (weak coupling)}
weak--coupling limit $\lambda \rightarrow 0$, $N\rightarrow +\infty $ leaving%
$\ \tau =\lambda N$ constant, then the free energy expansion can be
reorganized as follows
\[
F=\sum_{g}\left( \sum_{h}C_{g,h}\tau ^{2g-2+h}\right) N^{2-2g}=\sum_{g}%
\mathcal{F}_{g}\left( \tau \right) N^{\chi \left( g\right) }. \  \diamondsuit
\]

\begin{lemma}\label{ah}
\textit{Let}
\[
Z(S^{3})=\left( k+N\right) ^{-N/2}\sqrt{\frac{k+N}{N}}\prod_{j=1}^{N}\left\{
2\sin \left( \frac{j\pi }{k+N}\right) \right\} ^{N-j}
\]
\textit{be the Chern--Simons partition function, as  in }$\left(
\text{\ref {S3partition}}\right). $

\textit{ Set }$F(S^{3})=-logZ(S^{3})$ \textit{and}
\[
t=\frac{2\pi iN}{k+N},\ \ \ \lambda =\frac{2\pi }{k+N}
\]
\textit{as in }$\left( \text{\ref{lambda}}\right)
 $\textit{. The 't Hooft topological expansion  for large%
} $N$ $\left( \text{of equation \ref{top.exp.}}\right) $\textit{\ becomes, for small} $%
\lambda $%
\begin{equation}
F\left( \lambda ,t\right) =\sum_{g=0}^{+\infty }F_{g}\left( t\right) \lambda
^{-\chi \left( g\right) }  \label{top.exp.2}
\end{equation}
\textit{where} $F_{g}\left( t\right) =\tau ^{\chi \left( g\right) }\mathcal{F%
}_{g}\left( \tau \right) =\left( -1\right) ^{g+1}t^{\chi \left( g\right) }%
\mathcal{F}_{g}\left( -it\right) $\textit{. In particular: }
\begin{eqnarray}
F_{0}\left( t\right)  &=&\frac{i\pi ^{2}}{6}t-i\left( m+\frac{1}{4}\right)
\pi t^{2}+\frac{i}{12}t^{3}-\sum_{d=1}^{+\infty }d^{-3}\left(
1-e^{-dt}\right)   \nonumber \\
F_{1}\left( t\right)  &=&\frac{1}{24}t+\frac{1}{12}\log \left(
1-e^{-t}\right)   \label{CScoeff} \\
F_{g}\left( t\right)  &=&\frac{\left( -1\right) ^{g}B_{2g}}{2g\left(
2g-2\right) !}\left( \frac{B_{2g-2}}{\left( 2g-2\right) }%
+\sum_{d=1}^{+\infty }d^{2g-3}e^{-dt}\right) \quad \forall g\geq 2,
\nonumber
\end{eqnarray}
\textit{where} $m$ \textit{is an arbitrary integer coming from the
polydromic behavior of the complex logarithm and} $B_{h}$ \textit{is the} $%
h^{\text{th}}$ \textit{Bernoulli number defined by }
\[
\frac{x}{e^{x}-1}=\sum_{h=0}^{+\infty }B_{h}\frac{x^{h}}{h!}.
\]
\end{lemma}

\noindent  Note that in the physics literature, the
$2g^{\text{th}}$ Bernoulli number is often denoted by $B_{g}$
instead of $B_{2g}$.

The explicit computation of the expansion coefficients can be
performed either starting from $\mathcal{F}_{g}\left( \tau \right)
$ as in  \cite{Periwal 1993}
 (expansion for large $N$) or from $F_{g}\left( t\right) $ by following
Gopakumar and Vafa \cite{gvc} and \cite{gva}, \cite{gvb} (expansion for small $%
\lambda $). The key ingredient in expanding $F\left( S^{3}\right) $ is to
employ the Mittag--Leffler expansion for the logarithmic derivative of the
complex function $\sin \left( z\right) /z$. When $z=j\lambda /2$ we get the
following relation:
\[
\sin \left( \frac{j}{2}\lambda \right) =\frac{j}{2}\lambda
\prod_{d=1}^{+\infty }\left( 1-\frac{j^{2}\lambda ^{2}}{4\pi ^{2}d^{2}.}%
\right)
\]
Substituting in $\left( \text{\ref{S00}}\right) $ we get $\left( \text{%
\ref{top.exp.2}}\right) $. \vskip 0.2in

\noindent $\bullet$ {\bf The II--$A$ theory side.}

\vskip 0.1in
 We do not derive the perturbative expansion of the II--$A$
 theory; rather we take \ref{instanton} and \ref{gwyd} as its
 definition. See for example \cite{CoxKatz1999} for a discussion
 of these topics.
\begin{definition}
\label{instanton}Given a topological string theory whose target space is a
complex manifold $Y$, a \textit{world sheet instanton }(or simply \textit{%
instanton}) \textit{of genus }$g$ is a holomorphic map%
\index{instanton}
\index{instanton, of genus  $g$}
\index{instanton, correction}
\index{instanton, open}
\index{world sheet, of a instanton}
\[
\phi :\Sigma \longrightarrow Y
\]
from a Riemann surface of genus $g$. If the boundary $\partial
\Sigma $ is not empty $\phi $ is said to be \textit{open}, since a
similar instanton is typical of an open string.
\end{definition}
In our case  $Y=\mathcal{O}_{\P^{1}}\left( -1\right) \oplus
\mathcal{O}_{\P^{1}}\left( -1\right),$  and the only non-trivial
homology class is the exceptional $\Bbb P ^1$. The ``string
amplitude" ``counts" instantons with image the exceptional $\Bbb P
^1$:
\begin{definitionproposition}\label{gwyd} Let
\[
F^{\left( s\right) }\left( \lambda ,t\right) =\sum_{g=0}^{+\infty
}F_{g}^{\left( s\right) }\left( t\right) \lambda ^{-\chi \left(
g\right) }
\]be the perturbative expansion of the free energy (better:
``string amplitude'') of the type II--$A$(closed) string theory.
Then:
\begin{enumerate}
\item $\lambda$ is the string coupling constant and
$t$ \textit{ is interpreted
 as the K\"{a}hler modulus
of the exceptional locus }$S^{2}\cong \P^{1}$\textit{in} $Y$ (see
\ref{km} below).
\item  $F_{g}^{\left( s\right) }\left( t\right) $ is the contribution to
the string amplitude $F^{\left( s\right) }\left( \lambda ,t\right)
$ given by all the genus $g$ instantons and is called \textit{the
genus }$g$\textit{\ instanton correction.}
\item  $F^{(s)}_{g}\left( t\right) $ determines the genus $g$
Gromov--Witten invariants of $Y$, associated with maps of Riemann
surfaces with image the homology class of the exceptional locus
$\P^{1}\cong S^{2}\subset Y.$

\end{enumerate}
\end{definitionproposition}
\begin{definition}\label{km}
Let $\left( X,g\right) $ be a K\"{a}hler manifold; fix a closed
2--form $B$ on $X$ and denote by $J\in H^{2}\left( X,\R\right) $
the K\"{a}hler class of the Hermitian metric $g$. The cohomology
class of the form $\omega
=B+iJ$ is called the%
\index{K\"{a}lher, complexified class} \index{K\"{a}lher modulus}
\textit{complexified K\"{a}hler class associated
with }$g$. The \textit{K\"{a}hler modulus }of a given real 2--cycle $%
Z\subset X$ is defined by the period
\[
\int_{Z}\omega \in \C
\]
of the complexified K\"{a}hler class on it. By Stokes' theorem and
K\"{a}hler condition on $J$ it is well defined for the entire
homology class of $Z$.
\end{definition}

\begin{theorem}\label{gwy}
\textit{Let}
\[
F^{\left( s\right) }\left( \lambda ,t\right) =\sum_{g=0}^{+\infty
}F_{g}^{\left( s\right) }\left( t\right) \lambda ^{-\chi \left(
g\right) }
\]
 \textit{be the string amplitude of the type II--A
theory, as  in } $\left( \text{\ref{gwyd}}\right)$. \textit{ Then:
}
\begin{eqnarray}
F_{0}^{\left( s\right) }\left( t\right) &=&\frac{i\pi ^{2}}{6}t-ia\pi t^{2}+%
\frac{i}{12}t^{3}-\sum_{d=1}^{+\infty }d^{-3}\left(
1-e^{-dt}\right)   \nonumber \\
F_{1}^{\left( s\right)}\left( t\right)
&=& \frac{1}{24}t+\frac{1}{12}\log \left( 1-e^{-t}\right)   \\
F_{g}^{\left( s\right)}\left( t\right) &=&
\frac{|B_{2g}|}{2g\left( 2g-2\right) !}\left( \left( -1\right)
^{g+1}\frac{|B_{2g-2}|}{\left( 2g-2\right) }%
-\sum_{d=1}^{+\infty }d^{2g-3}e^{-dt}\right) \quad \forall g\geq
2. \nonumber
\end{eqnarray}

\end{theorem}

\noindent \textit{Sketch of the proof of \ref{gwy}:} $%
F_{0}^{\left( s\right) }\left( t\right) $ can be found in
\cite{CdOGP}. The coefficient $a$ does not have a direct
topological interpretation on $Y$.

\noindent The computation of $F_{1}^{\left( s\right) }\left( t\right) $ and $%
F_{2}^{\left( s\right) }\left( t\right) $ can be found in
\cite{bcov} and in \cite{bcov2} respectively: in our situation
they match exactly $F_{1}\left( t\right) $ and $F_{2}\left(
t\right) $.

\noindent In \cite{fp} Faber and Pandharipande compute
$F_{g}^{\left( s\right) }\left( t\right) $ for every genus $g\geq
2$. In particular for $g\geq 2$ we can write
\begin{equation}
F_{g}^{\left( s\right) }\left( t\right) =-\left\langle
1\right\rangle _{g,0}^{Y}-\sum_{d=1}^{+\infty }C\left( g,d\right)
e^{-dt} \label{g--instantons}
\end{equation}
where $\left\langle 1\right\rangle _{g,0}^{Y}$ is the genus $g$,
degree 0 Gromov--Witten invariant of our Calabi--Yau $Y$ giving
the instanton correction due to constant maps. On the other hand
the series on the right
gives, for every $d$, the instanton correction due to maps realizing a $d$%
--covering with genus $g$ of the exceptional $\P^{1}$. Theorem 3
in \cite{fp} gives
\[
C\left( g,d\right) =\left| \chi \left( \overline{M}_{g}\right) \right| \frac{%
d^{2g-3}}{\left( 2g-3\right) !}
\]
where $\chi \left( \overline{M}_{g}\right) $ is the orbifold Euler
characteristic of the coarse moduli space $\overline{M}_{g}$.
$\overline{M}_{g}$ denotes the compactified moduli space of
projective,
connected, nodal, Deligne--Mumford stable curves of arithmetic genus $%
g$; if $g \geq 2$,  $\overline{M}_{g}$ is an irreducible variety
of dimension $3g-3$. $\overline{M}_{g}$ has orbifold singularities
if regarded as an ordinary coarse moduli space, it is smooth if
regarded as a Deligne--Mumford stack: see \cite{fup} and chapter 7
in
\cite{CoxKatz1999} for general reference.\\

$\chi \left( \overline{M}_{g}\right) $ can be
expressed in terms of Bernoulli numbers by means of the following %
\index{Harer--Zagier formula}Harer--Zagier formula:
\[
\chi \left( \overline{M}_{g}\right) =\frac{B_{2g}}{2g\left(
2g-2\right) }
\]
Therefore we get
\begin{equation}
C\left( g,d\right) =\frac{\left| B_{2g}\right| d^{2g-3}}{2g\left(
2g-2\right) !}  \label{C(g,d)}
\end{equation}
Note that the genus 0 case is the Aspinwall--Morrison formula%
\index{Aspinwall--Morrison formula}
\[
C\left( 0,d\right) =d^{-3}
\]
(see \cite{am}, \cite{Manin1995}, \cite{Voisin}) and it is easy to
recover its contribution in the series comparing in $F_{0}\left(
t\right) $. For the genus 1 case see \cite{Graber--Pandharipande
1999}: in our particular situation, it turns out that the
non--constant instanton correction is  given by $\left( 1/12 \right) \log \left(
1-e^{-t}\right) $.

 Theorem 4 of \cite{fp} computes $\left\langle
1\right\rangle _{g,0}^{Y}$ in $\left(
\text{\ref{g--instantons}}\right) $. Let  ${\mathbb E}\rightarrow
\overline{M}_{g}$ be the {\it Hodge bundle}, that is, the  rank $g$
vector bundle whose fiber over the  curve $C \in \overline{M}_g$
is given by $H^{0}\left( C,\omega _{C}\right) $ (here $ \omega
_{C}$ is the dualizing sheaf of $C$).

 If $c_{j}\left( {\mathbb E}\right) $ is the
$j^{\text{th}}$ Chern class of ${\mathbb E}$ then
\[
c_{g-1}^{3}\left( {\mathbb E}\right) :=c_{g-1}\left( {\mathbb
E}\right) \wedge c_{g-1}\left( {\mathbb E}\right) \wedge
c_{g-1}\left( {\mathbb E}\right)
\]
is a top form over $\overline{M}_{g}$. A result in \cite{gp}
applied to our Calabi-Yau $Y$ gives
\begin{equation}
\left\langle 1\right\rangle _{g,0}^{Y}=\left( -1\right) ^{g}\int_{\overline{M%
}_{g}}c_{g-1}^{3}\left( {\mathbb E}\right) .
\label{const.instanton}
\end{equation}
Faber--Pandharipande \cite{fp} then show that
\[
\int_{\overline{M}_{g}}c_{g-1}^{3}\left( {\mathbb E}\right)
=\frac{\left| B_{2g}\right| }{2g}\frac{\left| B_{2g-2}\right|
}{2g-2}\frac{1}{\left( 2g-2\right) !}. \  \diamondsuit
\]
\vskip 0.2in

\noindent $\bullet$ {\bf Matching of the free energies.}
 \vskip 0.1in
\begin{theorem}
\label{string exp}
 \textit{Let} $F^{(s)}_{g}\left( t\right) $
\textit{ be as in Definition/Proposition (\ref{gwyd}) } and
$F_{g}\left( t\right)$ \textit{ as in equation (\ref{top.exp.2}).
With the identification
 }$\lambda $ \textit{ and} $t$ \textit{as in} Lemma \ref{ah}
\textit{ we have}
\[
F_{g}^{\left( s\right) }\left( t\right) =F_{g}\left( t\right) ,\ \ \forall g.
\]
 That is, the perturbative expansion (for large $N$) of the free energy
 of $SU(N)$ Chern--Simons theory on $S^3$ is equal
 to the perturbative expansion of the closed II--$A$ theory on
 $Y$.
\end{theorem}

\noindent \textit{Proof of \ref{string exp}:}\\
 In \cite{hkty} it is argued
that $a=1/4$ giving the matching with $F_{0}\left( t\right) $ when
$m=0$. This takes care of $g=0$; the case $g=1$ is immediate.
For $g\geq 2$, since $\left| B_{2g}\right| =\left( -1\right) ^{g+1}B_{2g}$, the relations $%
\left( \text{\ref{g--instantons}}\right) $, $\left( \text{\ref{C(g,d)}}%
\right) $ and $\left( \text{\ref{const.instanton}}\right) $ imply
\begin{eqnarray*}
F_{g}^{\left( s\right) }\left( t\right)  &=&\left( -1\right) ^{g+1}\int_{%
\overline{M}_{g}}c_{g-1}^{3}\left( {\mathbb E}\right) -\sum_{d=1}^{+\infty }%
\frac{\left| B_{2g}\right| d^{2g-3}}{2g\left( 2g-2\right) !}e^{-dt} \\
&=&\frac{B_{2g}}{2g\left( 2g-2\right) !}\left( \frac{\left|
B_{2g-2}\right| }{\left( 2g-2\right) }+\left( -1\right)
^{g}\sum_{d=1}^{+\infty }d^{2g-3}e^{-dt}\right) =F_{g}\left(
t\right) . \ \diamondsuit
\end{eqnarray*}

\phantom{xx}
 \subsection{The matching of expectation values. }~
\vskip 0.1in

Here we
discuss the \textit{matching of the expectation values of observables }in
the two theories of conjecture \ref{GVconjecture}. The conjecture would
be proved if the expectation values for any observable would coincide.
Unfortunately it is not known how to produce a similar ``universal
comparison theorem'' but a general set up to compare some kinds of observables
has been performed and the matching of expectation values has been proved in
some particular cases. In this section we present this strategy and its
striking mathematical consequences.

\noindent The basic idea was already suggested in \cite{gvc} and
then developed in \cite{ov}, \cite{lm}, \cite{kl} and \cite{ls}.
In Chern--Simons theory observables are assigned by
Wilson lines or products of them whose correlation functions are given by $%
\left( \text{\ref{unnorm.exp.value}}\right) $ and $\left( \text{\ref
{correlation}}\right) $ respectively. It is not clear a priori what these
functions correspond to on the topological closed string theory side, but
there are some leads.

\noindent First, Witten's open string interpretation of Chern--Simons theory
also gives the  translation of correlation functions of Wilson observables in
terms of instanton contributions:

\begin{proposition}
\label{open istanton} \textit{An observable in }$SU\left( N\right) $\textit{%
 Chern--Simons gauge theory represented by a link }$\mathcal{L}$\textit{\
corresponds in the Witten open string theory interpretation to the
Lagrangian submanifold }$\mathcal{C}_{\mathcal{L}}$\textit{\ given by the
conormal bundle in }$T^{*}S^{3}|_{\mathcal{L}}$.

\noindent \textit{The non--constant instanton contributions of a type II--}$A
$\textit{\ open string theory with non--compact D--branes wrapped on }$%
\mathcal{C}_{\mathcal{L}}$\textit{\ give a string theory interpretation of
the correlation function of }$\mathcal{L}$\textit{.}
\end{proposition}

\begin{definition}\label{defCK}
Let ${\mathcal{K}}$ be a knot in $S^3$,  parametrized by  $\mathbf{%
q=q}\left( s\right) $ for $s\in \left[ 0,2\pi \right) $. For any
$s$ consider the plane $\pi _{s}\subset \R^{4}\left(
\mathbf{p}\right) $ given by the equations
\begin{eqnarray*}
\sum_{j=1}^{4}q_{j}\left( s\right) p_{j} &=&0 \\
\sum_{j=1}^{4}\dot{q}_{j}\left( s\right) p_{j} &=&0
\end{eqnarray*}
The 3--dimensional submanifold $\mathcal{C}_{\mathcal{K}}:=\coprod_{s}\pi
_{s}$  is called
{\it the conormal bundle} of ${\mathcal{K}}$.
\end{definition}

 \begin{lemma}\label{CK}$\mathcal{C}_{\mathcal{K}}$ \textit{\ is a Lagrangian
submanifold with respect to the symplectic structure induced on}
$T^{*}S^{3}$ \textit{\ by the differential of the Liouville form}
$\vartheta :=\sum_{j=1}^{4}p_{j}dq_{j}$ \textit{\ of} $\R^{8}$.
\end{lemma}

\noindent {\it Proof:}  Consider $%
T^{*}S^{3}$ as embedded in $\R^{8}=\R^{4}\left( \mathbf{q}\right)
\Bbb{\times R}^{4}\left( \mathbf{p}\right) $ by the equations $\left( \text{%
\ref{ctg-bundle}}\right) $. For any $s$ consider the plane $\pi
_{s}\subset \R^{4}\left( \mathbf{p}\right) $ given by the equations
\begin{eqnarray*}
\sum_{j=1}^{4}q_{j}\left( s\right) p_{j} &=&0 \\
\sum_{j=1}^{4}\dot{q}_{j}\left( s\right) p_{j} &=&0
\end{eqnarray*}
Then
\begin{equation}
\vartheta |_{\mathcal{C}_{\mathcal{K}}}=\sum_{j=1}^{4}\dot{q}_{j}\left(
s\right) p_{j}ds=0.  \ \ \diamondsuit \label{lagrangian}
\end{equation}

\noindent \textit{Sketch of the proof of Proposition \ref{open
istanton}:} In \cite{Witten1992} Witten shows that one can
reproduce the correlation function of a Chern--Simons observable
by introducing further D--branes wrapping around a suitable
Lagrangian submanifold of $\widehat{Y}=T^{*}S^{3}$ which is not
the base $S^{3}$ and considering the partition function of the
limit QFT.

\noindent In \cite{gvc}, and \cite{ov} a Wilson line observable
represented by a knot $\mathcal{K}\subset S^{3}$ is associated
with the total space $\mathcal{C}_{\mathcal{K}}$ of  the
``conormal bundle'' defined in definition \ref{defCK}. By lemma
\ref{CK} it is a Lagrangian submanifold with respect to the
symplectic structure induced on $T^{*}S^{3}$ by the differential
of the Liouville form $\vartheta :=\sum_{j=1}^{4}p_{j}dq_{j}$ of
$\R^{8}$ since (\ref{lagrangian}) holds. Then the open string
theory having $T^{*}S^{3}$ as target space and boundary
conditions represented by $M$\ topological D$6$--branes wrapped on $\mathcal{%
C}_{\mathcal{K}}$ is \textit{exactly equivalent to a }$SU\left( M\right) $%
\textit{\ Chern--Simons gauge theory}, since the boundary condition $%
\partial \phi \subset \mathcal{C}_{\mathcal{K}}$, which is the analogue of $%
\left( \text{\ref{bdr.cond.}}\right) $, is satisfied for ``every
bosonic field'' $\phi $. But globally we have now an
``$A$--twisted sigma model'' whose open sector also contains open
strings having one end on $S^{3}$ and the other on
$\mathcal{C}_{\mathcal{K}}$: the non--constant instantons
associated with their world sheet give a non--trivial contribution
to the string amplitude. This means that the low energy limit QFT
is a $SU\left( N\right) \otimes SU\left( M\right) $ gauge theory
which is no longer a
Chern--Simons theory but a deformation of it. Because $\mathcal{C}_{\mathcal{%
K}}\cong \mathcal{K}\times \R^{2}$ and $S^{3}$ is simply
connected, Witten's argument shows that \textit{this partition
function is strictly related with the correlation function of the
original observable associated
with }$\mathcal{K}$\textit{\ in the }$SU\left( N\right) $\textit{%
 Chern--Simons theory on }$S^{3}$. Precisely if $S\left( \mathcal{L}_{%
\mathcal{C}_{\mathcal{K}}}\right) $ is the Chern--Simons action of the $%
SU\left( M\right) $ gauge theory on $\mathcal{C}_{\mathcal{K}}$
defined  as in $\left( \text{\ref{CSaction}}\right) $ then the
partition function of the limit QFT is defined by a Feynman
integration of the following Chern--Simons deformed action:
\begin{equation}
S\left( \mathcal{L}_{\mathcal{C}_{\mathcal{K}}}\right) -\frac{i}{2\pi k}%
\sum_{d}\eta _{d}\log \left( \text{ tr}_{R}\left(
h_{\mathcal{K}}^{d}\right) \right)   \label{CSdeformation}
\end{equation}
where $h_{\mathcal{K}}$ is the holonomy operator on $\mathcal{K}$
with
respect to a connection $\widetilde{A}$ of the $SU\left( M\right) $%
 principal bundle over $\mathcal{C}_{\mathcal{K}}$ and $\eta
_{d}=\pm 1$ for any $d$ (see Corollary \ref{unknot limit} in
Appendix \ref{Witten}).

\noindent The statement of proposition \ref{open istanton} follows
by repeating this construction for every knot in\textit{\
}$\mathcal{L}.\ \ \diamondsuit $

\vskip 0.2in

Then we can try to understand how the
conifold transition acts on those instantons:

\begin{theorem}
\label{open GW}\textit{\ For a suitable link} $\mathcal{L}$, \textit{the
correlation function of the related observable in }$SU\left( N\right) $
\textit{Chern--Simons gauge theory corresponds, on the II--}$A$\textit{\
string theory on }$\mathcal{O}_{\P^{1}}\left( -1\right) \oplus \mathcal{%
O}_{\P^{1}}\left( -1\right) ,$\textit{\ to ``open Gromov--Witten
invariants'' of maps from Riemann surfaces with boundary on
}$\P^{1}$ \textit{ determined by } $\mathcal{L}$.
\index{Gromov--Witten invariants, open}
\end{theorem}

\noindent The class of ``suitable'' links $\mathcal{L}$ in the
statement includes torus knots.

\begin{lemma}\textit{\ Any suitable link (as above)} $\mathcal{L}$,
\textit{determines through the transition a Lagrangian
submanifold}  $\widetilde{\mathcal{C}} \subset Y.$
\end{lemma}
\begin{remark}The
construction in the above lemma has been generalized to all knots by
C. Taubes in \cite{CT}.
\end{remark}

\noindent \textit{Sketch of the Proof of the Lemma for $
{\mathcal{L}}= {\mathcal{K}}$, the un-knot:} We now fix
a knot $\mathcal{K}$, consider the conormal Lagrangian submanifold $\mathcal{%
C}_{\mathcal{K}}$\textit{\ } and study its image, through the conifold
transition, on $Y=\mathcal{O}_{\P^{1}}\left( -1\right) \oplus \mathcal{O%
}_{\P^{1}}\left( -1\right) .$ Such a procedure can easily be
realized when $\mathcal{K}$ is the unknotted knot. Consider in
fact the involution of $\C^{4}\left( x,y,z,t\right) $ given by
\[
\left( x,y,z,t\right) \longmapsto \left( \overline{x},\overline{y},-%
\overline{z},-\overline{t}\right)
\]
Recalling now the chain of transformations given by $\left( \text{\ref{trans}%
}\right) $, $\left( \text{\ref{affinity}}\right) $ and $\left( \text{\ref
{local-smooth}}\right) $ we see that such an involution acts on $\R%
^{8}\left( \mathbf{q,p}\right) $ as follows:
\begin{equation}
\left( q_{1},q_{2},q_{3},q_{4},p_{1},p_{2},p_{3},p_{4}\right) \longmapsto
\left( q_{1},-q_{2},-q_{3},q_{4},-p_{1},p_{2},p_{3},-p_{4}\right)
\label{involution}
\end{equation}
We have then the following three properties:

\begin{enumerate}
\item  $T^{*}S^{3}$ turns out to be fixed by the involution $\left( \text{%
\ref{involution}}\right) $ as follows by its embedding equations
$\left( \text{\ref{ctg-bundle}}\right) $ in $\R^{8}$,

\item  the symplectic form
\[
\omega =d\vartheta =\sum_{j=1}^{4}dp_{j}\wedge dq_{j}
\]
changes its sign under $\left( \text{\ref{involution}}\right) $,

\item  the set of fixed points of $\left( \text{\ref{involution}}\right) $
is given by
\[
\mathcal{F}:=\left\{ \left( \mathbf{q,p}\right)
:q_{2}=q_{3}=p_{1}=p_{4}=0\right\}
\]
\end{enumerate}

\noindent These properties imply that $\mathcal{C}:=\mathcal{F}\cap
T^{*}S^{3}$ is a Lagrangian submanifold with respect to the symplectic
structure induced by $\omega $ on $T^{*}S^{3}$ whose equation in $\R%
^{8}\left( \mathbf{q,p}\right) $ turns out to be
\begin{eqnarray}
q_{1}^{2}+q_{2}^{2}-1 &=&q_{2}=q_{3}=0  \label{Cequations} \\
p_{1} &=&p_{4}=0  \nonumber
\end{eqnarray}
Hence topologically $\mathcal{C}\cong S^{1}\times \R^{2}$ and $\mathcal{%
K}:=\mathcal{C}\cap S^{3}$ is an equator of $S^{3}$ i.e. it is the unknotted
knot on $S^{3}$ and $\mathcal{C}=\mathcal{C}_{\mathcal{K}}$. Recall now
that, by Clemens' theorem \ref{Clemens thm}, the conifold transition can be
locally realized like a surgery by means of the diffeomorphism on boundaries
$\alpha $ represented in $\left( \text{\ref{alfa}}\right) $ whose equations
are
\[
\begin{array}{l}
q_{j}=\frac{u_{j}}{\sqrt{\sum_{i}u_{i}^{2}}} \\
p_{j}=v_{j}\sqrt{\sum_{i}u_{i}^{2}}
\end{array}
\]
Hence the image of $\mathcal{C}$ in the blow up
\[
Y=\mathcal{O}_{\P^{1}}\left( -1\right) \oplus \mathcal{O}_{\P%
^{1}}\left( -1\right) \longrightarrow \overline{Y}
\]
is the strict transform $\widetilde{\mathcal{C}}$ of the subvariety
described in $\overline{Y}$ by conditions $\left( \text{\ref{Cequations}}%
\right) $. Recall that $\overline{Y}$ has local equations $\left( \text{\ref
{real-node}}\right) $ in $\R^{8}\left( \mathbf{u,v}\right) $. Then $%
\widetilde{\mathcal{C}}$ is the strict transform of the 3--dimensional
degenerate hyperquadric of rank 4
\[
u_{1}^{2}+u_{4}^{2}-v_{2}^{2}-v_{3}^{2}=u_{2}=u_{3}=v_{1}=v_{4}=0
\]
Restrict the diffeomorphism $\left( \text{\ref{p-blow}}\right) $
to this hyperquadric: outside of the exceptional fibre it is then
topologically equivalent to $\left( \R_{>0}\times S^{1}\right)
\times S^{1}$. By extending $\left( \text{\ref{p-blow}}\right) $
over the exceptional locus as in $\left(
\text{\ref{local-blow}}\right) $ we get the following topological
interpretation of the strict transform
\[
\widetilde{\mathcal{C}}\cong \R^{2}\times S^{1}
\]
where the second factor $S^{1}$ is an equator of the exceptional locus $%
S^{2} $. Note that  $\widetilde{\mathcal{C}} \cap S^2 = S^1, $
the equator in the exceptional locus $%
S^{2}. \ \diamondsuit $

\vskip 0.2in
 \noindent \textit{Sketch of
the Proof of Theorem \ref{open GW} for $ {\mathcal{L}}=
{\mathcal{K}}$, the un-knot:} \noindent Ooguri and Vafa in
\cite{ov} argue that the Chern--Simons deformation $\left(
\text{\ref{CSdeformation}}\right) $ due to the Wilson line
associated with the unknot can be analytically continued,
for large $N$, to  $-i\Phi \left( \lambda ,t,\mathcal{K}%
\right) $ where
\begin{equation}
\Phi \left( \lambda ,t,\mathcal{K}\right)
=\sum_{d}\frac{\text{tr}_{R}\left( h_{\mathcal{K}}^{d}\right)
+\text{ tr}_{R}\left( h_{\mathcal{K}}^{-d}\right) }{2d\sin \left(
d\lambda /2\right) }e^{-dt/2}
\end{equation}\label{ogw}
$t$ is as in $\left( \text{\ref{lambda}}\right)$, and
$h_{\mathcal{K}}$ the holonomy operator around ${\mathcal{K}}$.
 (This is
formula (3.22) in \cite{ov}, the analytic continuation of (3.14).)
The computation requires a framing of the knot; in \cite{ov} the
trivial framing is chosen.
 Then, using $\mathbf{M}$--theory
duality (see \cite{gva},
\cite{gvb}), they show that  $-i\Phi \left( \lambda ,t,\mathcal{K}%
\right) $ is also the open topological string amplitude on $Y$,
with $D$-branes wrapped around $\widetilde{\mathcal{C}}$ (Section
4.2 and formula (4.4) in \cite{ov}).
 The latter should ``count" holomorphic non--constant
instantons sending  Riemann surfaces with boundary $\Sigma _{g,h}$
onto either the upper or the lower hemisphere of the
exceptional $S^{2}$, with  boundary on $\widetilde{\mathcal{C}}%
\cap S^{2}$.
The terms of the series on the right can be thought of as a sort of%
\index{Gromov--Witten invariants} \textit{Gromov--Witten
invariants of maps from Riemann surfaces with boundary to the
disk.} $\diamondsuit$

\vskip 0.2in
\begin{remark} In theory, the geometric set up that we have presented
for the unknot can
be generalized to every knot or link. In practice the associated
Chern--Simons deformation and the corresponding open instanton
corrections in closed string theory become very intricate and
difficult to compute. In \cite{lm} such a computation is carried
out in the highly non--trivial case of torus knots, again
\textit{showing the conjectured matching of quantities. }The same
result is obtained for other knots and links in
\cite{Ramadevi-Sarkar 2001} and \cite{lmv}.
\end{remark}
It is then natural to ask if one can define mathematically these
open  Gromov--Witten invariants and if they agree with the physics
results mentioned in the above remark.
  At this moment the answer to the first question is not known,
  but, under various assumptions, some results have been obtained
  regarding the second.
  The key observation in
 Katz and Liu \cite{kl} and Li and Song \cite{ls} is that
 $Y$ has a torus action, with nice fixed locus.
  They then assume that the action lifts to the
  compact ``moduli space of maps of open Riemann surfaces"
  and that  localization theorems as in
  \cite{Graber--Pandharipande 1999}, following \cite{MK} hold.
  Then Katz and Liu \cite{kl} and Li and Song \cite{ls} showed that
\[
\Phi \left( \lambda ,t,y\right) =\sum_{d}\frac{y^{d}}{2d\sin
\left( d\lambda /2\right) }e^{-dt/2}
\]
\noindent computes the open Gromov-Witten potential, and that it
is in fact the multiple cover formula of the disc. (Here $t/2$ is
the relative homology class of the (upper) hemisphere with
orientation represented by $y$.)

It turns out that in the ``open" case different torus actions give
rise to different Gromov-Witten potential:  \cite{AKV} showed that
this ambiguity should be expected and that it is related to the
choice of framing on the Chern-Simons side.

\cite{AKV} appeared at the time when these lectures were given.
Many relevant papers have been published since; we do not discuss
them here, as the notes follow closely the lectures.

\section{\label{mt}{ Lifting to }$\mathbf{M}$--theory}

We describe a geometrical construction which gives another
striking evidence for the Gopakumar--Vafa conjecture and reduces to
the conifold geometry by  a ``dimensional reduction''. The main
references for this construction are \cite{ Acharya00, amv} and the more
extensive \cite{aw}.

The geometric construction is suggested by the physical ``lift''
of II--$A$ theories with branes (resp. fluxes), to
$\mathbf{M}$--theory. In our situation, $\mathbf{M}$--theory is
then compactified on 7-dimensional, singular spaces $X_{-r},\
X_{r}$ with special ($G_{2}$) holonomy:
\[
\begin{array}{ccc}
X_{-r} & {\dashleftarrow \dashrightarrow } & X_{r} \\
\downarrow &  & \downarrow \\
\mathbb R^{4}\times S^{2} & <-\mathit{conifold}-> & S^{3}\times \mathbb R^{3}
\end{array}
\]
The vertical maps are essentially Hopf fibrations, the singularities on $%
X_{-r}$ and $X_{+r}$ are related to the presence of branes ( resp.
fluxes) and the special holonomy is needed to preserve the $N=1$
supersymmetry
condition. The conifold transition is lifted to a map between $7$%
-dimensional manifolds (the ``$\mathbf{M}$--theory flop'')%
\index{$\mathbf{M}$--theory, flop} \index{flop,
$\mathbf{M}$--theory}. The physics statement in \cite{Acharya00}, \cite{amv} and
\cite{aw} is that the theory does not go through a singularity
under the $\mathbf{M}$--theory flop: this implies the Gopakumar-Vafa
conjecture for the conifold transition.

In the following subsection we discuss Riemannian holonomy groups; next we
introduce the geometrical construction of the lift for $N=1$ branes. We will
check later its physical consistency with the $\mathbf{M}$--theory lift of $%
II$--A with $N$ branes. Some basic properties of such lifts are stated in
section (\ref{bml}).

\phantom{xx} \subsection{Riemannian Holonomy, $G_2$ manifolds and
Calabi--Yau, revisited.} ~ \vskip 0.1in

The purpose of this section is to fix some notation and basic properties;
details and proofs can be found, for example, in \cite{Joyce 2000}.

\noindent Let $\nabla $ be the Levi-Civita connection on the tangent bundle $%
TM$ of a Riemannian manifold $(M,g)$ and let $p\in M$:

\begin{definition}
The group $%
\text{Hol}_{p}(g)$
\begin{equation}
\text{Hol}_{p}(g):=\text{Hol}_{\nabla }(p)  \label{Hol(g)}
\end{equation}
is the
\index{holonomy, Riemannian}Riemannian \textit{\ holonomy group of }$g$ at $%
p\in M$; Hol$_{\nabla }(p)$ was defined in (\ref{Hol(m)}).
\end{definition}

\noindent It can be seen that when $M$ is connected the holonomy group $%
\text{Hol}(g)$ is a subgroup of $O(\dim M)$, fixed up to conjugation. If $M$
is orientable then $\text{Hol}(g)\subset SO(\dim M)$. If $(M,g,J)$ is a
K\"{a}hler manifold of dimension $2m$, then $\text{Hol}(g)\subset U(m)$.

\begin{theorem}
\textit{A compact K\"{a}hler manifold} $(M,g,J)$ \textit{of complex dimension%
} $m\geq 3$ \textit{is a Calabi--Yau variety if and only if }$\text{Hol}%
(g)=SU(m)$ (for a proof see \cite{Joyce 2000}).
\end{theorem}

\noindent In particular such a $(M,g,J)$ is always projective
algebraic. The following definition, often used in the physics
literature, is then equivalent for $m\geq 3$ to the one given in
(\ref{cym}):

\begin{definition}
\label{cyf} (Calabi--Yau, revisited) A compact Calabi--Yau manifold%
\index{Calabi--Yau, manifold} is a compact K\"{a}hler manifold of dimension $%
2m$, $m \geq 2$, and $%
\text{Hol}(g)=SU(m)$.
\end{definition}

\noindent From the point of view of physics it is the condition $\text{Hol}%
(g)\subseteq SU(m)$ which is relevant, as it preserves the
required supersymmetry. On a $7$-dimensional manifold, the needed
condition becomes $\text{Hol}(g)=G_{2},$ where $G_2$ is defined
below:

\begin{definition}
Let $(x_{1},\ldots ,x_{7})$ be coordinates on $\mathbb R^{7}$ and set%
\[
d\mathbf{x}_{i_{1}\dots i_{r}}=dx_{i_{1}}\wedge \ldots \wedge dx_{i_{r}.}
\]
$G_{2}$ is the Lie subgroup of $GL(7,\mathbb R)$ preserving the $3$--form
\[
\varphi _{0}:=d\mathbf{x}_{123}+d\mathbf{x}_{145}+d\mathbf{x}_{167}+d\mathbf{%
x}_{246}-d\mathbf{x}_{257}-d\mathbf{x}_{347}-d\mathbf{x}_{356}.
\]
\end{definition}

\begin{proposition}
\textit{The following holds:}

\begin{enumerate}
\item  $G_{2}$ \textit{fixes the 4--form} $*\varphi _{0}$ \textit{( }$*$%
\textit{\ is the Hodge star), the Euclidean metric }$g_{0}:=%
\sum_{i=1}^{7}dx_{i}^{2}$\textit{\ and the orientation on }$\mathbb R^{7}$%
\textit{. In particular }$G_{2}\subset SO(7)$\textit{.}

\item  $G_{2}$ \textit{is compact, connected, simply connected and semisimple.}

\item  $\dim G_{2}=14$.
\end{enumerate}
\end{proposition}

\begin{definition}
Let $M$ be an oriented manifold with $\dim M=7$. A 3--form $\varphi _{p}\in
\Lambda ^{3}T_{p}^{*}M$ is \textit{positive} at $p$ if there exists an
oriented isomorphism $T_{p}^{*}M\cong \mathbb R^{7}$ identifying $\varphi
_{p}\text{ with }\varphi _{0}$. Set
\[
\Lambda _{+}^{3}T_{p}^{*}M:=\{\varphi _{p}\in \Lambda ^{3}T_{p}^{*}M\text{
such that }\varphi _{p}\text{ is positive }\}
\]
A 3--form $\varphi $ on $M$ is \textit{positive}
\index{form, positive}if $\varphi |_{p}$ is positive for every point $p\in M$%
; set
\[
\Omega _{+}^{3}(M):=\{\varphi
\text{ such that }\varphi _{p}\in \Lambda _{+}^{3}T_{p}^{*}M,\ \forall p\in
M.\}
\]
\end{definition}

\noindent Note that by definition
\[
\Lambda _{+}^{3}T_{p}^{*}M\cong GL_{+}(7,\mathbb R)/G_{2}
\label{+forms}
\]
A dimensional computation implies
immediately that it is a non--empty open subset of $\Lambda ^{3}T_{p}^{*}M$.
Then a positive 3--form on $M$ is a global section of the open subbundle $%
\Omega _{+}^{3}M$. Fix a positive 3--form $\varphi $ on a Riemannian
7--manifold $(M,g)$. We will write
\[
\text{Hol}(g)\subseteq _{\varphi }G_{2}
\]
when for any $p\in M$ we get
\[
\Phi _{p}\circ \left( \text{Hol}_{p}({g})\right) \circ \Phi
_{p}^{-1}\subseteq G_{2}
\]
where $\Phi _{p}$ is an oriented isomorphism $T_{p}^{*}M\cong \mathbb R^{7}$
representative of the class in $GL_{+}(7,\mathbb R)/G_{2}$ associated with $%
\varphi |_{p}$ via the isomorphism (\ref{+forms}). Since $G_{2}$ is
invariant under conjugation, for any two positive forms $\varphi ,\psi $
\[
\text{Hol}({g})\subseteq _{\varphi }G_{2}\Longleftrightarrow \text{Hol}({g}%
)\subseteq _{\psi }G_{2}
\]
Without loss of generality we then write $\text{Hol}\left( g\right)
\subseteq G_{2}$.

\begin{definition}
$(M,g)$ has%
\index{holonomy, special $G_2$} a $G_{2}$ \textit{holonomy metric} if $%
\text{Hol}({g})=G_{2}$.
\end{definition}

The following properties assure that supersymmetry is preserved:

\begin{proposition}
\textit{Let} $(M,g)$ \textit{be a Riemannian 7--manifold with }$G_{2}$%
\textit{\ holonomy metric}. \textit{Then}

\begin{enumerate}
\item  $g$ \textit{is Ricci flat }

\item  $M$\textit{\ is an orientable spin manifold}

\item  $(M,g)$ \textit{\ has a non-zero covariant spinor.}
\end{enumerate}

(See for example, \cite{Joyce 2000} for a proof of these statements.)
\end{proposition}

\noindent The existence of manifolds with $G_{2}$ holonomy metric
was firstly studied in \cite{Bryant1987} and then solved in
\cite{BS} and
in \cite{Gibbons1990} for non--compact manifolds. Compact manifolds with $%
G_{2}$ holonomy metric were then constructed in \cite{Joyce 1996}. See also Chapter
11 in \cite{Joyce 2000}.

\phantom{xx}
\subsection{Branes and $\M$--theory lifts}~
 \vskip 0.1in

\label{bml}

 II--$A$ string theory may be  regarded  as a
dimensional reduction of an $\mathcal{N}=1$ supersymmetric Lorentz
invariant theory in 11 dimensions:
$\mathbf{M}$--theory\index{$\mathbf{M}$--theory}. (See
\cite{Sen1997}, section 7, for a quick review and references cited
there for details on the argument.) $\mathbf{M}$ --theory was
first proposed  in \cite{Townsend1995} and \cite{Witten1995}, who
observed that the low energy limit of a type II--$A$ string
theory, i.e. a type II--$A$ supergravity theory, can be obtained
by ``Kaluza--Klein''
dimensional reduction \index{Kaluza--Klein, dimensional reduction}of a $%
\mathcal{N}=1$ supersymmetric gravity theory in 11 dimensions.
The reduction is along an $S^{1}$,
 called the $11^{\text{th}}$
\textit{circle}.

When $\mathbf{M}$--theory and II--$A$ are ``compactified'' on
manifolds $M$ and $Y$ respectively, the ``Kaluza--Klein''
dimensional reduction induces an $S^{1}$ fibration $h:M\to Y$.

 If  $\  \su(N)$--branes are ``wrapped'' on a (lagrangian)
submanifold $L\subset Y$, $M$ is singular along $h^{-1}(L)$; the
type of singularity is determined by the group $\su(N)$ (see
Appendix (\ref{sings})) and $h$ is a singular Hopf fibration.
Furthermore, in order to preserve the $\mathcal{N}=1 $
supersymmetry of the theory, $M$ must be a manifold with $G_{2}$
holonomy.

 For
a survey on these topics see, for example, \cite{Johnson 1998} and
\cite{Johnson 2000}\textbf{.} \phantom{xx}

\phantom{xx}
 \subsection{The geometry of the lift for $N=1$ branes}~ \vskip 0.1in

The geometric construction for $N=1$ branes presented here is the
first step towards the $M$-theory lift explained in the following
section. The  equivalence in $M$-theory, and the relations between
parameters stated in Theorem 4.11,  is in fact valid only for $N >
> 0$.

\begin{lemma}
\label{bh} \textit{Fix} $r\ $\textit{in}$\ \R_{>0}$, $\ \C^{4}$
\textit{with coordinates} $(z_{1},z_{2},z_{3},z_{4})$ \textit{and
set}
\begin{eqnarray*}
M_{r} &:&=\{\mathbf{z}\in \C^{4}:\left| z_{1}\right| ^{2}+\left|
z_{2}\right| ^{2}-\left| z_{3}\right| ^{2}-\left| z_{4}\right| ^{2}=r\}, \\
M_{-r} &:&=\{\mathbf{z}\in \C^{4}:\left| z_{1}\right| ^{2}+\left|
z_{2}\right| ^{2}-\left| z_{3}\right| ^{2}-\left| z_{4}\right|
^{2}=-r\}.
\end{eqnarray*}
\textit{Then, topologically:}
\begin{eqnarray*}
M_{r} &\cong &S^{3}\times \C_{(z_{3},z_{4})}^{2}\cong S^{3}\times %
\mathbb R^{4} \\
M_{-r} &\cong &\C_{(z_{1},z_{2})}^{2}\times S^{3}\cong \mathbb %
R^{4}\times S^{3}.
\end{eqnarray*}
\end{lemma}

\noindent The proof of this lemma is presented after the proof of the
following proposition.

\begin{proposition}
\label{gl} \textit{There exists the following geometric lift of the conifold
transition}
\begin{equation}
\begin{array}{ccccccc}
M_{-r} & \cong  & \mathbb R^{4}\times S^{3} & {\leftarrow \cdots \rightarrow
} & S^{3}\times \mathbb R^{4} & \cong  & M_{r} \\
&  & _{h_{-}}\downarrow  &  & \downarrow _{h_{+}} &  &  \\
&  & \mathbb R^{4}\times S^{2} & <-\mathit{conifold}-> & S^{3}\times \mathbb %
R^{3} &  &
\end{array}
\label{lift}
\end{equation}
\textit{where:}

\begin{enumerate}
\item  $h_{-}$ \textit{is the identity on the first factor and the Hopf
fibration on $S^{3}$,}

\item  $h_{+}$ \textit{is the identity on the first factor and the
non-differentiable extension to} $\mathbb R^{3}$ \textit{of the Hopf
fibration on} $S^{3}$.
\end{enumerate}

\noindent \textit{Furthermore} $\mathbb R^{4}\times S^{3}$
\textit{admits a
$G_{2}$ holonomy metric.} \\
\textit{ Note also that   $SU(1)$
singularities are smooth points.}
\end{proposition}

\noindent \textit{Proof of proposition }\ref{gl}\textit{:} The key geometric
observation of the following argument is that $M_{-r}$ and $M_{r}$ are
resolutions of real cones over $S^{3}\times S^{3}$, while $\mathbb %
R^{3}\times S^{3}$ and $S^{2}\times \mathbb R^{4}$ are resolutions of a real
cone over $S^{2}\times S^{3}$. Furthermore the Hopf fibration maps $S^{3}\to
S^{2}$.

\noindent Clemens' Theorem \ref{Clemens thm} describes the conifold
transition as surgery between topological spaces with the same boundary.
This surgery is expressed by the morphism $\alpha $, which is the identity
on $S^{3}\times S^{2}$ (see (\ref{alfa})):
\[
\alpha :\left( \R^{4}\setminus \left\{ \mathbf{0}\right\} \right)
\times S^{2}{\cong }S^{3}\times \left( \R^{3}\setminus \left\{ \mathbf{0%
}\right\} \right) \text{.}
\]
Since:
\[
\begin{array}{l}
\left( \R^{4}\setminus \left\{ \mathbf{0}\right\} \right) \times
S^{2}\cong \R_{>0}\times S^{3}\times S^{2} \\
S^{3}\times \left( \R^{3}\setminus \left\{ \mathbf{0}\right\}
\right) \cong S^{3}\times S^{2}\times \R_{>0}
\end{array}
\]
we can re-write $\alpha $ as
\begin{equation}
\begin{array}{cccc}
\alpha : & \R_{>0}\times S^{3}\times S^{2} & \longrightarrow  &
S^{3}\times S^{2}\times \R_{>0} \\
& \left( \rho ,\mathbf{u},\mathbf{v}\right)  & \longmapsto  & \left( \mathbf{%
u},\mathbf{v},\rho \right) .
\end{array}
\label{alpha}
\end{equation}
As in the previous lemma, we embed $S^{3}\subset \C%
_{(z_{i},z_{i+1})}^{2}$ and consider the compatible%
\index{Hopf fibration} \textit{Hopf fibration:}
\begin{equation}
\begin{array}{cccccc}
h: & S^{3} & \longrightarrow  & \P_{\C}^{1}\cong S^{2} &  &  \\
& \left( z_{i},z_{i+1}\right)  & \longmapsto  & \left[ z_{i},z_{i+1}\right] =
& \left[ \lambda z_{i},\lambda z_{i+1}\right] , & \ \lambda \in \mathbb %
C^{*}.
\end{array}
\label{hopf}
\end{equation}
Then the following diagram:
\begin{equation}
\begin{array}{ccc}
\R_{>0}\times S^{3}\times S^{3} & \stackrel{%
\widetilde{\alpha }}{\longrightarrow } & S^{3}\times S^{3}\times
\R_{>0}
\\
_{h_{3}}\downarrow  &  & \downarrow _{h_{2}} \\
\R_{>0}\times S^{3}\times S^{2} & \stackrel{\alpha
}{\longrightarrow } & S^{3}\times S^{2}\times \R_{>0}
\end{array}
\label{lift1}
\end{equation}
commutes, where
\[
\begin{array}{l}
h_{3}:=\text{Id}_{\R_{>0}}\times \text{Id}_{S^{3}}\times h \\
h_{2}:=\text{Id}_{S^{3}}\times h\times \text{Id}_{\R_{>0}} \\
\widetilde{\alpha }\left( \rho ,\mathbf{u},\mathbf{u^{\prime }}\right)
:=\left( \mathbf{u},\mathbf{u^{\prime }},\rho \right) .
\end{array}
\]
Note that while $h_{3}$ can be smoothly extended to a fibration
\[
h_{-}:=\text{Id}_{\mathbb R^{4}}\times h:\mathbb R^{4}\times
S^{3}\longrightarrow \mathbb R^{4}\times S^{2},
\]
this is not true for $h_{2}$. There is however a topological extension $%
h_{+}$ of $h_{2}$. The extensions $h_{-}$ and $h_{+}$ then give the diagram (%
\ref{lift}) in the statement.

\noindent \cite{BS} and \cite{Gibbons1990} explicitly describe a $%
G_{2}$ holonomy metric on $M:=S^{3}\times \mathbb R^{4}.$

The metric in \cite{Gibbons1990} is a smooth extension of the
metric on the cone over $S^3 \times S^3$. Bryant and Salamon
\cite{BS} consider $SU(2) \cong S^3$, and the quaternions $\mathbb
H \cong \mathbb R^4$ as a cone over $SU(2)$.
Then $S^3 \times \mathbb R^4 \cong (SU(2) \times SU(2) \times %
\mathbb H)/ SU(2)$, with $SU(2)$ acting on the right, is a rank four vector
bundle on $SU(2)$. With this latter representation, it is evident that there are
other two resolutions of the cone over $S^3 \times S^3$:
\[
(\mathbb H \times SU(2) \times SU(2))/ SU(2) \cong \mathbb R ^4 \times S^3
\text{, }( SU(2) \times \mathbb H \times SU(2) )/ SU(2).
\]
The third manifold fibers, via the Hopf fibration,
 to the ``flopped" local Calabi--Yau $Y_+$ of the resolved conifold $Y$
 (see \ref{cyflop}); we have then that third branch in  Figure \ref{como1}
 (see also \cite{Manzoni 1842}).  $\diamondsuit $
 \begin{figure}[h]
     \begin{center}
     \scalebox{1}{\includegraphics{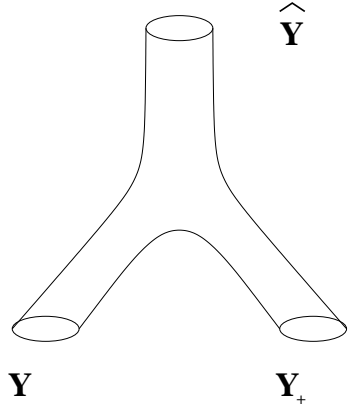}}
     \caption{The three branches of the moduli.
({\it ``Quei rami del lago di Como..."})}\label{como1}
  \end{center}
     \end{figure}

\vskip 0.2in

\noindent \textit{Proof of Lemma }\ref{bh}\textit{: \ } Let $%
(z_{1},z_{2},z_{3},z_{4})$ be coordinates in $\mathbb C^{4}$; for every
positive real number $r$ set:
\[
M_{r}:=\left\{ \mathbf{z}\in \C^{4}:\left| z_{1}\right|
^{2}+\left| z_{2}\right| ^{2}-\left| z_{3}\right| ^{2}-\left|
z_{4}\right| ^{2}=r\right\} .
\]
Then,
\[
\begin{array}{cccc}
\phi _{+}: & M_{r} & \longrightarrow  & S^{3}\times \C^{2} \\
& \left( z_{1},z_{2},z_{3},z_{4}\right)  & \longmapsto  & (\frac{z_{1}}{\rho
_{+}},\frac{z_{2}}{\rho _{+}},{z_{3}}\cdot {\rho_+},{z_{4}}\cdot {\rho_+})
\end{array}
\]
is an isomorphism, where $\rho _{+}:=\sqrt{\left| z_{1}\right| ^{2}+\left|
z_{2}\right| ^{2}}=\sqrt{r+\left| z_{3}\right| ^{2}+\left| z_{4}\right| ^{2}}%
.$

\noindent Similarly for $M_{-r}$. $\diamondsuit $

\newpage

 \subsection{$\M$--theory lifts and $\M$--theory flops}~
\vskip 0.1in

\begin{theorem}
\cite{Acharya1999, Acharya00}, \cite{amv}, \cite{aw} {\it There
exists a commutative diagram}
\begin{equation}
\begin{array}{ccc}
M_{-r} &  & M_{r} \\
_{\pi _{-}}\downarrow  &  & \downarrow _{\pi _{+}} \\
X_{-} &  & X_{+} \\
_{h_{-}^{\left( N\right) }}\downarrow  &  & \downarrow
_{h_{+}^{\left(
N\right) }} \\
\mathbb R^{4}\times S^{2} & <-\mathit{conifold}-> & S^{3}\times \mathbb %
R^{3}.
\end{array}
\label{lift-N}
\end{equation}
 \textit{where, }
 \begin{enumerate}
 \item $M_{-r}$ {\it and }$M_r$ {\it  are as in Proposition
 \ref{gl}.}
\item $X_{-}$\textit{\ and }$X_{+}$\textit{\ are} $G_{2}$ \textit{%
holonomy spaces.}
\item  $(S^{3},0)\subset X_{+}$ \textit{is a
locus of} $A_{N-1}$ \textit{singularities.}
\item {\it The  diagram  is physically consistent, for large} $N$ {\it with the
}
 $\mathbf{M}$--\textit{theory
lift to } $X_-$  ({\it resp. } $X_+$) \textit{\ of } $N$ {\it RR
fluxes on } $\co_\pl (-1)\oplus \co_\pl (-1) $  ({\it resp.}
$\su(N)$ \textit{ branes on } $ T^*S^3$).
\item  \textit{The surjections }$h_{-}^{\left( N\right) },h_{+}^{\left(
N\right) }$\textit{\ give rise to the fluxes and branes,
respectively, for the type II--}$A$\textit{\ string theories
obtained by dimensional reduction on the two sides of the conifold
transition.}
\item $\mathbf{M}$\textit{--theory compactified on} $X_{-}$
\textit{is equivalent to }$\mathbf{M}$\textit{--theory on }
$X_{+}$.
\end{enumerate}
{\it Thus, there is no ``phase'' transition between} $X_-$ {\it and
} $X_+$ {, exactly as when II--}$A$ {\it is compactified on
Calabi-Yau varieties related by a ``flop'' (see
\cite{Witten1993}). \\
Hence the term }$\M$--{\it {\bf theory  flop.}}
{\it This physics description is valid only for large} $N$.
\end{theorem}

The physics statement in \cite{Acharya00}, \cite{amv} and
\cite{aw} is that the theory does not go through a singularity
under the $\mathbf{M}$--theory flop: this implies the Gopakumar-Vafa
conjecture for the conifold transition.

\vskip 0.1in

\noindent \textit{Sketch of the proof:}
 At the time of this lecture the works  \cite{Acharya1999, Acharya00}, \cite{amv}
 were in print, while the main results of \cite{aw}
had just been recently announced. The geometric lift $\left(
\text{\ref {lift}}\right) $ gives an $\mathbf{M}$--theory lift of
II--$A$ string theories when $N=1$. The singularity of the map
$h_{+}$ denotes the presence of branes.

\noindent To get the $\mathbf{M}$--theory lift with $N$ D--branes wrapped on $%
S^{3}\times \{0\}\subset S^{3}\times \mathbb R^{3}$ we need to
introduce corresponding singularities on $M_{r}$ (see Section
\ref{bml}). We do so by defining a suitable action of the group of
$N^{\text{th}}$ roots of unity on $\mathbb C^{4}$: the induced
action on $M_{-r}$ will give $N$ units of RR flux on $\mathbb
R^{4}\times S^{2}$.

\noindent Let $\Gamma _{N}:=\Bbb{Z}/N\Bbb{Z}$ act on $\C^{4}$ as
\begin{equation}
\begin{array}{lll}
\Gamma \times \C^{4} & \longrightarrow  & \C^{4} \\
\left( n,\mathbf{z}\right)  & \longmapsto  & \left( z_{1},z_{2},\xi
_{n}z_{3},\xi _{n}z_{4}\right)
\end{array}
\label{action}
\end{equation}
where $\xi _{n}:=\exp \left( 2\pi in/N\right) .$ The complex plane $%
F=:\{z_{3}=z_{4}=0\}$ is the fixed locus of $\Gamma $. Recall that $%
M_{-r}\cong {{\mathbb C^{2}}_{(z_{1},z_{2})}\times S^{3}}$ and $M_{r}\cong {%
S^{3}\times {\mathbb C^{2}}_{(z_{3},z_{4})}}$. Then:
\[
F\cap M_{-r}=\emptyset ,\ \ F\cap M_{r}=S^{3}\times \left\{ \mathbf{0}%
\right\} .
\]
\noindent The quotient

\[
M_{-r}\cong {\mathbb C^{2}}_{(z_{1},z_{2})}\times S^{3}\longrightarrow
M_{-r}/\Gamma \cong {\mathbb R^{4}}\times (S^{3}/\Gamma ):=X_{-}
\]
is smooth; $(S^{3}/\Gamma )$ is called a \textit{lens space} and is denoted by $L\left( N\right) .$
Furthermore, since the $\Gamma $%
--action restricts to the fiber of the Hopf fibration, the map $h_{-}$ in $%
\left( \text{\ref{lift}}\right) $ can be factorized through the canonical
projection $\pi _{-}$ as follows
\[
\begin{array}{lll}
M_{-r} & \stackrel{h_{-}}{\longrightarrow } & \R^{4}\times S^{2} \\
_{\pi _{-}}\searrow  &  & \nearrow _{h_{-}^{\left( N\right) }} \\
& X_{-}. &
\end{array}
\]
On the other hand the quotient
\[
M_{r}\cong S^{3}\times {\mathbb C^{2}}_{(z_{3},z_{4})}\longrightarrow
M_{r}/\Gamma \cong S^{3}\times \left( {\mathbb R^{4}}/\Gamma \right) :=X_{+}
\]
contains an $S^{3}$ \textit{\ of singular points.} Furthermore, since the $%
\Gamma $--action restricts to the fiber of the Hopf fibration, the map $h_{+}
$ in $\left( \text{\ref{lift}}\right) $ can topologically be factorized
through the canonical projection as follows
\[
\begin{array}{lll}
M_{r} & \stackrel{h_{+}}{\longrightarrow } & S^{3}\times \mathbb R^{3} \\
_{\pi _{+}}\searrow  &  & \nearrow _{h_{+}^{\left( N\right) }} \\
& X_{+}. &
\end{array}
\]
$\mathbb R^{4}/\Gamma $ is an $A_{N-1}$ singularity, with gauge
group $SU(N)$
(see Appendix \ref{sings}). In fact with the change of coordinates $%
w_{3}=z_{3},w_{4}=\sqrt{-1}\cdot \overline{z}_{4}$, the action becomes: $%
(w_{3},w_{4})\to (\xi w_{3},\xi ^{-1}w_{4})$ as described in Appendix \ref
{sings}. This is the geometric incarnation of the $\M$--theory lift with $%
SU(N)$--branes wrapped on $S^{3}$ (see Section \ref{bml}).

\noindent Furthermore the non--singular $\Bbb{Z}_{N}$--quotient
(on the left of diagram $\left( \text{\ref{lift-N}}\right) $)
gives rise to $N$ units of RR flux. In fact, if $V(-r)$ is the
volume of $S^{3}\times \{\mathbf{0}\}$, then
$vol(S^{2})=vol(S^{3}/\Gamma )=V(-r)/N$.

 Recall that there exists a $G_{2}$ holonomy metric
(see \cite{BS}, \cite{Gibbons1990}) on $M:=S^{3}\times \mathbb R^{4}$%
. There is a precise description of the isometry group on $M$ and the action
of $\Gamma $ is included in this subgroup. Hence the quotients $X_{-}$ and $%
X_{+}$ are also $G_{2}$ holonomy spaces.

It is worth pointing out that the equivalence of the theory and
the relations between the physical parameters derived in
\cite{amv} {\it are only valid for large} $N$. The equivalence of
the theories also implies  the relations between K\"{a}hler
modulus of $Y$ and the parameters of the Chern-Simons theory
conjectured by Gopakumar and Vafa (see \cite{amv}).

 On the other hand, the asympotics of the $G_{2}$ metric is not
what  would be expected from the II--$A$ situation; based on this
observation Atiyah, Maldacena and Vafa conjectured the existence
of a deformation of the $G_{2}$ metric with such properties (see
\cite{amv}). This was later shown in \cite{bggg}. $\diamondsuit $

\section{\label{sings} Appendix: Some notation on singularities and their
resolutions}

Here we  adopt the same notation and terminology introduced in
\cite{Reid1980}, \cite{Reid1983} and \cite{Reid1987(a)}.

\begin{definition} A Weil divisor $D$ on a complex, normal and quasiprojective variety
$\overline{Y}$  is $\Bbb Q$-Cartier
if, for some $r\in \Bbb{Z}$, $rD$ is  a Cartier divisor(i.e. $D
\in Pic \left( \overline{Y}\right) \otimes \Bbb{Q}$).
\end{definition}
If $\overline{Y}$ is smooth then any Weil divisor is Cartier.
\begin{definition} A $\overline{Y}$ be a complex, normal and quasiprojective variety
is $\Bbb Q$-factorial if any Weil divisor is $\Bbb Q$-Cartier.
\end{definition}
\begin{definition} Let $\overline{Y}$ be a complex, normal and quasiprojective variety and
$K_{\overline{Y}} $ be its \textit{canonical divisor} which is in
general a Weil divisor. $\overline{Y}$ has \textit{canonical} (
respectively \textit{terminal])
singularities
}if:%
\index{canonical singulaty}
\index{terminal singularity}
\index{singularity, canonical/terminal}

\begin{enumerate}
\item[i)]  $K_{\overline{Y}}$ is $\Bbb Q$-Cartier.

\item[ii)]  given a smooth resolution $f:Y\longrightarrow \overline{Y}$ then
\[
rK_{Y}\equiv f^{*}K_{\overline{Y}}+\sum_{i}a_{i}E_{i}
\]
where $\equiv $ means ``linearly equivalent'', $E_{i}$ are all the
exceptional divisors of $f$ and $a_{i}\geq 0$ (respectively $a_{i}>0$).
\end{enumerate}
\end{definition}
\noindent The smallest integer $r$ for which such conditions hold is called%
\textit{\ }the\textit{\ (global) index of }$\overline{Y}$ and the
smallest $r^{\prime }$ for which $r^{\prime }K_{\overline{Y}}$ is
Cartier in a neighborhood of $P\in \overline{Y}$ is called\index{index,
of a variety}
\index{index, of a singularity} the \textit{index of the singularity }$P$%
\textit{.}

\noindent The divisor $\Delta :=\sum_{i}a_{i}E_{i}$ is called
\index{discrepancy, of a resolution} \index{crepant resolution}the
\textit{discrepancy of the resolution }$f$.

\noindent If
$\Delta \equiv 0$ then $f$ is called a \textit{crepant resolution of }$\overline{Y}$%
\textit{.}

We are interested in transitions of Calabi-Yau manifolds: in
particular, if at a point in the complex moduli space $\overline{Y}$ is
singular and $K_{\overline{Y}} \equiv 0$,  its birational resolution
should be crepant to preserve the Calabi-Yau condition on the
canonical bundle.

\begin{definition}
\label{moricone} (see for example, \cite{CKM})%
\index{cone, of Mori} \index{Mori cone}\newline
By $NE(Y)\subset \Bbb{R^{\ell }}$ we denote the cone generated (over $\Bbb{%
R}_{\geq {0}}$) by the effective cycles of (complex) dimension 1,
mod. numerical equivalence.\newline $\overline{NE(Y)}$ is the
closure of $NE(Y)\subset \Bbb{R^{\ell }}$ in the finite
dimensional real vector space $\Bbb{R^{\ell }}$ of all cycles of
complex dimension 1, mod. numerical equivalence.

Note that $\ell =rk(Pic(Y))$, and in the cases of Calabi-Yau manifolds, $%
\ell =b_{2}(Y)$, the second Betti number of $Y$.
\end{definition}

\begin{definition} A birational contraction $f: Y \to \overline{Y}$
is called {\it primitive extremal} if the numerical class of a
fiber of $f$ is on a ray of the Mori cone $NE(Y)$.
\end{definition}

\vskip 0.2in

\noindent {\bf Examples}

\noindent \textit{The surface case. }Let $X$ be a surface. It can be
proved that a point $P\in X$ is a terminal singularity if and only
if it is non-singular. Moreover the canonical (non--terminal)
singular
points are given by the%
\index{singularity, of Du Val}
\index{Du Val singularity}
\index{DV point} \textit{Du Val singularities }(DV points) which are
classified as follows in terms of their local equations
\begin{eqnarray*}
A_{n} &:&x^{2}+y^{2}+z^{n+1}=0\ ,\ n\geq 1 \\
D_{n} &:&x^{2}+y^{2}z+z^{n-1}=0\ ,\ n\geq 4 \\
E_{6} &:&x^{2}+y^{3}+z^{4}=0 \\
E_{7} &:&x^{2}+y^{3}+yz^{3}=0 \\
E_{8} &:&x^{2}+y^{3}+z^{5}=0
\end{eqnarray*}
\index{singularity, A-D-E}In particular each of them admits a
crepant resolution whose exceptional locus is composed of a set of
$\left( -2\right) $--curves (i.e. rational curves admitting
self--intersection index $-2$) whose configurations are dually
represented by the homonymous\index{Dynkin diagrams} Dynkin
diagrams: these are particular examples of Hirzebruch--Jung
strings (see \cite{Barth1984}, chapters I and III).

Note that an \index{node (ordinary double point)}
\index{singularity, ordinary double point}ordinary double point is
represented by $A_{1}$ and admits a crepant resolution whose
exceptional locus is given by a unique $\left( -2\right) $--curve.
 This equation is generalized to the threefold case
in definition \ref{a1t}.

Each of the above singularities can be  described as a quotient of
$\C^2$ by a discrete subgroup $\Gamma \subset SL(2)$. For $A_n$,
$\Gamma$ is the cyclic group of order $n+1$ generated by a
primitive $n$-th root of unity $\xi$; the action on $\C^2$ sends
$(w_1, w_2) \to (\xi w_1, \xi ^{-1} w_2)$ (see
\cite{Slodowy1980}).

\vskip 0.2in

\noindent \textit{The threefold case.} Let $X$ be a threefold and $P\in X$ be a
canonical singular point of index $r$. A first important fact is that there
exists a finite Galois covering $Y\longrightarrow X$ with group $\Bbb{Z}/r$
which is \'{e}tale in codimension 1 and such that $Y$ is locally canonical
of index 1 (see \cite{Reid1980}, corollary (1.9)).

\begin{definition} $P\in X$ is a
 \textit{compound Du Val singularity} (cDV point)
 if the restriction
to a surface section is a Du Val surface singularity.
\index{singularity, compound Du Val} \index{compound Du Val
singularity} \index{cDV points}
\end{definition}
The advantage of these kind of singularities is that they admit a
simultaneous small resolution, as studied by several authors (see
e.g. \cite{Reid1983}, \cite{Pinkham 1983}, \cite{Morrison 1985}, \cite{Friedman1986}).
The idea is that of thinking of an analytic neighborhood of an
isolated cDV point as the total space of a 1--parameter family
of deformations of the section over which we get a DV point. The
total space of the induced 1--parameter family of deformations of
a given resolution of such a DV point is then a small resolution
of the starting cDV point. One can now apply the theory of
simultaneous resolutions of DV points on surfaces \cite{brie1966},
\cite{brie1968}, \cite{Tyurina1970}.

The Main Theorem in \cite{Reid1983} states that:
\begin{enumerate}
\item[i)]  $P\in X$ is a terminal singularity of index $r$ if and only if
the local $r$--fold cyclic covering $Y\longrightarrow X$ has only
isolated \textit{compound Du Val singularities}.%

\item[ii)]  if $X$ admits at most canonical singularities then there exists
a crepant partial resolution $S\longrightarrow X$ such that $S$ admits at
most isolated terminal singularities.
\end{enumerate}

\noindent These results allow one to reduce the problem of resolving canonical
singularities to that of resolving cDV points, up to partial resolutions and
finite coverings.

\section{\label{Candelas'ex} Appendix: More on the Greene-Plesser
construction}

Here we will quickly sketch an example supporting the
Greene-Plesser construction explained in \cite{CdOFKM}, \cite{GP} and
\cite{Morrison 1999}.

\noindent Let $\overline{Y}_{1}$ be the degree 8 weighted hypersurface of $%
\P\left( 1,1,2,2,2\right) $ and $Y_{1}$ be the desingularization
induced by blowing up the singular locus of $\P\left( 1,1,2,2,2\right) $%
. Here $\phi $ is a primitive contraction of type $III$ and the
transition can be completed by considering the embedding of
$\P\left( 1,1,2,2,2\right) $ in $\P^{5}$ by means of the linear
system $\mathcal{O}\left( 2\right) $. The image of $\P\left(
1,1,2,2,2\right) $ is a rank 3 hyperquadric of $\P^{5}$. Hence the
image of $\overline{Y}_{1}$ is the complete
intersection of this hyperquadric with the generic quartic hypersurface of $%
\P^{5}$. By smoothing the hyperquadric we get $\widehat{Y}_{1}$.
Following the idea of \cite{GP} the mirror partners\index{mirror,
partners} may be found by taking the quotient with the subgroups
of automorphisms preserving the holomorphic 3--form. Since the
hypersurfaces' cohomology can be completely described by means of
Poincar\'{e} residues (see \cite{pg}) these subgroups are
respectively given by
\begin{eqnarray*}
G &:&=\left\{ \left( a_{0},\ldots ,a_{4}\right) \in \left( \Bbb{Z}%
_{8}\right) ^{2}\times \left( \Bbb{Z}_{4}\right) ^{3}:\sum a_{i}\equiv
0\left( 8\right) \right\}  \\
H &:&=\left\{ \left( b_{0},\ldots ,b_{5}\right) \in \left( \Bbb{Z}%
_{4}\right) ^{2}\times \left( \Bbb{Z}_{2}\right) ^{4}:b_{0}+b_{1}\equiv
b_{2}+\cdots +b_{5}\equiv 0\left( 4\right) \right\}.
\end{eqnarray*}
We denote by $a_{i},b_{j}$ the least non--negative integers
representing the associated congruence class in $\Bbb{Z}_{n}$. Hence the
mirror partner $\widehat{Y}_{2}$ of $\widehat{Y}_{1}$ may be
obtained by an $H$--invariant complete intersection of bidegree
$\left( 2,4\right) $ in $\P^{5}$ via the
desingularization of the quotient $\P^{5}/H$ where $H$ acts on $\P%
^{5}$ as follows
\[
\begin{array}{lll}
\left( H/\Delta _{H}\right) \times \P^{5} & \longrightarrow  & \P%
^{5} \\
\left( \mathbf{b},\mathbf{x}\right)  & \longmapsto  & \left( \beta
_{j}x_{j}\right)
\end{array}
\]
where
\[
\beta _j :=\left\{
\begin{array}{ll}
\exp \left( \frac{b_{j}\pi i}{4}\right)  & \text{ for }j=0,1 \\
\pm 1 & \text{ otherwise}
\end{array}
\right.
\]
and $\Delta _{H}$ is the subgroup of $H$ giving a trivial action on $\P%
^{5}$, i.e.
\[
\Delta _{H}:=\left\{ \left( 0,\ldots ,0\right) ,\left( 2,2,1,\ldots
,1\right) \right\}
\]
On the other hand the mirror partner $Y_{2}$ of $Y_{1}$ may be
obtained by a $G$--invariant hypersurface of degree 8 in $\P\left(
1,1,2,2,2\right) $ via the desingularization of the quotient
$\P\left( 1,1,2,2,2\right) /G$ where $G$ acts on $\P\left(
1,1,2,2,2\right) $ as follows
\[
\begin{array}{lll}
G/\Delta _{G}\times \P\left( 1,1,2,2,2\right)  & \longrightarrow
&
\P\left( 1,1,2,2,2\right)  \\
\left( \mathbf{a},\mathbf{x}\right)  & \longmapsto  & \left( \alpha
_{j}x_{j}\right)
\end{array}
\]
where
\[
\alpha _{j}:=\left\{
\begin{array}{ll}
\exp \left( \frac{a_{j}\pi i}{8}\right)  & \text{ for }j=0,1 \\
\exp \left( \frac{a_{j}\pi i}{4}\right)  & \text{ otherwise}
\end{array}
\right.
\]
and $\Delta _{G}$ is the diagonal subgroup of $G$, which is
\[
\Delta _{G}:=\left\{ \left( a,\ldots ,a\right) :0\leq a\leq 3\right\}
\]
It can be checked that there is a birational equivalence between $\widehat{Y}%
_{2}$ and $Y_{2}$ representing a mirror partner of our transition%
\index{transition, reverse}.

\section{\label{blackholes} Appendix: More on transitions in superstring
theory}

Strominger gave in \cite{Strominger1995} a physical explanation of how
to resolve the conifold
singularities of the moduli space of classical string vacua by means of
massless Ramond--Ramond (RR) black holes%
\index{black hole, RR}. More precisely, the possible compactifications of a
10--dimensional II--$B$ string theory to 4 dimensions on a Calabi-Yau manifold $%
Y $ may be parametrized by the choice of the complex structure
characterizing $Y$. Such a choice may be described by the periods of a
holomorphic 3--form $\Omega $ over a suitable symplectic basis of $%
H_{3}\left( Y,\Bbb{Q}\right) $ (see \cite{WLP} and \cite{Strominger1990}
for detailed notation in a $N=2$, 4--dimensional supergravity
theory and in special geometry) which can be considered as
projective coordinates of the moduli space $\mathcal{M}\left(
Y\right) $ of complex structures. The complex codimension 1 locus
defined in $\mathcal{M}$ by the vanishing of one of those periods
is composed of singular complex structures generically
geometrically realized by a conifold. In fact the generic
singularity is given by an ordinary double point. Note that the
associated vanishing cycle is represented by the 3--cycle of the
symplectic basis corresponding to the vanishing period.

\noindent Such singularities induce a polydromic behavior for the components
of the self--dual 5--form giving the classical field. Following an analogous
construction given in \cite{SeibergWitten} and applied in the completely
different context of $N=2$ supersymmetric Yang--Mills theory, Strominger
resolved this problem by means of a low--energy effective Wilsonian field
defined by including the light fields associated with extremal black
3--branes which can wrap around the vanishing 3--cycles and are always
contained in a 10--dimensional compactified type II--$B$ theory (see
\cite{Horowitz--Strominger 1991}). These 3--branes represent black holes whose mass
is proportional to the volume of the vanishing cycles they wrap around.
Hence they are massless at the conifold and by integrating out the smooth so
defined Wilsonian field we get exactly the polydromic behavior of the
classical field. This is enough to ensure that the theory may be smoothly
extended to the conifold.

\noindent On the other hand, in the case of a 10--dimensional
compactified type II--$A$ theory we get a similar picture by
taking the periods of a complexified K\"{a}hler form $\omega \in
H^{2}\left( Y,\C \right)
=H^{1,1}\left( Y\right) $ over a suitable basis of $H_{2}\left( Y,\Bbb{Q}%
\right) $ as projective coordinates of the moduli space $\mathcal{M}%
^{^{\prime }}\left( Y\right) $ of all possible K\"{a}hler structures on $Y$
(which parametrizes all the possible compactifications of a 10--dimensional
II--$A$ string theory to 4 dimensions on the Calabi-Yau manifold $Y$). We now get
black 2-branes (see \cite{Horowitz--Strominger 1991}) which can wrap around
vanishing 2--cycles and represent massless black holes at the conifold. Since in
this case these massless states are a result of large instanton corrections
the resolution of singularities can be obtained by passing to the dual II--$B$
compactification on a mirror model $Y^{\circ }$ of $Y$ and by proceeding
as before.

\section{\label{diffdef}{ Appendix: Principal bundles, connections etc}}

Here we review some terminology, concepts and properties from
differential geometry: for more details see, for example,
\cite{Helgason1978}, \cite{Poor 1981} and \cite{Warner}.

\begin{definition}
\label{LRactions}\label{principal bdl}Let $G$ be a Lie group. A left ( resp. right)
action of $G$ on a manifold $M$ is a homomorphism (resp. anti--homomorphism) to
the group of diffeomorphisms of $M$%
\[
L\ \left(
\text{resp. }R\right) :G\longrightarrow Diff\left( M\right)
\]
\end{definition}

\noindent In particular for every $\sigma ,\tau \in G$ we have $L\left(
\sigma \right) $ $\circ L\left( \tau \right) =L\left( \sigma \tau \right) $
(resp. $R\left( \sigma \right) $ $\circ R\left( \tau \right) =R\left( \tau
\sigma \right) $ ).

\begin{definition}
An action is \textit{free }if $id$ is the unique element of $G$ whose image in
$Diff(M)$ admits a fiexd point. Note that if the $G$--action is free then it is
an injection of $G$ into $Diff\left( M\right) $.
\end{definition}

\begin{definition}
A \textit{principal }$G$\textit{--bundle}
\index{bundle, principal}on a manifold $M$ is a manifold $P$ on which $G$
acts freely on the right together with a smooth, surjective map $\pi
:P\rightarrow M$ such that

\begin{enumerate}
\item  for every point $m\in M$ there is a \textit{local trivialization of }$%
P$ i.e. an open neighborhood $\{U_{a}\}$ and a local diffeomorphism $\varphi
_{U_{a}}:\pi ^{-1}\left( U_{a}\right) \stackrel{\cong }{\rightarrow }%
U_{a}\times G$ making the following diagram commutative
\begin{equation}
\begin{array}{ccc}
\pi ^{-1}\left( U\right)  & \stackrel{\varphi _{U}}{\longrightarrow } &
U\times G \\
^{\pi }\downarrow  & \stackrel{pr_{1}}{\swarrow } &  \\
U &  &
\end{array}
\label{triv}
\end{equation}

\item  $\pi $ is $G$--invariant i.e. for every $p\in P$ and every $\sigma
\in G$
\[
\pi \left( p\sigma \right) =\pi \left( p\right)
\]
where $p\sigma :=R\left( \sigma \right) p$.
\end{enumerate}
\end{definition}

\begin{remark}
\label{g-fibre}For a principal bundle $\left( P,\pi \right) $ the map $\pi $
is a submersion, implying that
\[
\mathcal{V}_{p}P:=\ker \left( d_{p}\pi \right) =T_{p}\pi ^{-1}\left( \pi
\left( p\right) \right)
\]
for every $p\in \pi ^{-1}\left( \pi \left( p\right) \right) $. Set $m:=\pi
\left( p\right) \in M$ and let $\left( U,\varphi _{U}\right) $ be a local
trivialization of $P$ near $m$. The commutative diagram (\ref{triv}) allows us to
define a diffeomorphism $\sigma _{m}^{U}$ such that
\[
\left( \sigma _{m}^{U}\right) ^{-1}:=\left( \varphi _{U}^{-1}\right) \mid
_{\left\{ m\right\} \times G}:G\stackrel{\cong }{\longrightarrow }\pi
^{-1}\left( m\right)
\]
Its differential gives the isomorphism
\[
d_{p}\sigma _{m}^{U}:T_{p}\pi ^{-1}\left( m\right) \stackrel{\cong }{%
\longrightarrow }T_{\sigma _{m}^{U}\left( p\right) }G
\]
On the other hand by differentiating the automorphism $r_{\sigma }$ of $G$,
given by right multiplication by $\sigma \in G$, we get the isomorphism
\[
d_{id}r_{\sigma }:\mathfrak{g}\cong T_{id}G\stackrel{\cong }{\longrightarrow }%
T_{\sigma }G
\]
where $\mathfrak{g}$ is the Lie algebra associated with $G$ whose
elements are all the left invariant vector fields on $G$. Hence
for every $p\in \pi ^{-1}\left( m\right) $ we get the isomorphism
\[
d_{p}\left( r_{\sigma _{m}^{U}\left( p\right) }^{-1}\circ \sigma
_{m}^{U}\right) :\ker \left( d_{p}\pi \right) \stackrel{\cong }{%
\longrightarrow }\mathfrak{g}
\]
This suffices to conclude that \textit{the vertical bundle }$\mathcal{V}P$%
\textit{\ associated with the principal }$G$--\textit{bundle
}$\left( P,\pi \right) $\textit{\ is a vector bundle whose
standard fibre is the Lie algebra $\mathfrak{g}$ associated with
}$G$.\textit{\ }In particular near a point $p\in P$ we have the
local trivialization $\left( \pi ^{-1}\left( U\right) ,\varphi
_{\pi ^{-1}\left( U\right) }\right) $ where
\[
\varphi _{\pi ^{-1}\left( U\right) }:\mathcal{V}P\mid _{\pi
^{-1}\left( U\right) }\stackrel{\cong }{\longrightarrow }\pi
^{-1}\left( U\right) \times \mathfrak{g}
\]
is the diffeomorphism defined by setting
\[
\varphi _{\pi ^{-1}\left( U\right) }\left( u\right) :=\left( q,d_{q}\left(
r_{\sigma _{\pi \left( q\right) }^{U}\left( q\right) }^{-1}\circ \sigma
_{\pi \left( q\right) }^{U}\right) \left( u\right) \right)
\]
for every $q\in \pi ^{-1}\left( U\right) $ and $u\in \mathcal{V}_{q}P$.
\end{remark}

\noindent Recall the definition \ref{def conn} of a connection
\index{connection}on a principal $G$--bundle $\left( P,\pi \right) $. It is
not difficult to show that every principal bundle on a paracompact manifold $%
M$ admits a connection (see e.g. \cite{Poor 1981}, theorems 2.35 and 9.3). Given a
connection $\mathcal{H}P\subset TP$ we can uniquely split a vector field $%
X:P\longrightarrow TP$ into a horizontal part $\mathcal{H}X:P\longrightarrow
\mathcal{H}P$ and a vertical part $\mathcal{V}X:P\longrightarrow $ $\mathcal{%
V}P$ such that for every $p\in P$
\begin{equation}
X_{p}=\mathcal{H}_{p}X+\mathcal{V}_{p}X  \label{split}
\end{equation}
Recalling definition \ref{conn.forms} let $A\in \Omega ^{1}\left( P,\mathfrak{g}%
\right) $ be the $\mathfrak{g}$--valued 1--form associated with the connection $%
\mathcal{H}P$ and$\ \Omega \in \Omega ^{2}\left(
P,\mathfrak{g}\right) $ be its curvature $\mathfrak{g}$--valued
2--form. \noindent These forms are related to each other
by the \textit{structure equation}%
\index{equation, structure}
\[
\Omega \left( X,Y\right) =dA\left( X,Y\right) +\left[ AX,AY\right]
\]
for any vector fields $X,Y$ on $P$. We can rewrite it in the following
shorter form
\begin{equation}
\Omega =dA+%
\frac{1}{2}\left[ A,A\right]  \label{struct.eqn}
\end{equation}
by setting $\left[ A,A\right] \left( X,Y\right) :=\left[ AX,AY\right]
-\left[ AY,AX\right] $.

\noindent Let $l_{\sigma }$ be the automorphism of $G$ given by
left multiplication by $\sigma \in G$. The dual vector space
$\mathfrak{g}^{*}$ of the Lie algebra $\mathfrak{g}$ can be
canonically identified with the vector space of all \textit{left
invariant }$1$\textit{--forms on} $G$ since all such forms assume
constant values on left invariant vector fields. The composition
\[
a_{\sigma }:=l_{\sigma }\circ r_{\sigma ^{-1}}:G\longrightarrow G
\]
is an automorphism of $G$ fixing $id\in G$. Therefore its differential
\begin{equation}
Ad_{\sigma }:=d_{id}a_{\sigma }  \label{Ad}
\end{equation}
may be thought as an automorphism of $\mathfrak{g}\cong T_{id}G$ and%
\index{$Ad$--representation } its codifferential $\delta
_{id}a_{\sigma }$ as an automorphism of $\mathfrak{g}^{*}$.

\begin{proposition}
\textit{Let us consider }$\theta \in \mathfrak{g}^{*}$ \textit{and
}$X,Y\in \mathfrak{g}$\textit{. Then for every} $\sigma \in G$
\begin{equation}
\left( \delta r_{\sigma }\right) \theta X=\left( \theta \circ Ad_{\sigma
}\right) X  \label{adj-transf}
\end{equation}
\textit{and they satisfy the Maurer--Cartan equation}\footnote{%
For this reason left invariant 1--forms are also called \textit{%
Maurer--Cartan forms.}}%
\index{equation, Maurer--Cartan}%
\index{Maurer--Cartan, form}
\index{Maurer--Cartan, equation}
\begin{equation}
d\theta \left( X,Y\right) =-\theta \left[ X,Y\right]   \label{MCeqn}
\end{equation}
\end{proposition}

\noindent \textit{Proof.\quad }To prove $\left(
\text{\ref{adj-transf}}\right) $ note that for every $\tau \in G$ left
invariance of $\theta $ gives
\[
\theta _{\tau \sigma }=\left( \delta _{\tau \sigma }l_{\sigma ^{-1}}\right)
\theta _{\sigma ^{-1}\tau \sigma }
\]
which implies
\[
\left( \delta _{\tau }r_{\sigma }\right) \theta _{\tau \sigma }=\left(
\delta _{\tau }r_{\sigma }\circ \delta _{\tau \sigma }l_{\sigma
^{-1}}\right) \theta _{\sigma ^{-1}\tau \sigma }=\left( \delta _{\tau
}a_{\sigma ^{-1}}\right) \theta _{\sigma ^{-1}\tau \sigma }=\theta _{\sigma
^{-1}\tau \sigma }\circ d_{\tau }a_{\sigma }
\]
To restrict this relation to a left invariant vector field $X\in
\mathfrak{g}$ means to choose $\tau =id$ and so to obtain just
$\left( \text{\ref {adj-transf}}\right) $. For $\left(
\text{\ref{MCeqn}}\right) $ let us observe that almost by
definition
\[
d\theta \left( X,Y\right) =X\theta Y-Y\theta X-\theta \left[ X,Y\right]
\]
Since $X,Y\in \mathfrak{g}$ left invariance of $\theta $ implies that both $%
\theta Y$ and $\theta X$ are constant functions. This suffices to
finish the proof. $\diamondsuit$
\smallskip \medskip

Given a point $p\in P$ let us now consider the codifferential
\[
\delta \lambda _{p}:T^{*}P\longrightarrow T^{*}G
\]
and let $A$ be the connection form of $\mathcal{H}P$. We can then
define the $\mathfrak{g}$--valued 1--form $\left( \delta \lambda
\right) A\in \Omega ^{1}\left( G,\mathfrak{g}\right) $ by setting
\begin{equation}
\left( \left( \delta \lambda \right) A\right) _{\sigma }:=\left( \delta
_{\sigma }\lambda _{p}\right) A_{p\sigma }  \label{MCconnection}
\end{equation}
for every $\sigma \in G$. This definition is not dependent on the choice of $%
p\in P$ since by $\left( \text{\ref{g-connection}}\right) $ we have for
every $v\in T_{\sigma }G$%
\[
\left( \delta _{\sigma }\lambda _{p}\right) A_{p\sigma }\left( v\right)
=A_{p\sigma }\left( \left( d_{\sigma }\lambda _{p}\right) v\right) =\left(
d_{id}\lambda _{p\sigma }\right) ^{-1}\left( \mathcal{V}_{p\sigma }\left(
d_{\sigma }\lambda _{p}\right) v\right)
\]
Since $\lambda _{p}$ is a diffeomorphism of $G$ onto the fiber $\pi
^{-1}\left( \pi \left( p\right) \right) $ it follows that $\left( d_{\sigma
}\lambda _{p}\right) v\in \mathcal{V}_{p\sigma }P$ and

\begin{equation}
\left( \delta _{\sigma }\lambda _{p}\right) A_{p\sigma }\left( v\right)
=\left( d_{id}\lambda _{p\sigma }\right) ^{-1}\left( \left( d_{\sigma
}\lambda _{p}\right) v\right) =d_{\sigma }\left( \lambda _{p\sigma
}^{-1}\circ \lambda _{p}\right) v=\left( d_{id}l_{\sigma }\right) ^{-1}v
\label{l.i.ext}
\end{equation}
where the last equality follows by differentiating the commutative diagram
\[
\begin{array}{lll}
G & \stackrel{\lambda _{p}}{\longrightarrow } & P \\
_{l_{\sigma ^{-1}}}\searrow &  & \swarrow _{\lambda _{p\sigma }^{-1}} \\
& G &
\end{array}
\]
The $\mathfrak{g}$--valued 1--form $\left( \delta \lambda \right)
A$ is actually left invariant since
\[
\delta _{\sigma }l_{\tau }\left( \left( \delta \lambda \right) A\right)
_{\tau \sigma }=\left( \delta _{\sigma }l_{\tau }\circ \delta _{\tau \sigma
}\lambda _{p}\right) A_{p\tau \sigma }=A_{p\tau \sigma }\circ d_{\sigma
}\left( \lambda _{p}\circ l_{\tau }\right)
\]
and given $v\in T_{\sigma }G$ we get
\begin{eqnarray*}
\delta _{\sigma }l_{\tau }\left( \left( \delta \lambda \right) A\right)
_{\tau \sigma }v &=&\left( d_{id}\lambda _{p\tau \sigma }\right) ^{-1}\left(
\mathcal{V}_{p\tau \sigma }d_{\sigma }\left( \lambda _{p}\circ l_{\tau
}\right) v\right) = \\
&&d_{\sigma }\left( \lambda _{p\tau \sigma }^{-1}\circ \lambda _{p}\circ
l_{\tau }\right) v=\left( d_{id}l_{\sigma }\right) ^{-1}v=\left( \left(
\delta \lambda \right) A\right) _{\sigma }v
\end{eqnarray*}
Therefore $\left( \delta \lambda \right) A\in \mathfrak{g}^{*}\otimes \mathfrak{g}%
\cong Hom\left( \mathfrak{g},\mathfrak{g}\right) $: call it%
\index{Maurer--Cartan, form of a connection} \textit{the
Maurer--Cartan form associated with the connection
}$\mathcal{H}P$\textit{. }By $\left( \text{\ref{l.i.ext}}\right) $
it is the $\mathfrak{g}$--valued\textit{\ }1--form which assigns
to each tangent vector to $G$ its left invariant extension: hence
its representative in $Hom\left( \mathfrak{g},\mathfrak{g}\right)
$ is the identity $id_{\mathfrak{g}}$ and the Maurer--Cartan
equation $\left( \text{\ref {MCeqn}}\right) $ gives
\[
d\left( \delta \lambda \right) A\left( X,Y\right) =-\left( \delta \lambda
\right) A\left[ X,Y\right] =-\left[ X,Y\right] =-\left[ \left( \delta
\lambda \right) AX,\left( \delta \lambda \right) AY\right]
\]
Then we get
\[
d\left( \delta \lambda \right) A+\frac{1}{2}\left[ \left( \delta \lambda
\right) A,\left( \delta \lambda \right) A\right] =0
\]
By defining $\left( \delta \lambda \right) \Omega $ just like we did for $%
\left( \delta \lambda \right) A$ in $\left( \text{\ref{MCconnection}}\right)
$ the structure equation $\left( \text{\ref{struct.eqn}}\right) $ and the
last one allows us to conclude that
\begin{equation}
\left( \delta \lambda \right) \Omega =0  \label{MCcurvature}
\end{equation}
Since $\delta _{id}\lambda _{p}$ realizes the isomorphism $\mathcal{V}%
_{p}^{*}P\cong \mathfrak{g}^{*}$ this actually means that
\textit{the curvature
2--form }$\Omega $\textit{\ vanishes on the tangent space to the fiber of }$%
P $\textit{.} Hence the structure equation $\left( \text{\ref{struct.eqn}}%
\right) $ can be rewritten as follows:
\[
dA=\Omega -\frac{1}{2}\left[ A,A\right]
\]
to give a decomposition of $dA$ into horizontal and vertical parts.

Let us now come back to consider the connection form $A$ of
$\mathcal{H}P$. It can be defined as in $\left(
\text{\ref{g-connection}}\right) $ since the connection
$\mathcal{H}P$ determines a splitting in the tangent bundle $TP$.
But also the converse is true and the connection $\mathcal{H}P$
may be obtained by the $\mathfrak{g}$--valued 1--form $A$ just
like the vector sub--bundle $\ker A$.

\begin{proposition}
\label{connections}\textit{If }$A$\textit{\ is the connection form of a
connection }$\mathcal{H}P$\textit{\ then}%
\index{connection}%
\index{form, connection}
\begin{eqnarray}
\forall p\in P,\forall u\in \mathcal{V}_{p}P\qquad \left( d_{id}\lambda
_{p}\right) A_{p}u=u  \label{conditions} \\
\forall \sigma \in G\qquad \delta R\left( \sigma \right) A=Ad_{\sigma
^{-1}}\circ A  \nonumber
\end{eqnarray}

\textit{Conversely, given a }$\mathfrak{g}$\textit{--valued
}$1$\textit{--form }$A $\textit{\ on }$P$\textit{\ satisfying
these conditions the vector sub--bundle }$\ker A\subset
TP$\textit{\ gives a connection on }$P$\textit{\
whose connection form is }$A$.\textit{\ Hence the set }$\mathcal{A}$\textit{$%
_{P}$ of all connections on }$P$\textit{\ can be identified
with the affine subspace of }$\Omega ^{1}\left(
P,\mathfrak{g}\right) $\textit{\ defined by conditions }$\left(
\text{\ref{conditions}}\right) $\textit{.}

\noindent \textit{Furthermore the curvature form }$\Omega \in \Omega
^{2}\left( P,\mathfrak{g}\right) $\textit{\ of }$\mathcal{H}$\textit{$P$ is a $%
\mathfrak{g}$--valued }$2$\textit{--form such that}
\begin{eqnarray}
\forall p\in P,\forall u,v\in \mathcal{V}_{p}P\qquad \Omega _{p}(u,v)=0
\label{curv-conditions} \\
\forall \sigma \in G\qquad \delta R\left( \sigma \right) \Omega =Ad_{\sigma
^{-1}}\circ \Omega   \nonumber
\end{eqnarray}
\end{proposition}

\noindent \textit{Proof.\qquad }The first equality in $\left( \text{\ref
{conditions}}\right) $ follows immediately by the definition of the
connection form $A$. For the second one note that
\[
\delta _{p}R\left( \sigma \right) A_{p\sigma }\left( u\right) =A_{p\sigma
}\left( d_{p}R\left( \sigma \right) u\right) =\left( d_{id}\lambda _{p\sigma
}\right) ^{-1}\mathcal{V}_{p\sigma }\left( d_{p}R\left( \sigma \right)
u\right)
\]
The condition $\left( \text{\ref{G-invariance}}\right) $ for the connection $%
\mathcal{H}P$ implies that $\mathcal{V}_{p\sigma }\left( d_{p}R\left( \sigma
\right) u\right) =d_{p}R\left( \sigma \right) \left( \mathcal{V}_{p}u\right)
$. On the other hand $\mathcal{V}_{p}u=d_{id}\lambda _{p}\left(
A_{p}u\right) $ and we can write
\[
\delta _{p}R\left( \sigma \right) A_{p\sigma }\left( u\right) =\left(
d_{id}\lambda _{p\sigma }\right) ^{-1}\circ d_{p}R\left( \sigma \right)
\circ d_{id}\lambda _{p}\left( A_{p}u\right) =Ad_{\sigma ^{-1}}\circ A\left(
u\right)
\]
where the last equality follows by the commutative diagram
\[
\begin{array}{ccc}
\pi ^{-1}\left( \pi \left( p\right) \right) & \stackrel{R\left( \sigma
\right) }{\longrightarrow } & \pi ^{-1}\left( \pi \left( p\sigma \right)
\right) \\
_{\lambda _{p}}\uparrow &  & \downarrow _{\lambda _{p\sigma }^{-1}} \\
G & \stackrel{a_{\sigma ^{-1}}}{\longrightarrow } & G
\end{array}
\]
For the converse it suffices to observe that the first equality in $\left(
\text{\ref{conditions}}\right) $ gives the splitting condition $\left( \text{%
\ref{splitting}}\right) $ and the second one ensures the $G$--invariance $%
\left( \text{\ref{G-invariance}}\right) $ for $\ker A$. Hence it is a
connection on $P$ whose connection form is clearly $A$.

\noindent Finally the first equality in $\left( \text{\ref{curv-conditions}}%
\right) $ follows by $\left( \text{\ref{MCcurvature}}\right) $ and the
second one by applying the second equality in $\left( \text{\ref{conditions}}%
\right) $ to the definition $\left( \text{\ref{g-curvature}}\right) $ of $%
\Omega $.\hfill
\smallskip \medskip

Let us recall that a \textit{gauge transformation of }$P$
\index{gauge transformation}is an automorphism $\varphi $ of $P$ which
induces the identity map on the base manifold $M$. Then it leaves every
fibre fixed and it makes sense to define the associated map
\begin{equation}
\sigma _{\varphi }:P\longrightarrow G  \label{gauge element}
\end{equation}
such that $\varphi \left( p\right) =p\sigma _{\varphi }\left( p\right) $. By
applying the \textit{Leibniz rule }to the connection form $A$ we get that
\[
\left( \delta _{p}\varphi \right) A_{\varphi \left( p\right) }=\delta
_{p}R\left( \sigma _{\varphi }\left( p\right) \right) A_{p\sigma _{\varphi
}\left( p\right) }+\left( \delta _{p}\sigma _{\varphi }\right) \left( \delta
\lambda \right) A_{\sigma _{\varphi }\left( p\right) }
\]
where $\left( \delta \lambda \right) A$ is the Maurer-Cartan form of the
given connection. The second equation in $\left(
\text{\ref{conditions}}\right) $ allows us to conclude that \textit{under a
gauge transformation} $\varphi $\textit{\ the connection form }$A$\textit{\
behaves as follows:}
\begin{equation}
\left( \delta \varphi \right) A=Ad_{\sigma _{\varphi }^{-1}}\circ A+\left(
\delta \sigma _{\varphi }\right) \left( \delta \lambda \right) A
\label{gauge on connection}
\end{equation}
\textit{If }$\Omega $\textit{\ is the associated curvature then by }$\left(
\text{\ref{MCcurvature}}\right) $ \textit{and }$\left( \text{\ref
{curv-conditions}}\right) $ \textit{\ it transforms under }$\varphi $\textit{%
\ as follows:}
\begin{equation}
\left( \delta \varphi \right) \Omega =Ad_{\sigma _{\varphi }^{-1}}\circ
\Omega  \label{gauge on curvature}
\end{equation}
Since gauge transformations on $P$ form a group $\mathcal{G}_{P}$ with
respect to the composition, $\left( \text{\ref{gauge on connection}}\right) $
defines an action of $\mathcal{G}_{P}$ on the affine space of connections $%
\mathcal{A}_{P}$.

Let us now consider the \textit{exponential map }$\exp :\mathfrak{g}%
\longrightarrow G$ \index{exponential map}which assigns to a left
invariant vector field $X\in \mathfrak{g}$ the point $\exp
_{X}\left( 1\right) \in G$ where $\exp _{X}\left( t\right) $ is
the unique 1--parameter group whose tangent vector at $0\in
\R$ is $X_{id}\in T_{id}G$. Since $Ad_{\sigma }\in $ Aut$\left( \mathfrak{g}%
\right) $, for every $\sigma \in G$, and the Lie algebra of Aut$\left( \mathfrak{%
g}\right) $ is End$\left( \mathfrak{g}\right) $ we get the
following commutative diagram:
\[
\begin{array}{ccc}
G & \stackrel{Ad}{\longrightarrow } &
\text{Aut}\left( \mathfrak{g}\right) \\
_{\exp }\uparrow &  & \uparrow _{\exp } \\
\mathfrak{g} & \stackrel{ad}{\longrightarrow } & \text{End}\left(
\mathfrak{g}\right)
\end{array}
\]
where $ad:=d\left( Ad\right) $.

\begin{definition}
\label{killing}For every $X,Y\in \mathfrak{g}$ the symmetric
bilinear form
\[
\left\langle X,Y\right\rangle :=\text{tr}\left( ad_{X}\circ ad_{Y}\right)
\]
is called the \textit{Killing form}%
\index{form, Killing}\textit{\ of the lie algebra
}$\mathfrak{g}$\textit{.}
\end{definition}

Given a point $m\in M$ recall the definition $\left(
\text{\ref{Hol(m)}}\right) $ of the\textit{\ holonomy group }Hol$_{\mathcal{H%
}P}\left( m\right) $\textit{\ of a connection }$\mathcal{H}$\textit{$P$ at }$%
m\in M$\textit{.} If the base manifold $M$ is connected all these groups are
isomorphic when $m $ varies in $M$ since we can send
\begin{equation}
h_{\gamma }\in \text{Hol}_{\mathcal{H}P}\left( m_{1}\right) \longmapsto
h_{\alpha *\gamma *\overline{\alpha }}\in \text{Hol}_{\mathcal{H}P}\left(
m_{2}\right)  \label{connect-hol}
\end{equation}
where $\alpha $ is a path from $m_{1}$ to $m_{2}$ and $\overline{\alpha }$
its reversed path. Then it make sense to define \textit{the holonomy group }%
Hol$_{\mathcal{H}P}$\textit{\ of the connection }$\mathcal{H}$\textit{$P$.}

\section{\label{Witten} Appendix: More on Witten's open string theory
interpretation of QFT}

\noindent \textit{Sketch of proof of Theorem \ref{Witten92}:} We have to
show that under the assumptions $\left( \text{\ref{ctg.hyp.}}\right) $ and $%
\left( \text{\ref{bdr.cond.}}\right) $ the weak coupling limit of the
abstract string Lagrangian reduces exactly to the Lagrangian of a QFT on $L$.

\noindent The low energy (or weak coupling) limit \index{limit,
low energy (weak coupling)}of a string theory is only approximated
by a QFT since the limit Lagrangian admits perturbative
corrections depending on the coupling constant and non--constant
instanton corrections (see definition \ref {instanton}). The
string theory analyzed in \cite{Witten1992} is a topological
theory given by an $A$--twisted sigma model. At first Witten
observes that this model does not depend on the coupling constant of
the theory, implying that there cannot be any perturbative
correction in the limit Lagrangian.

\noindent It remains then to show that \textit{all the non--constant
instanton contributions vanish}. Let $\sigma $ be the canonical symplectic
form on $\widehat{Y}=T^{*}L$. It is the differential of the Liouville form,
i.e. in local canonical coordinates $\sigma =d\vartheta $ where $\vartheta
:=\sum_{j=1}^{3}p_{j}dq_{j}$. The Liouville form vanishes on $L$ given by $%
p_{1}=p_{2}=p_{3}=0$. Note that the bosonic sigma model action reduces for
instantons to be
\[
I=\int_{\Sigma }\phi ^{*}\left( \sigma \right)
\]
Stokes' theorem and condition $\left( \text{\ref{bdr.cond.}}\right) $
suffice to conclude that
\begin{equation}
I\left( \phi \right) =0  \label{vanishing}
\end{equation}
for \textit{all} instantons $\phi $. On the other hand by its definition the
bosonic sigma model action $I$ vanishes \textit{only} for constant
instantons. Hence we can admit only constant instanton corrections and the
abstract string Lagrangian reduces exactly to the Lagrangian of the QFT on $L
$ realizing the low energy limit. In the $A$--twisted case such a limit
turns out to be exactly a Chern--Simons $U\left( N\right) $--gauge theory.

\noindent \textit{Dropping assumption }$\left( \text{\ref{ctg.hyp.}}\right) $%
. The main result of \cite{Witten1992} is more general than Theorem \ref
{Witten92}. In fact he analyzes (section 4.4) the low energy limit of an $A$%
--twisted topological open string theory whose target space is given by a %
\cy  threefold $\widehat{Y}$ admitting $L$ as a Lagrangian submanifold.

\begin{theorem}
\label{Witten92bis}\textit{Let }$\widehat{Y}$\textit{\ be a local
Calabi--Yau threefold and }$L\subset \widehat{Y}$ \textit{a Lagrangian
submanifold. Then there exist topological string theories with }$\widehat{Y}$%
\textit{\ as target space, such that their open sectors are equivalent to a
QFT on }$L$ \textit{up to the convergence of non--constant instanton
contributions. In the }$A$\textit{--twisted case the Lagrangian action of
the limit QFT is (if convergent) a }deformation \textit{of a Chern--Simons
action}.
\end{theorem}

\noindent This result follows by assuming the same boundary conditions as
above. But now $\left( \text{\ref{bdr.cond.}}\right) $ is no longer sufficient
to conclude the vanishing $\left( \text{\ref{vanishing}}\right) $ for
non--constant instantons: given $\phi $, its\index{instanton, number} \textit{%
instanton number }is
\[
q\left( \phi \right) :=\int_{\Sigma }\phi ^{*}\left( \omega \right)
\]
where $\omega $ is the symplectic form of $\widehat{Y}$. Instanton numbers
turn out to be non--negative. For any knot $\mathcal{K}\subset \phi \left(
\partial \Sigma \right) \subset L$ consider the Wilson line $W_{\mathcal{K}%
}^{R}$ constructed by holonomy on $L$. For a given connection $A$ on a $%
U\left( N\right) $--principal bundle Witten shows that the instanton
contribution of $\phi $ is given by
\[
-\frac{i\eta \left( \phi \right) e^{-\theta q\left( \phi \right) }}{2\pi k}%
\sum_{\mathcal{K}\subset \phi \left( \partial \Sigma \right) }\log \left(
\text{ tr}_{R}\left( h_{\mathcal{K}}\right) \right)
\]
where $\theta $ is a positive real parameter, $e^{-\theta q\left( \phi
\right) }$ a suitable weighting factor and $\eta \left( \phi \right) =\pm 1$%
. If $S\left( \mathcal{L}\left( A\right) \right) $ is the Chern--Simons
action on $L$ the limit action turns out to be
\begin{equation}
S^{\prime }=S\left( \mathcal{L}\left( A\right) \right) -\frac{i}{2\pi k}%
\sum_{\phi }\left[ \eta \left( \phi \right) e^{-\theta q\left( \phi \right)
}\sum_{\mathcal{K}\subset \phi \left( \partial \Sigma \right) }\log \left(
\text{ tr}_{R}\left( h_{\mathcal{K}}\right) \right) \right]
\label{deformedCS}
\end{equation}
Under suitable assumptions on the ``moduli space'' of instantons $\phi $ the sum
can be perturbatively evaluated for $\theta \gg 0$.

\begin{corollary}
\label{unknot limit}\textit{Assume that }$\widehat{Y}=T^{*}S^{3}$\textit{\
and }$L=\mathcal{C}$\textit{\ is the Lagrangian submanifold given by the
conormal bundle of the unknot knot in }$S^{3}$\textit{\ like in Proposition
\ref{open istanton}. Then the low energy limit QFT on }$\mathcal{C}$\textit{%
\ of the open sector of a type II--}$A$\textit{\ string theory with }$M$%
\textit{\ D--branes wrapped around }$\mathcal{C}$\textit{\ is a }$SU\left(
M\right) $\textit{--Chern--Simons gauge theory on }$\mathcal{C}$. \textit{%
Moreover the global open string theory with }$N$\textit{\ D--branes wrapped
around }$S^{3}$\textit{\ and }$M$ \textit{D--branes wrapped around }$%
\mathcal{C}$\textit{\ admits a low energy limit QFT whose action is the
following deformation of the }$SU\left( M\right) $ \textit{Chern--Simons
action on }$\mathcal{C}$:
\[
S^{\prime }=S\left( \mathcal{L}\right) -\frac{i}{2\pi k}\sum_{d}\eta
_{d}\log \left( \text{ tr}_{R}\left( h_{\mathcal{K}}^{d}\right) \right)
\]
\end{corollary}

\noindent The first part of the statement can be proved like Theorem
\ref {Witten92} since the Liouville form of $\R^{8}\supset
T^{*}S^{3}$ vanishes when restricted to $\mathcal{C}$, as in
$\left( \text{\ref {lagrangian}}\right) $. That is enough to
guarantee the vanishing $\left( \text{\ref{vanishing}}\right) $.

\noindent To prove the second part, note that the only non--trivial
non--constant contributions come from instantons $\phi $ such that $\phi
\left( \partial \Sigma \right) $ is a $d$--covering of the unknot in $S^{3}$%
. For these instantons $q\left( \phi \right) =0$ by Stokes' theorem and the
statement follows by $\left( \text{\ref{deformedCS}}\right) $.


\begin{thebibliography}{10}

\bibitem{Acharya1999}
 Acharya BS  1999 {\em M theory, Joyce orbifolds and super
Yang-Mills}  Adv.Theor.Math.Phys. {\bf 3} 227--248;
 {\urlfont hep--th/9812205}.

 \bibitem{Acharya00}
 Acharya BS  2000 {\em On realizing N=1 super Yang-Mills in M theory}
 1--16 {\urlfont hep--th/0011089}


\bibitem{AKV}
Aganagic M, Klemm A, Vafa C 2001 ``Disc instantons, mirror
symmetry and the duality web" Z. Naturforsch. {\bf A57} 1-28
(2002) ; hep--th/0105045.

\bibitem{alt1997}
 Altmann, K. 1997, {\em The versal deformation of an isolated
toric singularity}   Invent. Math.{\bf 128}, no. 3, 443--479.

\bibitem{am}
  Aspinwall PS and Morrison DR 1993  {\em Topological field
theory and rational curves}  Commun. Math. Phys. {\bf 151}
245--262.

\bibitem{Atiyah1989}
 Atiyah MF 1989 ``Topological quantum field theories" {\it
Publ. Math. I.H.E.S.} {\bf 68} 175--186.

\bibitem{Atiyah1990a} Atiyah MF 1990 {\it The geometry and physics of knots} Cambridge University
Press.

\bibitem{Atiyah1990b} Atiyah MF 1990 ``On framing of 3--manifolds" {\it Topology } {\bf 29 } 1--8.

\bibitem{AB}
 Atiyah MF and Bott R 1982 ``The Yang--Mills equations over Riemann surfaces"
{\it Phil. Trans. R. Soc. Lond.} {\bf A 308} 523--615.

\bibitem{amv}
Atiyah MF, Maldacena J and Vafa C 2001 {\em An M--theory flop as a
large $N$ duality} { J.Math.Phys.} {\bf 42} 3209-3220.

\bibitem{aw}
Atiyah MF and Witten E 2001 ``$M$--theory dynamics on a
manifold of $G_2$ holonomy" hep--th/0107177.

\bibitem{Avram1996}
Avram A, Candelas P, Jan\v {c}i\'{c} D and Mandelberg  1996 ``On
the connectedness of the moduli space of \cy manifolds" Nucl.
Phys. {\bf B 465} 458--472; hep--th/9511230.

\bibitem{Axelrod1991}
Axelrod S,
Della Pietra S and Witten E 1991 ``Geometric quantization of
Chern--Simons
gauge theory" {\it J. Diff. Geom.} {\bf 33} 787--902.

\bibitem{Barth1984}
 Barth W,
Peters C and Van de Ven A 1984 {\it Compact complex surfaces}
Springer--Verlag.

\bibitem{Batyrev1998}
Batyrev V, Ciocan-Fontanine I, Kim B and van Straten D 1998
``Conifold transitions and mirror symmetry for Calabi-Yau complete
intersections in Grassmannians" Nucl. Phys. {\bf B 514} 640--666;
alg--geom/9710022.

\bibitem{Berglund1995}
 Berglund P, Katz S and Klemm A 1995 ``Mirror
symmetry and the moduli space for generic hypersurfaces in
toric varieties" Nucl. Phys. {\bf B456} 153--204.

\bibitem{Berglund1997}
 Berglund P, Katz S, Klemm A
and Mayr P 1997 ``New Higgs transition between dual $N=2$ string
models" Nucl. Phys. {\bf B483} 209--228.

\bibitem{bcov}
Bershadsky M, Cecotti S, Ooguri H and Vafa C 1993 ``Holomorphic
anomalies in topological field theories" (appendix by Katz S)
Nucl. Phys. {\bf B 405} 279--304  and in {\it Mirror symmetry II}
AMS/IP Stud. Adv. Math. {\bf 1} 655--682; hep--th/9302103.

\bibitem{bcov2}
Bershadsky M, Cecotti S, Ooguri H and Vafa C 1994 ``Kodaira--Spencer theory
of gravity and exact results for quantum string amplitudes" {\it Commun. Math. Phys.}
{\bf 165} 311--428; hep--th/9309140.

\bibitem{bggg}
 Brandhuber A, Gomis J, Gubser SS and  Gukov S 2001 ``Gauge theory
at large $N$ and new $G_2$ holonomy metrics" Nucl. Phys. {\bf B
611} 179--204; hep--th/0106034.

\bibitem{brie1966}
 Brieskorn E 1966 ``\"{U}ber die
Aufl\"{o}sung gewisser Singularit\"{a}ten
von holomorphen Abbildungen''
{\it Math. Ann.} {\bf 166} 76--102.

\bibitem{brie1968}
 Brieskorn E 1968 ``Die
Aufl\"{o}sung der rationalitaten Singularit\"{a}ten
holomorpher
Abbildungen'' {\it Math. Ann.} {\bf 178} 255--270.

\bibitem{Bryant1987}
 Bryant RL 1987 ``Metrics with exceptional holonomy" {\it Ann. of
Math.} {\bf 126} 525--576.

\bibitem{BS} Bryant RL and Salamon SM 1989 ``On the construction of some
complete metrics with exceptional holonomy" {\it Duke Math. J.} {\bf 58}
829--850.

\bibitem{CdO} Candelas P, de la Ossa XC, 1990 ``Comments on conifolds"Nucl. Phys. {\bf B 342}, n.1,   246--268.

\bibitem{CdOGP} Candelas P, de la Ossa XC, Green PS and Parkes L 1991 ``A pair of \cy
manifolds as an exactly soluble superconformal theory" Nucl. Phys.
{\bf B 359} 21--74 and in {\it Essays on mirror manifolds}
International Press 31--95.


\bibitem{CdOFKM} Candelas P, de
la Ossa XC, Font A, Katz S and Morrison DR 1994 ``Mirror symmetry
for two parameter models, I" Nucl. Phys. {\bf B416} 481--538.

\bibitem{CS} Chern S, Simons J
1974 ``Characteristic forms and geometric invariants"
{\it Ann. Math.}
{\bf 99}, 48--69.

\bibitem{Chiang1996}
Chiang T--M, Greene B, Gross M and Kanter Y 1996
``Black hole condensation and the web of \cy manifolds" {\it Nucl. Phys. Proc. Supp.}
{\bf 46} 82--95; hep--th/9511204.

\bibitem{Clemens1983} Clemens C. H. 1983 ``Double Solids'' {\it Adv. in
Math.} {\bf 47} 107--230.

\bibitem{CKM}
 H. Clemens, J. Koll\'ar, S. Mori  1988 ``Higher Dimensional
Geometry"  {\it   ``Asterisque"} {\bf 166}


\bibitem{CoxKatz1999}
 Cox AD and Katz S 1999 {\it Mirror
Symmetry and Algebraic Geometry}
vol.~68  Math. Surveys and Monographs
(American Mathematical Society).

\bibitem{Deligne1999}
Deligne P, Etingof P, Freed DS, Jeffrey LC,
Kazhdan D, Morgan JW, Morrison D and Witten E eds. 1999
{\it Quantum fields and strings: a course for mathematicians I, II}
Americam Mathematical Society, IAS.

\bibitem{WLP}
de Wit B, Lauwers P and van
Proeyen A 1985 ``Lagrangians of $N=2$ supergravity--matter
systems" Nucl. Phys. {\bf B255} 569--608.

\bibitem{d}
 Donaldson SK 1983 ``A new proof of a theorem
of Narasimhan and Seshadri"
{\it J. Diff. Geom.} {\bf 18}
269--277.

\bibitem{dkv}
 Douglas M, Katz S and Vafa C 1997 ``Small instantons, del Pezzo surfaces
and type $I^{\prime}$ theory"  Nucl. Phys. {\bf B 497}  155--172.

\bibitem{fp}
Faber C and Pandharipande R 2000  ``Hodge integrals and Gromov--Witten
theory" {\it Invent. Math.} {\bf 139} 173--199; math.AG/9810173.

\bibitem{df}
 Freed DS 1995 ``Classical Chern--Simons theory, Part 1"
{\it Adv. Math.}
{\bf 113} 237--303.

\bibitem{fy}
Freyd P,  Yetter D, Hoste J,
Lickorish WBR, Millett K and Ocneanu A 1985
``A new polynomial invariant of
knots and links" {\it Bull. AMS} {\bf 12} 183--312.

\bibitem{Friedman1986}
Friedman R 1986
``Simultaneous resolution of threefold double points''
{\it Math. Ann.}
{\bf 274} 671--689.

\bibitem{fup}
 Fulton W and Pandharipande R 1995  ``Notes on stable maps and quantum
cohomology" in {\it Algebraic Geometry--Santa Cruz 1995}  Proc. Sym. Pure Math.
{\bf 62} Part 2, AMS,
45--96.

\bibitem{gp}
Getzler E and Pandharipande R 1998  ``Virasoro constraints and the Chern
classes of the Hodge bundle"  math.AG/9805114.

\bibitem{Gibbons1990}
Gibbons GW, Page DN and Pope CN 1990 ``Einstein metrics on $S^3,
\R ^3$ and $\R ^4$ bundles" {\it Commun. Math. Phys.} {\bf 127} 529--553.

\bibitem{gva}
Gopakumar R and Vafa C 1998 ``$M$--theory and topological
strings--I" hep--th/9809187

\bibitem{gvb}
 Gopakumar R and Vafa C 1998 ``$M$--theory and topological strings--II"
hep--th/9812127

\bibitem{gvc}
 Gopakumar R and Vafa C 1999 ``On the gauge
theory/geometry
correspondence" {\it Adv.Theor.Math.Phys.} {\bf 3}
1415--1443.

\bibitem{Graber--Pandharipande 1999} Graber T and Pandharipande R 1999
``Localization of virtual classes" {\it Invent. Math.} {\bf 135} 487-518.

\bibitem{Greene1995}
Greene B, Morrison DR and Strominger A 1995 ``Black hole
condensation and the unification of string vacua" Nucl. Phys. {\bf
B451} 109--120.

\bibitem{GP}
Greene B and Plesser MR 1990
``Duality in \cy moduli space" Nucl. Phys. {\bf B338} 15--37.

\bibitem{pg}
 Griffiths P
1969 ``On the periods of certain rational integrals I, II" {\it Ann. of
Math. } {\bf 90} 460--495, 498--541.

\bibitem{Gross1997(a)}
 Gross M 1997
 ``Deforming Calabi--Yau threefolds" {\it Math. Ann. 308} {\bf 2}, 187--220.

\bibitem{Gross1997(b)}
 Gross M 1997  ``Primitive Calabi-Yau threefolds" {\it J. Diff. Geom.}
{\bf 45} 288--318.

\bibitem{Helgason1978}
 Helgason S 1978 {\it
Differential geometry, Lie groups and symmetric spaces}
Academic Press,
New York.

\bibitem{Hitchin 1990} Hitchin NJ 1990 ``Flat connections and geometric
quantization" {\it Commun.
Math. Phys.} {\bf 131} 347--380.

\bibitem{Horowitz--Strominger 1991}
Horowitz G, Strominger A 1991 ``Black Strings and $p$--branes"
Nucl. Phys. {\bf B360} 197--209.

\bibitem{hkty}
 Hosono S, Klemm A, Theisen S and Yau ST 1995 ``Mirror
symmetry, mirror map and applications to complete intersection \cy
spaces" Nucl. Phys. {\bf B433} 501--554 and in {\it Mirror
symmetry II} AMS/IP Stud. Adv. Math. {\bf 1} 545--606;
hep--th/9406055.

\bibitem{Jones 1985} Jones V 1985 ``A polynomial invariant for knots via
von Neumann algebras"
{\it Bull. AMS} {\bf 12} 103--111.

\bibitem{Jones 1987} Jones V
1987 ``Hecke algebras representations of braid groups and link
plynomials"
{\it Ann. Math} {\bf 126} 335--388.

\bibitem{Johnson 1998} Johnson CV 1998   ``Etudes on D-Branes" hep--th/9812196.

\bibitem{Johnson 2000} Johnson CV 2000   ``D--Brane primer" in {\it Strings, branes and
gravity} TASI 99 World Sci. Publ., 129-350; hep--th/0007170.

\bibitem{Joyce 1996} Joyce, D 1996 ``Compact Riemann 7--manifolds with holonomy $G_2$
I, II" {\it J. Diff. Geom.} {\bf 43} 291--328, 329--375.

\bibitem{Joyce 2000} Joyce, D 2000 {\it
Compact manifolds with Special Holonomy} Oxford Science
Publications, Oxford--New York.

\bibitem{kl}
 Katz S and Liu CCM 2001  ``Enumerative geometry of stable maps with Lagrangian
boundary conditions and multiple covers of the disk", Adv. Theor. Math. Phys. {\bf 5}
1--49; math.AG/0103074.

\bibitem{kmp}
Katz S, Morrison DR and Plesser
M 1996 ``Enhanced gauge symmetry in type II string theory"  Nucl.
Phys. {\bf B477} 105--140.


\bibitem{MK} Kontsevich M, 1995 ``Enumeration of Rational Curves via Torus
Actions'', {\it The Moduli Space of Curves}, 335-368, Progr. Math.
{\bf 129}, Birkh\"auser Boston, MA

\bibitem{Labastida1999} Labastida JMF 1999  ``Chern--Simons gauge
theory: ten years after"
 e--print: hep--th/9905057

\bibitem{lmv}
Labastida JMF, Mari\~{n}o M and Vafa C 2000
``Knots, links and branes at large $N$" {\it JHEP} 0011:007;
 hep--th/0010102

\bibitem{lm}
 Labastida JMF and Mari\~{n}o M 2001  ``Polynomial invariants for torus knots and
topological strings" {\it Commun. Math. Phys.} {\bf 217} 423--449; hep--th/0004196.

\bibitem{ls}
 Li J and Song YS 2001  ``Open string instantons and relative stable morphisms"
Adv. Theor. Math. Phys. {\bf 5} 67-91; hep--th/0103100.

\bibitem{Manin1995}
 Manin Y 1995 ``Generating functions in algebraic geometry and sums over trees2" in
{\it The moduli space of curves} Birk\"{a}user 401--417.

\bibitem{Manzoni 1842} Manzoni, A 1842 ``I Promessi Sposi"
{\it The Bethroted}

\bibitem{mv2001} Mari\~{n}o M
and Vafa C 2001 ``Framed knots at large $N$"  hep-th/0108064

\bibitem{Moore-Sieberg 1988} Moore G and Sieberg N 1988  ``Polynomial equations for rational
conformal field
theories" {\it Phys. Lett.} {\bf B212} 451--466.

\bibitem{Morrison 1985} Morrison DR
1985  ``The birational geometry of surfaces with rational double
points''
{\it Math. Ann.} {\bf 271} 415--438.

\bibitem{Morrison 1999} Morrison DR 1999
``Through
the looking glass'' in {\it Mirror Symmetry III}
( American
Mathematical
Society and International Press)  263-277.

\bibitem{Morrison--Seiberg 1997} Morrison DR, Seiberg N. 1997 ``Extremal transitions and five-dimensional
 supersymmetric field theories" Nuclear Phys. B {\bf 483}  229--247,
    [hep-th/9609070].


\bibitem{Narasimhan--Seshadri 1965} Narasimhan MS and
Seshadri CS 1965 ``Stable and unitary vector bundles
on a compact Riemann
surface" {\it Ann. Math.} {\bf 82} 540--567.

\bibitem{Ok}
Okonek C, Schneider M and Spindler H 1980 {\it Vector bundles on complex
projective spaces} Progr. in Math. {\bf 3}, Birk\"{a}user.

\bibitem{ov}
 Ooguri H and Vafa C 2000 ``Knot invariants and topological strings" Nucl. Phys.
{\bf B 577} 419--438.

\bibitem{Periwal 1993} Periwal V 1993 ``Topological closed--string interpretation of Chern--Simons
theory" {\it Phys. Rev. Lett.} {\bf 71} 1295--1298.

\bibitem{Pinkham 1983} Pinkham H 1983
``Factorization of birational maps in dimension 3''  in
{\it
Singularities} Proc. Sym. Pure Math. {\bf 40}, AMS,
343--372.

\bibitem{Poor 1981} Poor
W 1981 {\it Differential geometric structures} McGraw--Hill, New
York

\bibitem{Quillen 1986} Quillen DG 1986 ``Determinants of Cauchy--Riemann operators
over a
Riemann surface" {\it Funct. Anal. Appl.} {\bf 19} 31.

\bibitem{Quinn 1995}
Quinn F 1995 ``Lectures on axiomatic topological quantum field theory''  in
{\it Geometry and Quantum Field Theory} IAS/Park City Math. Series {\bf 1}
AMS,
325--459.

\bibitem{Ramadevi-Sarkar 2001} Ramadevi P and Sarkar T 2001  ``On link invariants and topological string
amplitude" Nucl. Phys. {\bf B600} 487--511;  hep--th/0009188.

\bibitem{Reid1980}
 Reid M 1980
``Canonical 3--Folds'' in
{\it
Journ\'{e}es de g\'{e}om\'{e}trie
alg\'{e}brique d'Angers}
(Sijthoff \&
Noorddhoff) 273--310.

\bibitem{Reid1983}
 Reid M 1983  ``Minimal models of canonical
3--Folds'' in
{\it Algebraic
varieties and analytic varieties} vol.~1 Adv.
Stud. Pure Math.
(North
Holland) 131--180.

\bibitem{Reid1987(a)}
Reid M 1987
``Young person's guide to canonical singularities''  in
{\it Algebraic
Geometry, Bowdoin 1985} vol.~1 Proc. Sym. Pure Math. {\bf 46},
AMS,
354--414.

\bibitem{Reid1987}
Reid M 1987  ``The moduli space of 3--folds with $K=0$ may
neverthless be irreducible'' {\it Math. Ann.} {\bf 278} 329--334.

\bibitem{Reidemeister1933} Reidemeister K 1933  ``Z\"{u}r dreidimensionalen
topologie" {\it Abh. Math. Sem.
Univ. Hamburg} {\bf 9} 189--194.

\bibitem{Schlessinger1971} Schlessinger, M 1971
``Rigidity of quotient singularities"
{\it Invent. Math.} {\bf 14} 17--26.

\bibitem{Segal1988}
Segal G 1988 ``Two--dimensional conformal field theories and modular
functions" {\it Proc. Int. Congr. Math. Phys. } Swansea 22--37.

\bibitem{SeibergWitten}
Seiberg N and Witten E 1994 ``Electromagnetic duality, monopole
condensation and confinement in $N=2$ supersymmetric Yang--Mills
theory" Nucl. Phys. {\bf B426} 19--52.

\bibitem{Sen1997} Sen A 1997 ``An introduction to non--perturbative string theory"
in \textit{Duality and supersymmetric theories}, 297--413 Cambridge;
hep--th/9802051.

\bibitem{Silverman} Silverman JH 1994  {\it Advanced topics in the arithmetic of
elliptic curves} vol.~151
Graduate Texts in Mathematics (Springer--Verlag).

\bibitem{Slodowy1980} Slodowy  1980 {\it Simple Singularities and Simple Algebraic Groups}
 vol.~815
Lecture Notes in Mathematics (Springer--Verlag).

\bibitem{Strominger1990}  Strominger A
1990  ``Special geometry" {\it Commun. Math. Phys. }
{\bf 133}
163--180.

\bibitem{Strominger1995}  Strominger A 1995  ``Massless black holes and
conifolds in string theory" Nucl. Phys. {\bf B451} 97--109.

\bibitem{'t Hooft 1974} 't Hooft G 1974 ``A planar diagram theory for strong
interactions" Nucl. Phys. {\bf B72} 461--473.

\bibitem{CT}  Taubes C, 2001 ``Lagrangians for the Gopakumar--Vafa
conjecture", Adv. Theor. Math. Phys{\bf 5} 139-163; math.dg/0201219.

\bibitem{Townsend1995} Townsend P 1995 ``The eleven--dimensional supermembrane revisited"
Nucl. Phys. {\bf B 350} 184--187; hep-th/9501068.

\bibitem{Tyurina1970}
Tyurina G 1970
``Resolution of the singularities of flat deformations of rational
double
points'' {\it Funkts. Analiz i ego Prilozh.} {\bf 4} 77--83.

\bibitem{Verlinde1988}
Verlinde E 1988  ``Fusion rules and modular transformations in
$2d$ conformal field theory" Nucl. Phys. {\bf B300} 360--383.

\bibitem{Vafa2001} Vafa C 2001 ``Superstrings and topological strings at large $N$" \ JMP {\bf 42}
2798--2817.

\bibitem{Voisin} Voisin C 1996  ``A mathematical proof of a formula of Aspinwall and Morrison"
{\it Comp. Math.} {\bf 104} 135--151.

\bibitem{Warner} Warner F 1983
{\it Foundations of differentiable manifolds and Lie groups}
 vol.~94
Graduate Texts in Mathematics (Springer--Verlag).

\bibitem{WvG}
 Werner J and van
Geemen B  ``New examples of threefolds with $c_1=0$''
{\it Math. Z.} {\bf
203} 211--225.

\bibitem{Wilson1989} Wilson PMH 1989 ``\cy manifolds with large
Picard
number''
{\it Invent.Math.} {\bf 98} 139--155.

\bibitem{Wilson1992} Wilson PMH
1992
``The \ka cone on \cy threefolds''
{\it Invent.Math.} {\bf 107}
561--583.

\bibitem{Witten1989} Witten E 1989  ``Quantum field theory and the Jones
polynomial"
{\it Commun. Math. Phys. }  {\bf 121} 351--399.

\bibitem{Witten1992} Witten E
1992 ``Chern--Simons gauge theory as a string theory"  in
{\it The Floer memorial volume} Birkh\"{a}user 1995, 637--678;
hep-th/9207094.

\bibitem{Witten1993} Witten E 1993 ``Phase of ${\mathcal {N}}=2$ theories in two dimensions" Nucl. Phys. {\bf B 403}
159--222;  hep-th/9301042.

\bibitem{Witten1995} Witten E 1995 ``String theory dynamics in various dimensions" Nucl. Phys.
{\bf B 443} 85--126; hep-th/9503124.

\end{thebibliography}
\end{document}